\documentclass[fleqn]{article}
\newtheorem{theorem}{Theorem}[section]

\newtheorem{lemma}[theorem]{Lemma}

\newtheorem{proposition}[theorem]{Proposition}
\newtheorem{definition}[theorem]{Definition}

\newtheorem{corollary}[theorem]{Corollary}

\title{Stopping time convergence for processes associated with Dirichlet forms}
\author{J.R. Baxter and  M. Nielsen Hernandez}

\usepackage{amsmath}
\usepackage{amsfonts}
\usepackage{amssymb}

\usepackage[titletoc]{appendix}

\usepackage{hyperref}  

\numberwithin{equation}{section}

\begin{document}

\maketitle

\begin{abstract}
Convergence is proved for solutions $u_n$ of Dirichlet problems
in regions with many small excluded sets (holes), as the holes become smaller and more
numerous.  The problem is formulated in the context of Markov processes associated with
 general Dirichlet forms, for random and nonrandom excluded sets.
Sufficient conditions are given in Theorem~\ref{random_sets_example_theorem}
under which the sequence of entrance
times or hitting times of the excluded sets converges in the stable topology.  Convergence in the 
stable topology is a  strengthened 
form of convergence in distribution,  introduced by R\'{e}nyi.  Stable convergence of 
the entrance times implies convergence of 
the solutions $u_n$ of the corresponding Dirichlet problems.
Theorem~\ref{random_sets_example_theorem} applies to Dirichlet forms 
such that the Markov  process associated 
with the form has continuous paths and satisfies
an absolute continuity condition for occupation time measures
(equation~\eqref{kernel_density_cond_eqn}). 
Conditions for convergence are formulated
in terms of the sum of the expectations of the
 equilibrium measures for the excluded sets.
The proof of convergence uses the fact that
any martingale with respect to the 
natural filtration of the process must be continuous.
In the case that the excluded sets are iid random,
Theorem~\ref{random_sets_example_theorem}
strengthens previous results in the classical Brownian motion setting.

Mathematics classification numbers 60J45, 60K37, 35J25
\end{abstract}

 \section{Introduction}
 \label{intro_section}
Let $X$ be a Markov process
and let~$U$ be an open subset of its state space $E$.
The probabilistic solution~$v$ of the Dirichlet problem 
on $U$, with killing rate $\alpha \in [0, \infty)$, source term~$f$ and boundary value function~$\varphi$,
is given by
\begin{equation}
\label{general_prob_Dirichlet_no_holds_eqn}
v(x) = \mathbf{E}_{x} \left[  \int_0^{\sigma} e^{ - \alpha t} f\left( X_t\right) \, d t \right]
+ \mathbf{E}_{x} \left[  e^{ - \alpha \sigma} \varphi \left( X_{\sigma} \right) \right].
\end{equation}
Here $\sigma$ is the exit time of $U$, 
$f$ is a measurable function on $U$,
and $\varphi$ is a measurable function on $E -  U$.
The solution $v(x)$ is given for those $x \in U$ such that 
$\mathbf{E}_{x} \left[  \int_0^{\sigma} e^{ - \alpha t} f \left( X_t\right) \, d t \right]$
and $\mathbf{E}_{x} \left[  e^{ - \alpha \sigma} \varphi \left( X_{\sigma} \right) \right]$ exist,
and $\varphi$ is defined to be zero at the cemetery point $\partial$
of the process.
In some applications it is natural to
consider the limit of a sequence of solutions~$u_n$
for which a subset $\Lambda\left( n\right)$ of the region $U$ is
excluded for each~$n$, so $u_n$ solves a Dirichlet problem 
on $U - \Lambda\left( n\right)$.  Typically $\Lambda\left( n\right)$ is
the union of many small sets $\Lambda_{j}\left( n\right)$, which become smaller and more numerous as $n \rightarrow \infty$,
and in this case $U - \Lambda\left( n\right)$ is often referred  to as a region with many small holes.
Define $u_n(x)$ on $U$ by
\begin{equation}
\label{general_prob_Dirichlet_eqn}
u_n(x) = \mathbf{E}_{x} \left[  \int_0^{\tau_n \wedge \sigma} e^{ - \alpha t} f\left( X_t\right) \, d t \right]
+ \mathbf{E}_{x} \left[  e^{ - \alpha \tau_n \wedge \sigma} \varphi \left( X_{\tau_n \wedge \sigma} \right) \right],
\end{equation}
where $\tau_n$ is the entrance time of $\Lambda\left( n\right)$ and~$x$ is such that 
the expected values exist.
We study conditions under which the solutions $u_n$ converge
to a limit $u$, which could then be considered as an approximation
to $u_n$ when the holes are small.
Analytical formulations of~\eqref{general_prob_Dirichlet_eqn}
can of course be given, as in Lemma~\ref{analytical_from_limit}) below,
and there are many approaches to this problem.
In the classical case of Brownian motion and $\alpha = 0$, where $- \Delta u_n  = f$ holds in $U - \Lambda\left( n\right)$, 
the size of each small set $\Lambda_{j}\left( n\right)$ as a target for Brownian motion is measured by its 
capacity.  When convergence holds in this setting, 
the limit~$u$ often satisfies the equation $- \Delta u + q u = f$ in $U$, where $q$
is an appropriate limiting density for the capacities of the holes $\Lambda_{j}\left( n\right)$.
Convergence problems in regions with many small holes have been considered 
for many classes of equations, both linear and nonlinear, using a variety
of techniques.  
Early results in this area include
\cite{khruslov-1972}, \cite{kac}, 
\cite{rauch_taylor}, 
and~\cite{papanicolaou_varadhan}. 
In the setting of  Dirichlet forms, convergence properties for sequences of
solutions of Dirichlet problems
have been studied using variational methods 
(cf. \cite{biroli_mosco}, \cite{dalmaso_decicco_notarantonio_tchou}, \cite{biroli_tchou}, 
\cite{mataloni_tchou}),  extending a similar 
approach for elliptic equations
(cf. \cite{dalmaso_mosco}, \cite{dalmaso_capacities}, \cite{dalmaso_garroni}).
The present paper uses a probabilistic approach based on \cite{baxter_chacon_jain},
applied to the Markov process~$X$ which is associated with a Dirichlet form.
Theorem~\ref{random_sets_example_theorem} gives conditions 
on a sequence of random sets $\Lambda\left( n\right)$ under which 
the stopping times $\tau_n$ in~\eqref{general_prob_Dirichlet_eqn}
converge with respect to the stable topology introduced by R\'{e}nyi \cite{renyi}.
This implies convergence for the corresponding solutions of the Dirichlet problem
(Lemmas~\ref{convergence_for_probability_solutions_lemma} and~\ref{analytical_from_limit}).
Theorem~\ref{random_sets_example_theorem} holds 
for a wide class of Dirichlet forms, and also strengthens earlier results in
the Brownian motion case.
Precise statements are given in Section~\ref{main_results_section}.
For definitions and results concerning Dirichlet forms we will refer to~\cite{ma_roeckner}, but
some of the facts which are needed are discussed
in Section~\ref{Dirichlet_form_properties_section}. 
Properties of stable convergence are given in
Section~\ref{stable_convergence_facts_section}.

\section{Main results}
\label{main_results_section}
The main theorem is given below in Section~\ref{random_holes_subsection}, after
some preliminary definitions.
\subsection{Stable convergence definitions}
Stable convergence was defined by R\'{e}nyi in~\cite{renyi},
as a 
stronger form of convergence in distribution.
A general treatment of stable convergence
of random variables is given in \cite{hausler_luschgy}.
Here we consider stable convergence of randomized stopping times,
as in \cite{baxter_chacon_compactness}, \cite{edgar_millet_sucheston}, \cite{meyer}).
By definition, a
randomized stopping time~$\tau$
with respect to a filtration $\left( \mathcal{F}_{t}\right)_{t \geq 0}$
on a probability space $\left( \Omega, \mathcal{F}, \mathbf{P}\right)$,
 is simply a stopping time on 
$\left( \Omega \times (0,1), \mathcal{F} \times \mathcal{B}^{\ast}, \mathbf{P} \times \lambda^{\ast}\right)$,
where $\lambda^{\ast}$ is Lebesgue measure on the Borel sets $\mathcal{B}^{\ast}$ of $(0,1)$,
and we use the enriched 
filtration $\left( \mathcal{F}_{t} \times \mathcal{B}^{\ast}\right)_{ t \geq 0}$
on the randomized space $\left( \Omega \times (0,1), \mathcal{F} \times \mathcal{B}^{\ast}, \mathbf{P} \times \lambda^{\ast}\right)$.
An ordinary stopping time on $\left( \Omega, \mathcal{F}, \mathbf{P}\right)$ 
can be regarded as defined on $ \Omega \times (0,1)$, 
and so an ordinary stopping time is a special case of a randomized one.
In this setting, stable convergence for randomized stopping times
is defined as follows.
Let~$\mathcal{G}$ be a sub-$\sigma$-field of~$\mathcal{F}$.
Let $\tau_n, \tau$ be randomized stopping times, 
or more generally let $\tau_n, \tau$ be any $\mathcal{F} \times \mathcal{B}^{\ast}$-measurable
maps from $ \Omega \times (0,1)$ to $[0, \infty]$.
Then $\tau_n \rightarrow \tau$,  $\mathbf{P}, \mathcal{G}$-stably,
if for every set~$G \in \mathcal{G}$ with $\mathbf{P}\left(G\right) > 0$, 
\begin{equation}
\label{stable_dist_conv_sets_G_eqn}
\tau_n \, \big\vert_{\displaystyle G \times (0,1)}\Longrightarrow \tau \, \big\vert_{\displaystyle G \times (0,1)} \ \text{in distribution}
\end{equation}
with respect to the conditional probability measure
$\mathbf{P} \times \lambda^{\ast}\left(\,\,  \cdot\, \,  \vert G \times (0,1)\right)$. 
Here $\tau_n \, \big\vert_{\displaystyle G \times (0,1)}$ denotes the restriction of $\tau_n$
to the set $G \times (0,1)$. 
When~$\mathcal{G}$ is known from the context 
we may simply write $\tau_n \rightarrow \tau$, $\mathbf{P}$-stably.
The space of randomized stopping times associated with a given filtration is 
closed and compact 
with respect to stable convergence (\cite{baxter_chacon_compactness}, \cite{meyer}).
When stopping times for a Markov process $\left(X_t\right)$ are considered, 
 unless otherwise stated $\mathcal{G}$ will be 
$\sigma\left(\mathcal{F}_{0}, X_t, \, 0 \leq t < \infty\right)$.
A sequence $\tau_n$ may converge stably with respect 
to one probability 
measure~$\mathbf{P}$ and not converge with respect to another 
probability measure~$\mathbf{Q}$, although convergence is preserved 
if $\mathbf{Q} << \mathbf{P}$ for enough sets,
which can be useful when Girsanov's theorem is applicable.
Other properties of stable convergence are stated
in Section~\ref{stable_convergence_facts_section}.
If, as in the present paper, 
$\mathcal{G}$ is
countably generated modulo sets of $\mathbf{P}$-measure zero,
the topology for stable convergence with respect to $\mathbf{P}$ is metrizable.

For any randomized stopping time~$\tau$
with respect to the  filtration $\left(\mathcal{F}_{t} \times \mathcal{B}^{\ast}\right)_{t \geq 0}$
on $ \Omega \times  (0,1)$, let
$F^{ \tau}_{t}  = 
\mathbf{E} \left[\left.\mathbf{1}_{ \left\{\tau \leq t\right\} } \,\,\right\vert \mathcal{F} \times \left\{\emptyset, (0,1)\right\} \right]$
for $t \in [0, \infty]$, where here the conditional expectation is with respect to 
the probability measure $\mathbf{P} \times \lambda^{\ast}$.
A version of $F^{ \tau} $ will always be chosen 
such that $t \mapsto F^{ \tau}_{t} (\omega)$
nondecreasing and right continuous for each~$\omega$ and such that 
$F^{ \tau}_{t} $ is $\mathcal{F}_t$-measurable. 
The randomized stopping time~$\tau$ can be chosen
so that $\tau(\omega, \cdot)$ is left-continuous and 
nondecreasing on $(0,1)$, and we will always use such a version.
For convenience in stating formulas, we also define $S^{\tau}_{t} = 1 - F^{ \tau}_{t} $. 
$S^{\tau}$ and $F^{ \tau} $ 
describe the observable properties of $\tau$.
For $\omega \in  \Omega$, 
 $S^{\tau}_{t}(\omega)$ can be thought of as the fraction of~$\omega$
which is not yet stopped at time~$t$.
In the case of an ordinary stopping time, $S^{\tau}_{t}$ is either zero or one for each~$t$.
It is easy to check that for any randomized stopping time $\tau$, 
\begin{equation}
\label{randomized_st_from_pathwise_dist_fn_eqn}
\tau(\omega, r) 
= \inf \left\{ u: \ u \in [0, \infty], \ S^{\tau}_{u} (\omega) \leq 1 - r\right\}.
\end{equation}
Also, $\tau$ is an 
$\left(\mathcal{F}_{t} \times \mathcal{B}^{\ast}\right)$-stopping time if and only if 
$\tau( \cdot, u)$ is a $\left( \mathcal{F}_{t}\right)$-stopping time
for each $u \in (0,1)$.

\subsection{Rate measures for randomized stopping times}
We are interested in stopping times in the setting of a Markov process.
Let $E$ be a separable metric space with Borel $\sigma$-algebra $\mathcal{B}$.
Let $\left( \Omega, \mathcal{F}, \left(\mathcal{F}_{t}\right)_{t \geq 0}, \left(X_t\right)_{ t \in [0, \infty]},
\mathbf{P}_{z}\right)$, $z \in E_{\partial}$,
 be  a  Hunt process with state space~$E$,
 cemetery point~$\partial$ and lifetime~$\zeta$.
Unless otherwise stated  $\mathcal{F}_{t}$ is the natural filtration for~$X$, that is, 
an appropriate closure of the filtration generated by~$X$,
and $\mathcal{F} = \mathcal{F}_{\infty}$.
Let $\Lambda\left( n\right)$ be a sequence of closed subsets of $E$, 
and let $\tau_n$ be the entrance time $D_{\Lambda\left( n\right)}$ of $\Lambda\left( n\right)$ or the hitting time 
$ T_{\Lambda\left( n\right)}$ of $\Lambda\left( n\right)$.  Sufficient
conditions will be given under which stable convergence holds for the stopping time sequence
 $\tau_n$.  In the cases studied here, convergence will be proved for situations in which the
sets $\Lambda\left( n\right)$ are sparse enough that~$\tau_n$ can
converge to a randomized stopping time~$\tau$ which
is associated with a \emph{rate} of stopping. The rate of stopping
is expressed by $S^{\tau}_{t} = e^{ - A_t}$ for $t \in [0, \infty)$, 
where  $A$ is a 
positive continuous
additive functional.
For a path~$\omega$ and $t \in [0, \infty)$, the probability that the path has not yet stopped
by time~$t$ is equal to $e^{ - A_t(\omega)}$. 
If there is a nonnegative Borel function~$h$
on the state space $E$ such that $A_t = \int_0^t h\left(X_s\right) \, d s$, 
then one can say that the stopping time~$\tau$ results from stopping at a rate $h(x)$
when the process~$X$ is near the point~$x$. 
 More generally, 
let~$m$ be a fixed $\sigma$-finite excessive measure on~$E$. 
Any positive continuous additive functional $A$ for~$X$ 
is associated  with a Revuz measure $\mu_A$ on $E$ with respect to~$m$
(cf. Theorem~A.3.5 in \cite{chen_fu}).  We will refer to the Revuz
measure $\mu_A$ as the rate measure for the stopping time.

\subsection{Assumptions on $X$}
\label{assumptions_X_subsection}
Unless otherwise stated, from now on it is assumed
 that the Hunt process~$X$ is properly associated with
a quasi-regular Dirichlet form  $\left( \mathcal{E}, D\left(\mathcal{E}\right)\right)$
on $L^{2}\left( m\right)$, not necessarily symmetric, as defined in  IV.1.13 and
 IV.2.5 of~\cite{ma_roeckner}. Here~$m$
 is a $\sigma$-finite measure on~$E$.  
$E$ is assumed to be a metrizable Lusin space, 
which is defined to be the continuous one-to-one image of a Polish space,
or, equivalently, a space which is
homeomorphic to a Borel subset of a compact metric space.
Unless otherwise stated, $\mathcal{E}$ is assumed to have
the local property, so that 
$\mathbf{P}_{x}\left( t \mapsto X_t \text{ is continuous on } [0, \zeta)\right) =1$
for $\mathcal{E}$-q.e.~$x$.
  The special case in which 
$\left( \mathcal{E}, D\left(\mathcal{E}\right)\right)$
is a regular Dirichlet form on a locally compact separable 
metric space (cf. IV.4.3(a) in \cite{ma_roeckner})
will be referred to as the regular case.
Our main interest is in the regular case, and the transfer method (Chapter~VI in 
\cite{ma_roeckner}) allows one to obtain proofs in the quasi-regular case from 
results in the regular case.   However, 
the proofs here for the regular case do not seem significantly easier, 
so direct proofs will be given under the general quasi-regular conditions.

  Let~$\hat{\mathcal{E}}$ be the form defined by
 by $\hat{\mathcal{E}}\left(x,y\right) = \mathcal{E}\left(y,x\right)$.
 There is a special standard Markov process~$\hat{X}$
 properly associated 
 with~$\hat{\mathcal{E}}$.
  For any $\alpha > 0$ the 
form $\mathcal{E}_{\alpha}$ is defined by
$\mathcal{E}_{\alpha}\left(u,v\right) = \mathcal{E}\left(u,v\right) + \alpha \left( u, v \right)$, 
where $\left( ,  \right)$ is the usual inner product on $L^{2}\left( m\right)$.
The set $D\left(\mathcal{E}\right)$ 
with the symmetric inner product $\widetilde{\mathcal{E}}_{\alpha} = (1/2) \left( \mathcal{E}_{\alpha}
+ \hat{\mathcal{E}}_{\alpha} \right)$
is a Hilbert space.
Let $\left\lVert\cdot\right\rVert_{\mathcal{E}, \alpha}$ be the norm on 
this Hilbert space, 
so that 
$\left\lVert u\right\rVert_{\mathcal{E}, \alpha}
=  \widetilde{\mathcal{E}}_{\alpha}\left(u,u\right)^{1/2}
=  \mathcal{E}_{\alpha}\left(u,u\right)^{1/2}$.
Clearly all the norms $\left\lVert\cdot\right\rVert_{\mathcal{E}, \alpha}$, $\alpha >0$, 
are equivalent. Let $K_{\alpha}$
denote a continuity constant for the weak sector condition
(equation~I.2.(2.3) in \cite{ma_roeckner}), so that 
$\left\lvert\mathcal{E}_{\alpha}\left(u,v\right)\right\rvert\
\leq K_{\alpha} \left\lVert u\right\rVert_{\mathcal{E}, \alpha} \left\lVert v\right\rVert_{\mathcal{E}, \alpha}$
for all $u,v \in D\left(\mathcal{E}\right)$. 
We will say that a sequence $u_n \in D\left(\mathcal{E}\right)$ converges $\mathcal{E}$-weakly to $u$
if it converges weakly in the Hilbert space $D\left(\mathcal{E}\right)$ with 
inner product $\widetilde{\mathcal{E}}_{\alpha}$.
The arguments used in proving I.2.12 in \cite{ma_roeckner} 
show that a sequence $u_n$ converges $\mathcal{E}$-weakly
to $u$ if and only if $\mathcal{E}_{\alpha}\left(u_n,v\right) \rightarrow 
\mathcal{E}_{\alpha}\left(u,v\right)$ for every $v \in D\left(\mathcal{E}\right)$.

The bounded and 
nonnegative $\mathcal{B}$-measurable functions
on~$E$  will be denoted by 
$b\mathcal{B}$ and $\mathcal{B}^{+}$, respectively.
Let $R_{\alpha}$, $p_{t}$
and $\hat{R}_{\alpha}$, $\hat{p}_{t}$ be the resolvents 
and Markov operators associated with $X$
and $\hat{X}$.
$\mathcal{E}$-exceptional sets, $\mathcal{E}$-quasi-everywhere properties,
and $\mathcal{E}$-quasi-continuity 
are defined in III.2.1 and III.3.2 of \cite{ma_roeckner}.

A finite measure is said to be smooth if it does not charge 
$\mathcal{E}$-exceptional sets.  Conversely, one can show that
a set is $\mathcal{E}$-exceptional if every smooth measure gives
it measure zero.  General 
smooth measures are defined in VI.2.3 of \cite{ma_roeckner},
and by VI.2.4 in \cite{ma_roeckner}, for
any smooth measure~$\mu$ there is a
unique positive continuous
additive functional $A$
such 
that $\mu$ is the Revuz measure for $A$.
For a smooth probability measure $\mu$, any martingale with respect to 
 the natural filtration of $X$ has continuous paths $\mathbf{P}_{\mu}$-almost surely
(Proposition~\ref{bounded_var_mart_constant_prop}).  This property is used
in proving convergence of stopping times.

The potentials  $G_{\alpha}$ and 
$\hat{G}_{\alpha}$ for the forms $\mathcal{E}$ and $\hat{\mathcal{E}}$
are defined in I.2.8 of \cite{ma_roeckner}
and satisfy $R_{\alpha} \, f = G_{\alpha}\, f$,
$\hat{R}_{\alpha} \, f = \hat{G}_{\alpha}\, f$, $\mathcal{E}$-q.e.,
for all $f \in L^{2}\left( m\right)$, by IV.2.9 
and IV.3.3 of \cite{ma_roeckner}.
$\mathcal{E}$-quasi-continuity is defined in II.3 of \cite{ma_roeckner}.
Each $u \in D\left(\mathcal{E}\right)$ has
$\mathcal{E}$-quasi-continuous versions, any of which is denoted by~$\tilde{u}$.
A  $\sigma$-finite measure~$\mu$ on~$\mathcal{B}$ will be said to have finite energy
if $\mu$ does not charge  $\mathcal{E}$-exceptional sets
and  the map $u \mapsto \int \tilde{u} \, d \mu$ is bounded on 
$D\left(\mathcal{E}\right)$
with respect to $\left\lVert\cdot\right\rVert_{\mathcal{E}, \alpha}$-norm
for some (and hence all) $\alpha >0$.
This condition only depends on the symmetric part of~$\mathcal{E}$.
By I.2.7 in \cite{ma_roeckner}, there exist unique elements
$v, w \in D\left(\mathcal{E}\right)$
such that $\mathcal{E}_{\alpha}\left(v,u\right) = \int \tilde{u} \, d \mu
= \mathcal{E}_{\alpha}\left(u,w\right)$ for all $u \in 
D\left(\mathcal{E}\right)$.  
$G_{\alpha}\, \mu, \hat{G}_{\alpha}\, \mu$ 
are defined to be 
$v,w$ respectively.  The new definition for potentials is consistent
with the old, 
in the sense that 
if a measure~$\mu$ has a density $f \in L^{2}\left( m\right)$
with respect to~$m$, then~$\mu$ has finite energy, and
$G_{\alpha}\, \mu = G_{\alpha}\, f$.
We will choose
$\mathcal{E}$-quasi-continuous versions of 
$G_{\alpha}\, \mu, \hat{G}_{\alpha}\, \mu$ whenever
pointwise values are needed.  By 
VI.2.1 in \cite{ma_roeckner}, the measure~$\mu$
is uniquely determined by $G_{\alpha}\, \mu$.
When~$\mu$ and~$\nu$ are measures with finite energy,
\begin{equation}
\label{background_adjoint_pot_on_measures_eqn}  
 \mathcal{E}_{\alpha}\left(G_{\alpha}\, \nu,G_{\alpha}\, \mu\right)
 =
\int \left(G_{\alpha}\, \mu\right) \, d \nu
= \mathcal{E}_{\alpha}\left(G_{\alpha}\, \mu,\hat{G}_{\alpha}\, \nu\right)
= \int \left(\hat{G}_{\alpha}\, \nu\right) \, d \mu.
\end{equation}
The definitions imply that the
resolvent equation holds for potentials of measures, 
so that in  particular $G_{\alpha}$
and $G_{\beta}$ commute.

The $\alpha$-equilibrium measure for a closed set~$B$ is defined
as the unique measure $\gamma$ 
such that $G_{\alpha}\, \gamma = 1$ holds $\mathcal{E}$-q.e.\ on~$B$ and 
$\gamma \left( B^c \right) = 0$, when such a measure exists.
The $\alpha$-equilibrium measure for $B$ will exist
 if there is an $\mathcal{E}$-quasi-continuous $v \in  D\left(\mathcal{E}\right)$ such that $v \geq 1$
holds $\mathcal{E}$-q.e.\ on~$B$ (see Section~\ref{Dirichlet_form_properties_section}). 
Define the $\alpha$-capacity of $B$
by  $\mathsf{Cap}_{\alpha}\left(B\right) = \gamma(B) = \mathcal{E}_{\alpha}\left( G_{\alpha}\,  \gamma , G_{\alpha}\,  {\gamma} \right)
= \int \left( G_{ \alpha}\,  \gamma \right) \, d \gamma$.
  The collection of closed sets which have 
$\alpha$-equilibrium measures will be denoted by
$\mathsf{C}$. 
 When $\mathcal{E}$ is a regular Dirichlet form all compact sets 
are in $\mathsf{C}$.

\subsection{A convergence theorem for random holes}
\label{random_holes_subsection}
The main convergence result, Theorem~\ref{random_sets_example_theorem}, 
deals with random sets of the
following sort.
Let  $\kappa_{n}$, $n=1,2,\ldots$ be a sequence
of positive integers with $\kappa_{n} \nearrow \infty$, 
and for each~$n$ let  $\Lambda_{j}\left( n\right)$, $j=1, \ldots, \kappa_{n}$ 
be an independent  sequence of random variables
 (not necessarily identically distributed), whose values
are compact sets in $\mathsf{C}$.
 Since identical distributions are not assumed, nonrandom sets 
 are included as  a special case.

The $\Lambda_{j}\left( n\right)$ are assumed to be measurable as maps
into the space $F(E)$ of  compact subsets of~$E$, 
equipped with the Hausdorff metric and its Borel $\sigma$-algebra.
Let
 $\Lambda\left( n\right) = \Lambda_{1}\left( n\right) \cup \cdots \cup \Lambda_{\kappa_{n}}\left( n\right)$.
 Each random set $\Lambda\left( n\right)$ provides a random environment for
 the Markov process~$X$.  Let
 $\tilde{\mathbf{P}}_n$ and $\tilde{\mathbf{E}}_n$ denote probability and 
 expectation for the probability space $\tilde{\Omega}_n$ on which
$\Lambda\left( n\right)$ is defined.  The probability space for the environment
can depend on~$n$, but for convenience in stating results,
 we will usually assume that the $\Lambda\left( n\right)$ are 
 all defined on the same space $\tilde{\Omega}$, and write $\tilde{\mathbf{P}}_n$
and $\tilde{\mathbf{E}}_n$ as $\tilde{\mathbf{P}}$ and $\tilde{\mathbf{E}}$. 
For any probability measure $\pi$ on $E$, 
it will be assumed that the entrance time $D_{ \Lambda\left( n\right) }$
and hitting time $ T_{\Lambda\left( n\right)}$
are measurable as maps from the sample space to 
the space of randomized stopping times, 
when the space of randomized stopping times is given the topology 
of stable convergence with respect to $\mathbf{P}_{\pi}$.
This measurability will hold automatically in the regular case, since in the regular case
$B \mapsto D_{B}$ is a pointwise limit of certain maps $G_k$, 
where each $G_k$ is  constant
on the sets of a measurable partition of $F(E)$.
 
  Fix $\alpha \in (0, \infty)$. 
 For each~$n$ and each~$j$, 
 let  $\gamma^{n}_{j}$ be the $\alpha$-equilibrium measure
 for $\Lambda_{j}\left( n\right)$, as defined in Section~\ref{Dirichlet_form_properties_section}. 
 It is assumed that $\gamma^{n}_{j}(A)$
 is a measurable function of the random environment
for each $A \in \mathcal{B}$, 
and that $\int \left( G_{\alpha}\,  \gamma^{n}_{j} \right)
\left( G_{ \alpha}\,  \gamma^{n}_{i} \right) \, d m$ is measurable
for all $i, j$.  In the regular case this is automatically true.
 For each $n,j$, define the averaged measure $\bar{\gamma}^{n}_{j}$ by
 $\bar{\gamma}^{n}_{j}(A) = \tilde{\mathbf{E}}\left[ \gamma^{n}_{j}(A)\right]$
 for each $A \in \mathcal{B}$.
Then $\tilde{\mathbf{E}}\left[\int f \, d\gamma^{n}_{j} \right]
 = \int f \, d \bar{\gamma}^{n}_{j}$, 
for all  $n, j$ and  $f \in b\mathcal{B} \cup \mathcal{B}^{+}$.
Let $\gamma^{n} = \sum_{j} \gamma^{n}_{j}$
 and $\bar{\gamma}^{n} =  \sum_{j} \bar{\gamma}^{n}_{j}$.

It is assumed that $\left\lVert\gamma^{n}_{j}\right\rVert_{\text{tv}} \leq \chi_{n} \text{ for all }j = 1, \ldots, \kappa_{n}$, 
where $\chi_{n}$ is a deterministic sequence of numbers with $\chi_{n} \rightarrow 0$, 
and that   
 $\sup_n \left\lVert \bar{\gamma}^{n}\right\rVert_{\text{tv}} < \infty$,
where $\left\lVert \bar{\gamma}^{n}\right\rVert_{\text{tv}} = \bar{\gamma}^{n}\left( E\right)$ is 
the total variation norm of $\bar{\gamma}^{n}$.
Since $\left\lVert G_{\alpha}\, \gamma^{n}_{j} \right\rVert_{\mathcal{E}, \alpha}^2
= \int \left( G_{\alpha}\,  \gamma^{n}_{j } \right) \, d \gamma^{n}_{j}
= \left\lVert \gamma^{n}_{j}\right\rVert_{\text{tv}} \leq \chi_{n}$, $\left\lvert \int v \, d \bar{\gamma}^{n}_{j} \right\rvert
= \left\lvert \tilde{\mathbf{E}}\left[\int v \, d \gamma^{n}_{j}\right]\right\rvert \leq K_{ \alpha}
\left\lVert v\right\rVert_{\mathcal{E}, \alpha} \sqrt{\chi_{n} }$ for any $v \in D\left( \mathcal{E}\right)$. 
Hence $\bar{\gamma}^{n}_{j}$
has finite energy.

\begin{theorem}
\label{random_sets_example_theorem}
Suppose that 
for  $\alpha \in (0, \infty)$, 
\begin{equation}
\label{kernel_density_cond_eqn}
\delta_x \, R_{\alpha} << m \text{ for }  m\text{-a.e. } x.
\end{equation}
Let $\eta$ be a finite measure with finite energy, such that 
$G_{\alpha}\,  \bar{\gamma}^{n} 
\rightarrow G_{\alpha}\,  \eta$, $\mathcal{E}$-weakly, and 
let  $A$ be the positive continuous additive functional
with Revuz measure~$\eta$.
Let $\tau_n$ be the entrance time $D_{\Lambda\left( n\right)}$ or the hitting time $ T_{\Lambda\left( n\right)}$
 for~$\Lambda\left( n\right)$.
If
\begin{equation}
\label{unequal_inner_converges_eqn}
\limsup_{n \rightarrow \infty}
\sum_{ i \neq j} 
\int \left( G_{\alpha}\,  \bar{\gamma}^{n}_{j}\right) \, d \bar{\gamma}^{n}_{i}
\leq \int \left( G_{\alpha}\, \eta \right) \, d \eta,
\end{equation}
then for any smooth probability $\pi$,  the sequence
$\tau_n \wedge \zeta$ converges $\mathbf{P}_{\pi}$-stably 
to $\tau \wedge \zeta$
in $\tilde{\mathbf{P}}$-probability, where $\tau$ 
denotes the randomized stopping time with 
$S^{\tau}_{t} = e^{ - A_t}$.
If with $\tilde{\mathbf{P}}$-probability one all the sets $\Lambda\left( n\right)$
are contained in a single compact subset of~$E$,
$\tau_n$ converges $\mathbf{P}_{\pi}$-stably 
to~$\tau$ in $\tilde{\mathbf{P}}$-probability.
\end{theorem}
In the statement of the theorem, 
stable convergence  in $\tilde{\mathbf{P}}$-probability means convergence in 
$\tilde{\mathbf{P}}$-probability 
with respect to any metric for $\mathbf{P}_{\pi}$-stable convergence.
Theorem~\ref{random_sets_example_theorem} is proved in 
Section~\ref{random_obstacle_example_section}.

Stable convergence of $\tau_n \wedge \zeta$
 implies convergence of the solutions of 
the corresponding Dirichlet problems.
This is discussed in more detail in Section~\ref{dirichlet_problems_section}.  

It follows from the resolvent equation that~\eqref{kernel_density_cond_eqn}
implies that $\delta_x \, R_{\alpha} << m$
for $\mathcal{E}$-q.e.~$x$, and thus is similar to
condition $\text{(AC)}^{\prime}$
of Definition~A.2.16 in \cite{chen_fu}.

Condition~\eqref{unequal_inner_converges_eqn}  was 
introduced  in \cite{papanicolaou_varadhan} for deterministic sets $\Lambda_{i}\left( n\right)$,
in the Brownian motion case. This condition
 ensures that
the  sets $\Lambda_{j}\left( n\right)$ are rather evenly distributed.
The proof of Theorem~\ref{random_sets_example_theorem}
shows that even 
when~\eqref{unequal_inner_converges_eqn} is not satisfied, 
$\left( \tau_n \wedge \tau \right) \wedge \zeta$ converges 
$\mathbf{P}_{\pi}$-stably 
to $\tau \wedge \zeta$ 
in $\tilde{\mathbf{P}}$-probability, 
and $\tau_n \wedge \tau \rightarrow \tau$
when all the sets $\Lambda\left( n\right)$ are contained in a compact set.
Thus asymptotically $\tau_n \geq \tau$,
verifying the physical intuition that when the sets $\Lambda_{j}\left( n\right)$
are allowed to clump together, some of their capacity may be wasted.
When $ \bar{\gamma}^{n}\left( E\right) \rightarrow 0$,
so that $A = 0$, the same limit shows that
the bodies $\Lambda\left( n\right)$ have a negligible stopping effect for large~$n$.

In the regular case, the assumption that $G_{\alpha}\,  \bar{\gamma}^{n} 
\rightarrow G_{\alpha}\,  \eta$, $\mathcal{E}$-weakly, will hold
whenever $\sup_{n} \left\lVert G_{\alpha}\,  \bar{\gamma}^{n}\right\rVert_{\mathcal{E}, \alpha} < \infty$
 and $\bar{\gamma}^{n}$ converges vaguely to~$\eta$ as a sequence of measures
(Lemma~\ref{background_vague_gives_E_weak_lemma}).

Most of the assumptions on $\gamma^{n}$ in  Theorem~\ref{random_sets_example_theorem} 
only involve the average measure $\bar{\gamma}^{n}$.
This simplifies applications, especially in the iid case.

\begin{corollary}
\label{iid_corollary_random_sets_example_theorem}
Suppose that~\eqref{kernel_density_cond_eqn} holds, 
and that the sequences $\Lambda_{j}\left( n\right)$, $j = 1, \ldots, \kappa_{n}$, are iid for each~$n$.
Assume as above that $\left\lVert \gamma^{n}_{j} \right\rVert_{\text{tv}} \leq \chi_{n}$
and $\sup_n \left\lVert \bar{\gamma}^{n}\right\rVert_{\text{tv}} < \infty$.
Let $\eta$ be a finite measure with finite energy such that for some~$\alpha$,
$\lim_{n \rightarrow \infty} \left\lVert G_{\alpha}\,  \bar{\gamma}^{n} 
- G_{\alpha}\,  \eta \right\rVert_{\mathcal{E}, \alpha} = 0$.
Then the hypotheses
of Theorem~\ref{random_sets_example_theorem} are satisfied.
\end{corollary}
\medskip \par \noindent
        \mbox{\bf Proof}\ \   
Clearly $G_{\alpha}\,  \bar{\gamma}^{n} 
\rightarrow G_{\alpha}\,  \eta$, $\mathcal{E}$-weakly,
and  $\bar{\gamma}^{n}_{j} = (1/\kappa_{n}) \bar{\gamma}^{n}$.  Thus $\sum_{ i \neq j} 
\int \left( G_{\alpha}\,  \bar{\gamma}^{n}_{j}\right) \, d \bar{\gamma}^{n}_{i}
= (1 - 1/\kappa_{n}) \left\lVert G_{\alpha}\,  \bar{\gamma}^{n} \right\rVert_{\mathcal{E}, \alpha}^2
\rightarrow \left\lVert G_{\alpha}\,  \eta \right\rVert_{\mathcal{E}, \alpha}^2$, so
\eqref{unequal_inner_converges_eqn} holds.  \hfill \mbox{\raggedright \rule{0.1in}{0.1in}}

In the setting of the classical Dirichlet problem with iid random holes
in a subset of ${\mathbb R}^{d}$, Corollary~\ref{iid_corollary_random_sets_example_theorem} 
gives a more general form of
Theorem~4.2 of \cite{balzano} (Section~\ref{example_lemmas_section}), and similarly extends
 Theorem~4.2 of \cite{balzano_notarantonio}, on the Dirichlet problem for the Laplace-Beltrami operator
on a compact Riemannian manifold.

The proof of 
Theorem~\ref{random_sets_example_theorem}
is given in Section~\ref{random_obstacle_example_section},
based on a more general convergence criterion,  Theorem~\ref{convergence_from_limit_equilibrium_measures_theorem}.
Details of the proofs for Proposition~\ref{bounded_var_mart_constant_prop} and the examples
in  Section~\ref{example_lemmas_section}
are given in~\cite{baxter_nielsen_supp}.
Transformations which simplify applying Theorem~\ref{random_sets_example_theorem}
are
discussed in Section~\ref{transformations_section}.  In particular one can use Girsanov's theorem
to deal with a drift term.

\subsection{Relaxed Dirichlet problems}
In the setting of Dirichlet forms, the solution~$u$
of a Dirichlet problem can be defined
by specifying the boundary values for~$u$
and requiring that the equation $\mathcal{E}_{\alpha}\left(u ,v\right)  = \int v f \, d m$ 
must hold for all~$v$, 
where $f$ represents the source term in the equation, 
and~$v$ lies in a suitable class of test functions.  More generally, 
$u$ is said to solve a \emph{relaxed} Dirichlet problem, with penalty measure~$\eta$,
if $\mathcal{E}_{\alpha}\left(u ,v\right) + \int u v \, d \eta = \int v f \, d m$ for 
suitable~$v$, 
where the penalty measure $\eta$ is a measure which does not charge
sets of capacity zero, but may be infinite on some sets. 
Any ordinary Dirichlet problem on a region~$U$ can be represented as a relaxed 
Dirichlet problem using a suitable infinite penalty measure~$\eta$, 
so that convergence for solutions of ordinary Dirichlet problems can be formulated
as a special case of convergence of solutions of relaxed Dirichlet problems.
General properties of relaxed Dirichlet problems in the setting of Dirichlet forms
 have been studied in a number of papers, 
including  \cite{dalmaso_capacities}, \cite{biroli_mosco}, 
\cite{dalmaso_decicco_notarantonio_tchou}, \cite{biroli_tchou}, \cite{mataloni_tchou}. 
These papers deal with a class of Dirichlet forms
which have certain extra regularity properties.
A Dirichlet form $\mathcal{E}$ in this class is regular and strongly local,
and satisfies some additional assumptions,
in particular that 
a Poincar\'{e} inequality holds for $\mathcal{E}$
and $m$ has a doubling property.   
\cite{mataloni_tchou} proves some results for nonsymmetric Dirichlet forms, 
while the other references study symmetric forms.
In the present paper the extra regularity assumptions are replaced
by~\eqref{kernel_density_cond_eqn}.
For symmetric forms satisfying the extra regularity assumptions,
\eqref{kernel_density_cond_eqn} holds, since it
 follows from equation~(1.10) 
of \cite{biroli_mosco}, using the arguments
of~Theorem~4.2.7 in \cite{fu_symm}.
We do not know if equation~\eqref{kernel_density_cond_eqn} always
holds for nonsymmetric forms under the assumptions
in~\cite{mataloni_tchou}.

The results in \cite{dalmaso_capacities}, \cite{biroli_mosco}, 
\cite{dalmaso_decicco_notarantonio_tchou}, \cite{biroli_tchou}, \cite{mataloni_tchou}
include some necessary and sufficient conditions for convergence
of solutions of general nonrandom relaxed Dirichlet problems.
Thus these results are relevant to the problems  considered
here.  However, 
further arguments would be needed in order to verify the hypotheses of these 
results and obtain convergence in that way.
Such an approach was used in  the papers \cite{balzano} and \cite{balzano_notarantonio}
mentioned earlier, dealing with the Laplacian operator
and the Laplace-Beltrami operator.

\section{Stable convergence facts}
\label{stable_convergence_facts_section}
Let $\mathcal{B}_\infty$ denote
the Borel sets on $[0, \infty]$, and for
any randomized stopping time $\tau$, 
let $\Phi^{\tau} = \Phi^{\tau}_{\omega} $ be the random measure on $\mathcal{B}_\infty$ such that 
$\Phi^{\tau}_{\omega}\left( (t, \infty] \right) = S^{\tau}_{t}(\omega)$.  
The Functional Monotone Class Theorem shows that
for any bounded $\mathcal{F} \times \mathcal{B}_\infty$-measurable $Y$, 
\begin{equation}
\label{random_measure_integral_survival_function_eqn}
\int \int_0^1 Y \left( \omega,  \tau(\omega, u)  \right) \, d u \, \mathbf{P}\left( d \omega\right) 
= \int \int_0^\infty Y \left( \omega, s\right) \, \Phi^{\tau}_{\omega}\left(  d s \right) \mathbf{P}\left( d \omega\right).
\end{equation}

\begin{lemma}
\label{extended_stable_convergence_lemma}
Let $\left(\mathcal{G}_{t}\right)_{t \geq 0}$ 
 be a filtration on a probability space
 $\left(  \Omega, \mathcal{G}, \mathbf{P}\right)$, with $\mathcal{G}_{t} \subset \mathcal{G}$.
Let 
$\tau_n, \tau$ be
randomized $\left( \mathcal{G}_{t}\right)$-stopping times 
on $ \Omega \times (0,1)$, 
such that $\tau_n \rightarrow \tau$,  $\mathbf{P}, \mathcal{G}$-stably.
Let $\lambda^{\ast}$ denote Lebesgue measure on $(0,1)$.
\begin{description}
\item{(i)}
Let~$Y$ be a real-valued process on~$ \Omega \times [0, \infty]$
which is $\mathcal{G} \times \mathcal{B}_\infty$-measurable, 
where $\mathcal{B}_\infty$ is the collection of Borel subsets of 
$[0, \infty]$.  Let~$Z$ be a nonnegative random variable 
with $\mathbf{E}\left[Z\right] < \infty$, 
and such that for $\mathbf{P}$-a.e.~$\omega$, 
$\left\lvert Y\left(\omega, t\right)\right\rvert  \leq Z(\omega)$ 
for all~$t$.
Let $A \in \mathcal{G} \times \mathcal{B}_\infty$, such that~$A$
contains all pairs $(\omega, u)$
for which the map: $t \mapsto Y(\omega, t)$ is discontinuous
at~$u$. 
If $\mathbf{P} \times \lambda^{\ast} \left( 
\left\{ (\omega, s): \ (\omega, \tau(\omega, s) ) \in A\right\}\right)= 0$,
then  
\begin{equation*}
\lim_{n \rightarrow \infty}
\int Y_{\tau_n} \, d \mathbf{P} \, d \lambda^{\ast} = \int Y_\tau \, d \mathbf{P} \, d \lambda^{\ast}
\end{equation*}
\item{(ii)}
Let~$\xi:  \Omega \times [0, \infty] \rightarrow [-\infty, \infty]$
be $\mathcal{G} \times \mathcal{B}_\infty$-measurable, 
c\`{a}dl\`{a}g\ and
quasi-left continuous with respect 
to $\mathbf{P}$ on the closed interval $[0, \infty]$.
Let $H \in L^{ 1}\left(  \mathcal{G}, \mathbf{P}\right)$ 
be such that $\mathbf{E}\left[ \left\lvert H \right\rvert\sup_{t} \left\lvert \xi_t \right\rvert \right] < \infty$.
Then
\begin{equation}
\label{xi_integral_converge_truncated_eqn}
\lim_{n \rightarrow \infty}
\int  H \xi_{\tau_n } \, d \mathbf{P} \, d \lambda^{\ast} =
\int H \xi_{\tau}  \, d \mathbf{P} \, d \lambda^{\ast}.
\end{equation}
\end{description}
\end{lemma}
Statement~(i) of Lemma~\ref{extended_stable_convergence_lemma}
follows from the Corollary to Theorem~7 in \cite{meyer}.
Statement~(ii) follows from Theorem~(1.10) in \cite{baxter_chacon_compactness}
or Theorem~5 in \cite{meyer}.  Statement~(i) is the main tool in applying
stable convergence.  The next two lemmas are simple observations based on 
the definitions.

\begin{lemma}
\label{more_extended_stable_convergence_lemma}
Let  $\left( \Omega, \mathcal{G}, \mathbf{P}\right)$
be a probability space. 
Let $\tau_n, \tau, \sigma$ be randomized times, 
meaning $\mathcal{F} \times \mathcal{B}^{\ast}$-measurable
maps from $ \Omega \times (0,1)$ to $[0, \infty]$.
\begin{description}
\item{(a)}
Suppose that 
$\tau_n  \rightarrow \tau$,
$\mathbf{P}, \mathcal{G}$-stably.
Then $\tau_n \wedge \sigma \rightarrow \tau \wedge \sigma$, 
 $\mathbf{P}, \mathcal{G}$-stably.
\item{(b)}
Let $\sigma_k$ be a sequence of randomized times
such that $\sigma_k \nearrow \sigma$ and 
$\tau_n \wedge \sigma_k \rightarrow \tau \wedge \sigma_k$,
$\mathbf{P}, \mathcal{G}$-stably as $n \rightarrow \infty$, for each~$k$.
Then 
$\tau_n \wedge \sigma \rightarrow \tau \wedge \sigma$,
$\mathbf{P}, \mathcal{G}$-stably.
\end{description}
\end{lemma}
\medskip \par \noindent
        \mbox{\bf Proof}\ \   
(i) follows at once from Lemma~\ref{extended_stable_convergence_lemma} 
and the definitions (see the corollary to Lemma~3.1 in 
\cite{baxter_dalmaso_mosco}).
To prove (ii), by choosing a subsequence we may assume
that $\tau_n$ converges stably to some limit $\psi$.
By~(i), $\tau_n \wedge \sigma $ converges stably
to $\psi \wedge \sigma$, and also
$\tau_n \wedge \sigma_k$ converges stably to $\psi \wedge \sigma_k$
for each~$k$.  Hence $\psi \wedge \sigma_k = \tau \wedge \sigma_k$
for all~$k$, 
and so $\psi \wedge \sigma = \tau \wedge \sigma$, 
independent of the choice of subsequence.  \hfill \mbox{\raggedright \rule{0.1in}{0.1in}}

It is usually sufficient to have convergence for $\tau_n \wedge \zeta$ rather than for 
$\tau_n$, as in Lemma~\ref{convergence_for_probability_solutions_lemma}.
The two forms of convergence are sometimes equivalent,
 as the next result shows.
\begin{lemma}
\label{convergence_to_time_made_infinite_lemma}
Let $\tau_n, \tau$ be randomized $\left( \mathcal{G}_{t}\right)$-stopping times.
Suppose that for any $\mathbf{P}, \mathcal{G}$-stable limit point $\tau$
of the sequence $\tau_n$,  
$\mathbf{P} \times \lambda^{\ast} \left(  \left\{ \zeta \leq \tau < \infty\right\} \right)= 0$.
Suppose also that there is a randomized 
stopping time $\sigma$
such that $\tau_n \wedge \zeta \rightarrow \sigma \wedge \zeta$, $\mathbf{P}, \mathcal{G}$-stably.
Then $\tau_n$ converges  $\mathbf{P}, \mathcal{G}$-stably to  the randomized
stopping time $\hat{\tau}$ defined by 
$\hat{\tau} = \sigma $ if $\sigma < \zeta$, $\hat{\tau} = \infty$ otherwise.
\end{lemma}

\medskip \par \noindent
        \mbox{\bf Proof}\ \   
Suppose $\tau_{n_k} \rightarrow \tau$, $\mathbf{P}, \mathcal{G}$-stably.
By Lemma~\ref{more_extended_stable_convergence_lemma} (a), 
$\tau_{n_k} \wedge \zeta \rightarrow \tau \wedge \zeta$, 
so $\tau \wedge \zeta = \sigma \wedge \zeta$ modulo $\mathbf{P} \times \lambda^{\ast}$.
Hence modulo $\mathbf{P} \times \lambda^{\ast}$ we have $\tau = \sigma$ if $\tau \wedge \zeta < \zeta$,
and if $\tau \wedge \zeta \geq \zeta$, then $\tau \geq \zeta$, so $\tau = \infty$
modulo $\mathbf{P} \times \lambda^{\ast}$.
Thus $\tau = \hat{\tau}$ modulo $\mathbf{P} \times \lambda^{\ast}$.
\hfill \mbox{\raggedright \rule{0.1in}{0.1in}}

The final lemma in this section
is only used in proving analytical consequences of stable convergence 
(Lemma~\ref{convergence_for_probability_solutions_lemma}).
\begin{lemma}
\label{stable_from_a_point_lemma}
Let $X$ be a Markov process satisfying the assumptions of Section~\ref{assumptions_X_subsection}.
Let $\mathcal{G}$ be a $\sigma$-algebra with $\mathcal{F}_{t} \subset \mathcal{G}$
for all~$t$.  Let  $A_t$ be a positive continuous additive functional, 
 and let $\tau$ be the randomized stopping time with
$ S^{\tau}_{t} = e^{ - A_t}$.
Let $\sigma$ be a terminal time, and 
let $\tau_n$ be a sequence of terminal times,
such that $\tau_n \wedge \sigma  \rightarrow \tau \wedge \sigma $, 
$\mathbf{P}_{\pi}, \mathcal{G}$-stably
for some probability measure $\pi$.
Suppose that for $\pi$-a.e.~$x$, 
$\tau_n \rightarrow \tau$, $\mathbf{P}_{\mu}, \mathcal{G}$-stably for $\mu$ with 
$\mu <<  \delta_x \, R_{\alpha}$-stably.
 (This last condition
will be true, for example,  if
$\pi$ is smooth and $\tau_n \wedge \sigma  \rightarrow \tau \wedge \sigma $, 
$\mathbf{P}_{\mu}, \mathcal{G}$-stably whenever $\mu << m$.)
Then there exists a subsequence $n_k$ such that 
for $\pi$-a.e~$x$,  $\tau_{n_k} \wedge \sigma  \rightarrow \tau \wedge \sigma  $, 
$\mathbf{P}_{x}, \mathcal{G}$-stably.
\end{lemma}
This follows from the proof of Theorem~7.1 in~\cite{nielsen}, and holds for the general 
case of stopping times such that $S^{\tau_n}_{t}$ and $\tau$ are multiplicative functionals.
For Brownian motion a different proof was given in Theorem~1.3 of \cite{baxter_jain}.

 \section{A convergence criterion}
\label{a_general_convergence_criterion_section}

Theorem~\ref{random_sets_example_theorem} will be derived from a more general convergence
result,
given below as Theorem~\ref{convergence_from_limit_equilibrium_measures_theorem}.
 The conditions for convergence in this theorem are motivated by the intuitive idea
that the limiting rate 
of hitting $\Lambda\left( n\right)$ within a neighborhood of a point 
should be determined by the 
size (i.e.\ capacity) of that part of the  $\Lambda\left( n\right)$ which lies
the neighborhood (cf. \cite{baxter_chacon_jain},  \cite{baxter_jain}
for the Brownian motion case).
Expressing this picture in terms of potentials leads to the following.
\begin{definition}
\label{limit_equilibrium_measure_short_def}
 Let~$\eta$ be a measure with finite energy.  Let $\alpha \in (0, \infty)$. 
For any sequence $\gamma_n$ of measures with finite energy, 
we will write $ \gamma_n \overset{\mathcal{E}}{\sim} \eta$ if $G_{ \alpha}\, \gamma_n
\rightarrow G_{ \alpha}\,  \eta $, $\mathcal{E}$-weakly, 
 and also 
$\lim_{n \rightarrow \infty} 
\int h \left\lvert G_{\alpha}\, \gamma_n 
- G_{\alpha}\, \eta\right\rvert   \, d m = 0$ for all $h \in L^{2}\left(  m\right)$.

Let $\Lambda\left( n\right)$ be a closed set, $n=1,2, \ldots$.
Let $\tau_n = D_{\Lambda\left( n\right)}$, the entrance 
time of $\Lambda\left( n\right)$.
Let $ \nu_n$ be a sequence of measures with finite energy,
such that  $\sup_{n} \nu_{n} \left( E\right) < \infty$,
and $\nu_n \overset{\mathcal{E}}{\sim} \eta$.  
If
\begin{equation}
\label{alpha_upper_sense_eqn}
\lim_{ n \rightarrow \infty}
\mathbf{E}_{\pi_0} \left[ 
e^{ - \alpha \tau_{n}}\left( 1 + G_{\alpha}\, \eta\left(X_{\tau_n}\right)
- G_{\alpha}\, \nu_{n}\left(X_{\tau_n} \right) \right)^{+}  \right]= 0
\end{equation}
for every smooth probability measure~$\pi$, then we will say that
$\eta$ is 
 $\alpha$-bounded from above, for the sequence $\left(\Lambda\left( n\right)\right)$, 
and $\nu_n$ will be called an $\alpha$-upper sequence for $\Lambda\left( n\right), \eta$.

If there exists a sequence of measures $ \mu_n$ with finite energy,
with  $\sup_{n} \mu_{n} \left( E\right) < \infty$ 
and
$\mu_{n} \left(\Lambda\left( n\right)^{c}\right) = 0$ for
all~$n$, 
such that $\mu_n \overset{\mathcal{E}}{\sim} \eta$ and
\begin{equation}
\label{alpha_lower_sense_eqn}
\lim_{ n \rightarrow \infty}
\mathbf{E}_{\pi} \left[ 
e^{ - \alpha \tau_{n}}\left( 1 + G_{\alpha}\, \eta\left(X_{\tau_n}\right)
 - G_{\alpha}\, \mu_{n}\left(X_{\tau_n}\right)  \right)^{-} \right]
= 0
\end{equation}
for every smooth probability measure~$\pi$, 
then~$\eta$ will be said to be $ \alpha$-bounded from below 
for the sequence $\left(\Lambda\left( n\right)\right)$, 
and $\mu_n$ will be called an $\alpha$-lower sequence for $\Lambda\left( n\right), \eta$.
\end{definition}

If~\eqref{kernel_density_cond_eqn} holds,  Lemma~\ref{background_weak_gives_L_1_lemma}
shows that
$ \gamma_n \overset{\mathcal{E}}{\sim} \eta$ automatically holds whenever
$G_{ \alpha}\, \gamma_n
\rightarrow G_{ \alpha}\,  \eta $, $\mathcal{E}$-weakly.

The proof of Theorem~\ref{random_sets_example_theorem} in Section~\ref{random_obstacle_example_section}
will show that in the setting of that theorem the conditions in 
Definition~\ref{limit_equilibrium_measure_short_def} are satisfied
in the nonrandom case, with $\Lambda\left( n\right)$ equal to the union of many small bodies $\Lambda_{i}\left( n\right)$, 
and $\nu_n, \mu_n$ each approximately equal 
to the sum of the $\alpha$-equilibrium measures
of the $\Lambda_{i}\left( n\right)$.  The same is true in the general case with probability one 
for a subsequence.

\begin{theorem}
\label{convergence_from_limit_equilibrium_measures_theorem}
Let $\Lambda\left( n\right)$ be a closed set, $n=1,2, \ldots$.
Let $\alpha \in (0, \infty)$,
and let~$\eta$ be $\alpha$-bounded from above
for $\left(\Lambda\left( n\right)\right)$.
Let $\tau$ be the randomized stopping 
time with rate measure~$\eta$.
Let  $\tau_n = D_{ \Lambda\left( n\right)}$
or $\tau_n =  T_{\Lambda\left( n\right)}$,
and let $\pi$ be a smooth probability measure.
Then $\tau_n \wedge \tau \wedge \zeta \rightarrow \tau \wedge \zeta$, 
$\mathbf{P}_{\pi}$-stably.
In the case that all the sets $\Lambda\left( n\right)$
are contained in a single compact subset of~$E$,
$\tau_n \wedge \tau \rightarrow \tau$, 
$\mathbf{P}_{\pi}$-stably.

Suppose~$\eta$ is also $ \alpha$-bounded from below
for $\left( \Lambda\left( n\right)\right)$.  Then
$\tau_n \wedge \zeta \rightarrow \tau \wedge \zeta$, 
$\mathbf{P}_{\pi}$-stably.  If
the sets $\Lambda\left( n\right)$
are contained in a single compact subset of~$E$,
$\tau_n  \rightarrow \tau$, 
$\mathbf{P}_{\pi}$-stably.
\end{theorem}
Although Theorem~\ref{convergence_from_limit_equilibrium_measures_theorem}
does not require the absolute continuity condition~\eqref{kernel_density_cond_eqn}, 
it may not be easy to show that the hypotheses are satisfied
without using condition~\eqref{kernel_density_cond_eqn}.  
Before proving the theorem a few auxiliary facts are needed.

\begin{lemma}
\label{obtaining_killing_equation_lemma}
Let $\Lambda\left( n\right)$ be a closed set, $n=1,2, \ldots$, 
and let $\tau_n = D_{\Lambda\left( n\right)}$.
Let $\alpha \in (0, \infty)$, 
and let $\pi$ be a smooth probability measure.
Let~$\eta$ be $\alpha$-bounded from above for $\left( \Lambda\left( n\right)\right)$.
Let~$A$ be the positive 
continuous additive functional with Revuz measure~$\eta$.
Let~$\tau$ be a randomized stopping time on $ \Omega \times (0,1)$ 
which is 
a $\mathbf{P}_{\pi}$-stable limit point of the sequence $\left(\tau_n\right)$.
Then $\mathbf{P}_{\pi}\left( \tau = 0 \right) = 0$. Furthermore, 
for any $a, b \in [0, \infty)$ with $a < b$,
and any bounded $\mathcal{G}_{a}$-measurable function~$H$, 
\begin{equation}
\label{killing_time_potential_integrated_ineq_eqn}
\int H \mathbf{1}_{\left\{a  < \tau \leq b \right\}} e^{- \alpha\tau} \, d \mathbf{P}_{\pi} \, d \lambda^{\ast}
 \leq
\int H   \mathbf{1}_{ \left\{ a < \tau\right\}}\int_{(a, \tau \wedge b ]}e^{- \alpha r} \, d A_r  \, d \mathbf{P}_{\pi} \, d \lambda^{\ast}.
\end{equation}
If~$\eta$ is $\alpha$-bounded from below 
for $\left(\Lambda\left( n\right)\right)$, then 
equality holds in~\eqref{killing_time_potential_integrated_ineq_eqn}.
\end{lemma}

\medskip \par \noindent
        \mbox{\bf Proof}\ \   
By approximating $\pi$, we may assume that $\pi$ has finite energy.

Let $\nu_n$ be a $\pi, \alpha$-upper sequence for $\eta, \Lambda\left( n\right)$.

We may assume that the random variable~$H$
in~\eqref{killing_time_potential_integrated_ineq_eqn}
satisfies $0 \leq H \leq 1$, 
and also,  by 
passing to a subsequence and relabelling,  that $\tau_n \rightarrow \tau$,
$\mathbf{P}_{\pi}$-stably.
For $t \in [0, \infty)$, let us say that $t$
is good if $\mathbf{P}_{ \pi} \times \lambda^{\ast} \left(  \tau = t\right) = 0$
and also $\lim_{n \rightarrow \infty} \int \left\lvert G_{\alpha}\,  \nu_n
- G_{\alpha}\,  \eta \right\rvert \, d \left(  \pi \, p_{t} \right) = 0$.
Since $\pi$ has finite energy it is easy to show
from the definitions that $ \pi \, R_{\alpha}$ has a density with respect 
to $m$ which is in $L^{2}\left(  m\right)$.
Hence $\lim_{n \rightarrow \infty} \int \left\lvert G_{\alpha}\,  \nu_n
- G_{\alpha}\,  \eta \right\rvert \, d \left(  \pi \, R_{\alpha} \right) = 0$.
Thus $\lim_{n \rightarrow \infty}
\int_0^\infty e^{ - \alpha t} \int \left\lvert G_{\alpha}\,  \nu_n
- G_{\alpha}\,  \eta \right\rvert \, d \left(  \pi \, p_{t} \right) \, \lambda( d t)= 0$,
where $\lambda$ be Lebesgue measure on $\mathbb R$.
Replacing $\nu_n$ by a suitable subsequence, we may assume that $\lambda$-a.e.\ 
$t \in [0, \infty)$ is good.

Let $a, b$ be given. For the moment assume
that $b$ is good and that either $a$ is good or $a = 0$.

Let $Y_n = \mathbf{E}_{\pi} \left[  H \mathbf{1}_{ \left\{ a \leq \tau_n\right\} } 
e^{ - \alpha \tau_n \wedge b} G_{\alpha}\, \nu_n 
\left( X_{\tau_n \wedge b}\right) \right]$, 
$K_n = \mathbf{E}_{\pi} \left[  H \mathbf{1}_{ \left\{ a \leq \tau_n\right\} } 
e^{ - \alpha a} G_{\alpha}\, \nu_n 
\left( X_{a}\right) \right]$ for $n = 1, 2, \ldots$.
By~\ref{continuous_excessive_gives_submartingale_lemma},
 $e^{ - \alpha t}  
G_{\alpha}\,  \nu_n \left( X_{ t} \right)$
is a supermartingale,  and so $Y_n \leq K_n$.

Let 
$V_n =  \mathbf{E}_{\pi} \left[  H \mathbf{1}_{ \left\{ a \leq \tau_n\right\} } 
e^{ - \alpha a} G_{\alpha}\, \eta 
\left( X_{a}\right) \right]$,
for $n = 1, 2, \ldots$,  $V_\infty =  \int  H \mathbf{1}_{ \left\{ a \leq \tau\right\} } 
e^{ - \alpha a} G_{\alpha}\, \eta 
\left( X_{a}\right) \, d \mathbf{P}_{\pi} \, d \lambda^{\ast}$.
If  $a$ is good, so that $\mathbf{P}_{\pi} \times \lambda^{\ast} \left(  \tau = a\right) = 0$, 
$\lim_{n \rightarrow \infty}
V_n
= V_{\infty}$,  by
Lemma~\ref{extended_stable_convergence_lemma} (i).
If $a = 0$, $V_n = V_\infty$ for all~$n$.

Let $\gamma$ be the distribution of $X_0$ with respect to $H \mathbf{P}_{\pi}$, 
so that for $a = 0$ we have $K_n = \int \left( G_{\alpha}\, \nu_n \right) \, d \gamma
= \mathcal{E}_{\alpha}\left( G_{ \alpha}\,  \nu_n, \hat{G}_{ \alpha}\,  \gamma \right)$
and $V_n = \int \left( G_{\alpha}\, \eta \right) \, d \gamma
= \mathcal{E}_{\alpha}\left( G_{ \alpha}\,  \eta, \hat{G}_{ \alpha}\,  \gamma \right)$.
Then 
$\lim_{n \rightarrow \infty}  \left( K_n - V_n \right) = 0$
if $a = 0$ since $\nu_n \overset{\mathcal{E}}{\sim} \eta$,
and the same limit holds by the definition of good if $a$ is good.
By~\eqref{background_additive_func_pot_energy_eqn} 
and the additive functional property,
$V_{\infty} =  
\int  H \mathbf{1}_{ \left\{ a \leq \tau\right\} } e^{ - \alpha a} \mathbf{E}_{X_{a} } \left[ \int_0^{\infty} e^{ - \alpha s} d A_s \right]   \, d \mathbf{P}_{\pi} \, d \lambda^{\ast}
=  \int  H \mathbf{1}_{ \left\{ a \leq \tau\right\} } \int_a^{\infty} e^{ - \alpha r} d A_r \, d \mathbf{P}_{\pi} \, d \lambda^{\ast}$.

For $n = 1, 2, \ldots$, let
$Z_n = \mathbf{E}_{ \pi} \left[  H \mathbf{1}_{ \left\{ a \leq \tau_n \right\} } 
e^{ - \alpha \tau_n \wedge b} G_{\alpha}\,  \eta \left( X_{\tau_n \wedge b} \right) \right]
= 
\mathbf{E}_{\pi} \left[ 
H  \mathbf{1}_{ \left\{ a \leq \tau_n \right\} }\int_{\tau_n \wedge b}^{\infty} 
 e^{ - \alpha r} \, d A_r \right]$,
where the second equality follows from~\eqref{background_additive_func_pot_energy_eqn} 
and the additive functional property.
We have
$\int_{\tau_n \wedge b}^{\infty} e^{ - \alpha v} \, d A_v
\leq 
\int_{0}^{\infty} e^{ - \alpha v} \, d A_v$ for all~$n$, and 
$\mathbf{E}_{\pi} \left[ 
 \int_{0}^{\infty} 
 e^{ - \alpha r} \, d A_r \right] = \int \left( G_{\alpha}\,  \eta\right) \, d  \pi
= \mathcal{E}_{\alpha}\left( G_{\alpha}\,  \eta,\hat{G}_{ \alpha}\,   \pi\right) < \infty$.
Thus 
Lemma~\ref{extended_stable_convergence_lemma} (i)
is applicable and $\lim_{n \rightarrow \infty} Z_n =  \int 
H  \mathbf{1}_{ \left\{ a \leq \tau \right\} }\int_{\tau \wedge b}^{\infty} 
 e^{ - \alpha r} \, d A_r \, d \mathbf{P}_{\pi} \, d \lambda^{\ast} \equiv Z_{\infty}$.

 Let  
$W_n = \mathbf{E}_{ \pi} \left[  H \mathbf{1}_{ \left\{ a \leq \tau_n  \leq b\right\}} e^{ - \alpha \tau_n } \right]$, 
for $n = 1, 2, \ldots$, $W_\infty = \int  H \mathbf{1}_{ \left\{ a \leq \tau 
 \leq b\right\}} e^{ - \alpha \tau } \, d \mathbf{P}_{ \pi} \, d \lambda^{\ast}$. 
By Lemma~\ref{extended_stable_convergence_lemma} (i),
 $\lim_{n \rightarrow \infty} W_{n } = W_{\infty}$.
For  $n= 1, 2, \ldots$ we have
\begin{multline*}
 W_n + Z_n - Y_n
= \mathbf{E}_{ \pi} \left[ H \mathbf{1}_{ \left\{ a \leq \tau_n \leq b \right\} }e^{ - \alpha \tau_n}  \left( 1 + G_{\alpha}\,  \eta -G_{\alpha}\,  \nu_n \right)\left( X_{\tau_n}\right)   \right]\\
+  \mathbf{E}_{ \pi} \left[ H \mathbf{1}_{ \left\{ a \leq \tau_n \right\} } e^{ - \alpha b}\mathbf{1}_{ \left\{ \tau_n > b\right\}}\left( G_{\alpha}\,  \eta -G_{\alpha}\,  \nu_n\right)\left( X_{b}\right)  \right] \\
\leq \mathbf{E}_{  \pi} \left[  e^{ - \alpha \tau_n} \left( 1 + G_{\alpha}\,  \eta- G_{\alpha}\,  \nu_n \right)^{+} \left( X_{\tau_n}\right) \right]
+ \int   \left\lvert G_{\alpha}\,  \eta- G_{\alpha}\,  \nu_n \right\rvert \, d \left(  \pi \, p_{b} \right).
\end{multline*}
By~\eqref{alpha_upper_sense_eqn} 
and the fact that $b$ is good, 
$\limsup_{n \rightarrow \infty} \left( W_n + Z_n - Y_n\right)
\leq  0$.  Hence $\limsup_{n \rightarrow \infty} \left(  W_n + Z_n - K_n \right) \leq  0$, 
so $ W_{\infty} \leq V_{\infty} - Z_{\infty}$.
Taking $a = 0$ and $H = 1$ gives $\int  \mathbf{1}_{ \left\{ \tau 
 \leq b\right\}} e^{ - \alpha \tau } \, d \mathbf{P}_{ \pi} \, d \lambda^{\ast} \leq \int  
\int_0^{\tau \wedge b} e^{ - \alpha r} d A_r \, d \mathbf{P}_{\pi} \, d \lambda^{\ast} 
\leq \mathbf{E}_{\pi} \left[ \int_0^{ b} e^{ - \alpha r} d A_r \right] $.  Letting $b \searrow 0$ gives
$\mathbf{P}_{ \pi} \times \lambda^{\ast} \left(  \tau = 0\right) = 0$.  For general $H$ and good $a, b$, 
$W_{\infty} \leq V_{\infty} - Z_{\infty}$ gives~\eqref{killing_time_potential_integrated_ineq_eqn}.
Both sides of~\eqref{killing_time_potential_integrated_ineq_eqn} are right continuous in~$a$
and~$b$, so~\eqref{killing_time_potential_integrated_ineq_eqn} holds for all $a, b$.

If  $\eta$ is $\alpha$-bounded from below,
let $\mu_n$ be an $\alpha$-lower sequence.  
Define quantities $Y_n, K_n, V_n, Z_n, W_n$ as above, except with 
$\nu_n$ replaced by $\mu_n$.
By Lemma~\ref{martingale_up_to_closed_support_lemma},
 $e^{ - \alpha t \wedge \tau_n }  
 G_{\alpha}\,  \mu_n \left( X_{ t \wedge \tau_n} \right)$
is a martingale,  so now $Y_n = K_n$.  All the remaining arguments
work as before, with inequalities reversed.
This gives ~\eqref{killing_time_potential_integrated_ineq_eqn} with the inequality reversed, 
so the equality case of~\eqref{killing_time_potential_integrated_ineq_eqn} holds
when $\eta$ is also $\alpha$-bounded from below
for $\left(\Lambda\left( n\right)\right)$.
\hfill \mbox{\raggedright \rule{0.1in}{0.1in}}

\begin{proposition}
\label{bounded_var_mart_constant_prop}
Let  $\mathcal{G}_{t} = \sigma\left( X_s, \, s \leq t\right)$.
There exists an exceptional set $K \in \mathcal{B}$ such that 
for any probability measure $\pi$ on $E$ with $\pi \left( K \right) = 0$, 
if  $M = \left(M_t\right)_{t \in [0, \infty)}$ is a right continuous 
$\left(\mathcal{G}_{t}\right)$-martingale
with respect to $\mathbf{P}_{\pi}$, then
\[
\mathbf{P}_{\pi}\left( t \mapsto M_t \text{ is continuous on } [0,\zeta)\right) = 1.
\] 
\end{proposition}
The proof of Proposition~\ref{bounded_var_mart_constant_prop}
is similar to the proof for the Brownian motion case.
Let $V$ denote the set of all  $Y \in L^{1}\left(  \pi\right)$, 
such that there exists a right continuous
$\left( \mathcal{F}_{t} \right)$-martingale $N^Y_t$ with 
$N^Y_t = \mathbf{E}_{\pi} \left[\left.Y \,\,\right\vert \mathcal{F}_{t}  \right]$
for each $t \in [0, \infty)$, 
and such that with $\mathbf{P}_{\pi}$-probability one, $t \mapsto N^Y_t$ is 
continuous on $[0, \zeta)$.  
One can show that $V$ is closed in $L^{1}\left(  \pi\right)$
and contains a dense class, consisting of
all 
$Y = R_{ \alpha} \,  { h}\left( X_{u}  \right)$, 
where $\alpha, u \in {\mathbb R}$ with 
$\alpha > 0$, $u \geq 0$.  Thus $V = L^1(\pi)$,
and this implies the result.  Details of the proof are given in \cite{baxter_nielsen_supp}.

\begin{lemma}
\label{cannot_converge_to_lifetime_lemma}
Let $\left(  \Omega, \mathcal{G}, \mathbf{P}\right)$ be a probability space
and let $\left( \mathcal{G}_{t}\right)_{t \geq 0}$ be a filtration with 
$\mathcal{G}_{t} \subset \mathcal{G}$ for all~$t$.
Let $E$ be a set and let $\partial$ be a point not in $E$, 
and assume that $E \cup \left\{ \partial\right\}$ has a metrizable topology.
Let $\left( X_t \right)_{ t\in [0, \infty]}$ be any $\mathcal{G} \times \mathcal{B}_\infty$-measurable
 process with 
$X_t \in E$ for $t \in [0, \infty)$ and $X_\infty = \partial$.
Assume that~$X$ is c\`{a}dl\`{a}g\ and quasi-left continuous on $[0, \infty)$
with respect to~$\mathbf{P}$.
Let $\tau_n, \tau$ be randomized $\left( \mathcal{G}_{t}\right)$-stopping times.
Let $F \subset E$ be closed in $E_{\partial}$ and 
let $B_F \in \mathcal{G}$ be such that 
$\mathbf{P} \times \lambda^{\ast} \left(  \left\{ X_{\tau_n} \notin F, \tau_n < \infty\right\} \cap \left( B_F \times (0,1) \right)\right) = 0$
for each~$n$, where $\lambda^{\ast}$ is Lebesgue measure on $(0,1)$.  
Then for any $\mathbf{P}, \mathcal{G}$-stable limit point $\tau$
of the sequence $\tau_n$,  
$\mathbf{P} \times \lambda^{\ast} \left(  \left\{ X_\tau \notin F,  \tau  < \infty\right\} \cap \left( B_F \times (0,1) \right)\right)= 0$.
\end{lemma}

\medskip \par \noindent
        \mbox{\bf Proof}\ \   
$\tau_{n_k} \rightarrow \tau$,  $\mathbf{P}, \mathcal{G}$-stably, for some 
$n_k$.
Let $s \in (0, \infty)$. Let $\varphi$ be a continuous function on $E_{\partial}$ such 
that $0 \leq \varphi \leq \mathbf{1}_{  E_{\partial} - F}$, and let $\psi$ be a continuous 
function on $[0, \infty]$ such that $0 \leq \psi \leq \mathbf{1}_{ [0,s]}$. 
Let $\xi_t = \varphi\left( X_t\right) \psi( t)$ on $B_F$, $\xi = 0$ otherwise.  
Since $\xi$ is c\`{a}dl\`{a}g\ and quasi-left continuous on $[0, \infty]$, 
$\lim_{k \rightarrow \infty} 
\int  \mathbf{1}_{B_F} \varphi\left( X_{\tau_{n_k} } \right) \psi\left( \tau_{n_k}\right) \, d \mathbf{P} \, d \lambda^{\ast}
= \int  \mathbf{1}_{B_F} \varphi\left( X_{\tau} \right) \psi\left( \tau\right) \, d \mathbf{P} \, d \lambda^{\ast}$
by Lemma~\ref{extended_stable_convergence_lemma}~(ii).
$\int  \mathbf{1}_{B_F} \varphi\left( X_{\tau_{n_k}} \right) \psi\left( \tau_{n_k}\right) \, d \mathbf{P} \, d \lambda^{\ast}
= 0$ for all~$k$, so $ \int  \mathbf{1}_{B_F} \varphi\left( X_{\tau} \right) \psi\left( \tau\right) \, d \mathbf{P} \, d \lambda^{\ast} = 0$
for all $\varphi, \psi$, which gives
$\mathbf{P} \times \lambda^{\ast} \left(  \left\{ X_\tau \notin F,  \tau  < \infty\right\} \cap \left( B_F \times (0,1) \right)\right)= 0$.
\hfill \mbox{\raggedright \rule{0.1in}{0.1in}}

\medskip \par \noindent
        \mbox{\bf Proof of 
Theorem~\ref{convergence_from_limit_equilibrium_measures_theorem}}\ \  

Suppose $\tau_n = D_{\Lambda\left( n\right)}$.  Let $\alpha \in (0, \infty)$.
Let $\eta$ be $\alpha$-bounded from above for
for $\left(\Lambda\left( n\right)\right)$.

Let $\sigma$ be any randomized stopping time which is a $\mathbf{P}_{\pi}$-stable
limit point of a subsequence $\tau_{n_i}$.  
Let $\pi$ be a smooth probability measure.
We will show that
with $\mathbf{P}_{\pi}$-probability one,
$S^{\sigma}_{ t} \geq e^{ - A_t}$ for $t \in [0, \zeta)$.
By relabelling we may assume  that the full sequence $\tau_n$
converges $\mathbf{P}_{\pi}$-stably 
to~$\sigma$.

Applying Lemma~\ref{obtaining_killing_equation_lemma}
shows that~\eqref{killing_time_potential_integrated_ineq_eqn}
holds.
Let $\Phi^{\sigma}$ be the random measure on $[0, \infty]$ such that 
$\Phi^{\sigma}{\left( (t, \infty] \right)} = S^{\sigma}_{t}$.
The left side of~\eqref{killing_time_potential_integrated_ineq_eqn} is
$\mathbf{E}_{\pi} \left[ H\int_{(a, b] } e^{- \alpha u} \,  \Phi^{\sigma}{\left(  d u  \right)}  \right]$.
The right side of~\eqref{killing_time_potential_integrated_ineq_eqn} is:
\begin{multline*}
\mathbf{E}_{\pi} \left[ H \int \int \mathbf{1}_{ \left\{ (a, \infty)\right\}}(z) 
\mathbf{1}_{(a, z \wedge b ]}(u) e^{- \alpha u}
 \, d A_u  \,  \Phi^{\sigma}{\left(  d z \right)} \right]\\
=\mathbf{E}_{\pi} \left[ H \int \int 
\mathbf{1}_{(a, b ]}(u) e^{- \alpha  u}
\mathbf{1}_{[u, \infty)}(z)
  \,  \Phi^{\sigma}{\left(  d z  \right)} \, d A_u \right]
= \mathbf{E}_{\pi} \left[ H \int_{(a, b ]} 
e^{- \alpha  u} S^{\sigma}_{u} \, d A_u \right].
\end{multline*}
Hence by~\eqref{killing_time_potential_integrated_ineq_eqn},
for all $a, b \in [0, \infty)$ with $a \leq b$, 
\begin{equation}
\label{survival_function_integrated_killing_eqn}
\mathbf{E}_{\pi} \left[ H\int_{(a, b]} e^{- \alpha u} \, \Phi^{\sigma}{\left( d u \right)}  \right]
\leq \mathbf{E}_{\pi} \left[ H  \int_{(a, b ]}
 e^{- \alpha  u} S^{\sigma}_{u} \, d A_u \right]
\end{equation}
for any bounded $\mathcal{G}_{a}$-measurable function~$H$.
Then $N_t \equiv \int_{(0, t ]} e^{- \alpha u} \, \left(\Phi^{\sigma}{\left(  d u \right)} 
- S^{\sigma}_{u} d A_u\right)$ is a $\mathcal{G}_{t}$-supermartingale.

 Let~$F$ be a compact subset of $E$, and
as usual let $D_{ E_{\partial} - F}$ be the entrance time of $E_{\partial}  - F$. 
Then 
 by Lemma~\ref{cannot_converge_to_lifetime_lemma}, with $B_F = \left\{ D_{E_{\partial} - F} \geq \zeta \right\}$, 
\begin{equation}
\label{no_finite_stop_ge_lifetime_eqn}
 \int_{ \left\{ D_{ E_{\partial} - F} \geq \zeta  \right\}  }
 \Phi^{\sigma}{\left(  [\zeta, \infty) \right)} \, d \mathbf{P}_{\pi}
= 0.
\end{equation}
Let $N^F_t = N_{ t \wedge D_{ E_{\partial} - F}}$.
It is easy to check using equation~\eqref{background_additive_func_pot_energy_eqn}
 that $N^F$ is of class (DL), so by Doob-Meyer,
$N^F_t = M_t - G_t$, where~$M$
is a right continuous martingale and $G$ 
is a natural right-continuous increasing process,
and $M$ and $G$ are unique.
Thus $M_t = M_{t \wedge D_{E_{\partial} - F} }$
and $G_t = G_{ t \wedge D_{ E_{\partial}  - F}}$.
Also $M = N^F + G$.  By Proposition~\ref{bounded_var_mart_constant_prop}, 
$t \mapsto M_t$ is continuous on $[0, \zeta)$ with probability one.

The formula for~$N$ shows that the only discontinuities in 
the paths of~$N$ and $N^F$ occur in the form of positive jumps.
The same is true for~$G$, so if a positive jump for $N_{\cdot}^F(\omega)$ occurs at time~$t$
then 
a positive jump for $M_{\cdot}(\omega)$ must also occur at time~$t$.  
With probability one, since~$M$ has continuous paths
on $[0, \zeta)$,  $N^F$ must also have continuous paths on $[0, \zeta)$.
Hence $t \mapsto N^F_t$ is continuous on $[0, \infty)$ for $\mathbf{P}_{\pi}$-a.e.\ path
 in  $\left\{ \zeta = \infty\right\} \cup \left\{ D_{ E_{\partial} - F} < \zeta  \right\}$, and also 
for $\mathbf{P}_{\pi}$-a.e.\ path in $\left\{ D_{ E_{\partial} - F} \geq \zeta  \right\} \cap \left\{ \zeta < \infty\right\}$,
by~\eqref{no_finite_stop_ge_lifetime_eqn}.
Hence $t \mapsto N^F_t$
is continuous on $[0, \infty)$  with probability one.
Since $N^F$ is continuous, $M_{\zeta - } - M_{\zeta} = G_{\zeta} - G_{\zeta-} \geq 0$.
Since $G$ is natural, $\mathbf{E}_{\pi} \left[  \left( M_{\zeta - } - M_{\zeta} \right)
\left(  G_{\zeta} - G_{\zeta-}  \right)  \right]= 0$, so $t \mapsto M_t$ is continuous at $\zeta$,
and $M$ is continuous on $[0, \infty)$.
Since the paths of $M$ have finite variation, 
they are constant on $[0, \infty)$ with probability one.
Hence $t \mapsto N^F_t$ is continuous and monotonic decreasing with probability one.
Thus $t \mapsto N_t$ is continuous and monotonic decreasing on $[0, D_{E_{\partial} - F} )$
with probability one.

Let $\left( F_k \right)$ be an $\mathcal{E}$-nest, with $F_k$ compact for each $k = 1, 2, \ldots$.
By~IV.5.30 in \cite{ma_roeckner},
$\mathbf{P}_{\pi}\left( \lim_{k \rightarrow \infty} D_{ E_{\partial} - F_k} <
 \zeta\right)= 0$.  
It follows that with probability one, $t \mapsto N_t$ is continuous and 
monotonic decreasing on $[0, \zeta)$, 
i.e.\  $\Phi^{\sigma}{\left(  d t \right)} - S^{\sigma}_{t} d \,A_t \leq 0$ 
on $[0, \zeta)$. By Lemma~\ref{obtaining_killing_equation_lemma}, 
$\mathbf{P}_{\pi}\left( S^{\sigma}_{0} = 1\right) = 1$.
Multiplying by $e^{A_t}$ shows  
$d \left( e^{ A_t} S^{\sigma}_{t}\right) \geq 0$ on $[0, \zeta)$, 
and so with $\mathbf{P}_{\pi}$-probability one,
$S^{\sigma}_{ t} \geq e^{ - A_t}$ for $t \in [0, \zeta)$.
By~\eqref{randomized_st_from_pathwise_dist_fn_eqn},
if $\sigma( \omega, r) <c < \zeta(\omega)$ then there exists $u$ such that $S^{\sigma}_{u}  \leq 1 -r$.
But then $e^{ A_u} \leq 1 - r$, so $\tau(\omega, r) <c $ by~\eqref{randomized_st_from_pathwise_dist_fn_eqn}.
It follows  that $\sigma \wedge \zeta \geq \tau \wedge \zeta$
with $\mathbf{P}_{\pi}$-probability one, 
and hence $\sigma \wedge\tau \wedge \zeta = \tau \wedge \zeta$.
Since this is true for any stable limit point $\sigma$, and compactness holds for stable convergence, 
we have proved that $\tau_n \wedge \tau \wedge \zeta \rightarrow \tau \wedge \zeta$, 
$\mathbf{P}_{\pi}$-stably.

In the special case that for some compact set~$F$ we have
$\Lambda\left( n\right) \subset F$
for all~$n$,   $\mathbf{P} \times \lambda^{\ast} \left(  \left\{ X_\sigma \notin F,  \sigma  < \infty\right\}\right)= 0$
by Lemma~\ref{cannot_converge_to_lifetime_lemma}
with $B_F =  \Omega$.
Obviously $\mathbf{P} \times \lambda^{\ast} \left(  \left\{ X_\sigma \in F, \zeta \leq  \sigma  < \infty\right\}\right)= 0$.
Hence $ \int_{ E } 
 \Phi^{\sigma}{\left(  [\zeta, \infty) \right)} \, d \mathbf{P}_{\pi}
= 0$, so  with probability one $t \mapsto S^{\sigma}_{t}$ is constant
on $[\zeta, \infty)$ and has no discontinuity at $t = \zeta$.  Hence
 with $\mathbf{P}_{\pi}$-probability one, 
$S^{\sigma}_{ t} \geq e^{ - A_t}$ for $t \in [0, \infty)$.
Using~\eqref{randomized_st_from_pathwise_dist_fn_eqn}
now shows that $\sigma \geq \tau$, so $\sigma \wedge \tau = \tau$, 
and hence $\tau_n \wedge \tau \rightarrow \tau$, 
$\mathbf{P}_{\pi}$-stably.

Up to this point we have assumed that $\eta$ is $\alpha$-bounded from above.
Now let~$\eta$ be  $\alpha$-bounded from below for
$\left( \Lambda\left( n\right)\right)$ as well, and let $\sigma$ be a $\mathbf{P}_{\pi}$-stable limit point.
Lemma~\ref{obtaining_killing_equation_lemma} now
gives the equality case of~\eqref{killing_time_potential_integrated_ineq_eqn}, 
and we can repeat the same arguments as before, replacing inequalities by
equalities at each step.
Thus the process~$N$ is now a martingale.
As before,  but omitting the Doob-Meyer decomposition step, 
one obtains $S^{\sigma}_{t} =   e^{-A_{ t} } \text{ for all } t \in [0, \zeta)$,
or
$S^{\sigma}_{t} =  e^{-A_{ t} } \text{ for all } t \in [0, \infty)$
if all the sets $\Lambda\left( n\right)$
are contained in some compact set.
By~\eqref{randomized_st_from_pathwise_dist_fn_eqn}, $\tau \wedge \zeta = \sigma \wedge \zeta$, 
and $\tau = \sigma$ if all the sets $\Lambda\left( n\right)$
are contained in some compact set.
Since $\sigma$ was any limit point,
it follows that the full sequences $\tau_n \wedge \zeta$ and $\tau_n$ converge as claimed.
This completes the proof of the convergence statements in the case that $\tau_n =D_{\Lambda\left( n\right)}$.
We know by right continuity of the process that $\mathbf{P}_{\pi}\left( \zeta = 0\right) = 0$.
Thus
 $\limsup_{ n\rightarrow \infty}
\pi\left( \Lambda\left( n\right)\right) \leq \limsup_{n \rightarrow \infty} \mathbf{P}_{\pi}\left(\tau_n = 0\right)
\leq  \limsup_{n \rightarrow \infty} \mathbf{P}_{\pi}\left(\tau_n \wedge \zeta  = 0\right)
\leq \mathbf{P}_{\pi} \times \lambda^{\ast} \left( \tau \wedge \zeta  = 0\right) = 0$, 
by what has already been proved.
Consequently $\mathbf{P}_{\pi}\left( D_{ \Lambda\left( n\right)} \neq  T_{\Lambda\left( n\right)}\right)
\rightarrow 0$, so the conclusion of the theorem remains true when 
$\tau_n =  T_{ \Lambda\left( n\right) }$.

\hfill \mbox{\raggedright \rule{0.1in}{0.1in}}

 \section{Proof of 
 Theorem~\ref{random_sets_example_theorem}} \nopagebreak
 \label{random_obstacle_example_section}

Let all the assumptions for 
Theorem~\ref{random_sets_example_theorem}
other than~\eqref{unequal_inner_converges_eqn} hold.


 We  can choose  our sample space for the 
 random variables $\Lambda_{j}\left( n\right)$ to be a product space.  That is, 
 for each $n,j$ let $\Lambda_{j}\left( n\right)$ be defined as a random variable on
some probability space  $\left( \tilde{\Omega}^n_j, q^n_j\right)$.  Let 
 \[
 \tilde{\Omega}^n = \prod_{j=1}^{\kappa_{n}} \tilde{\Omega}^n_j, 
 \ q^n = \prod_{j=1}^{\kappa_{n}} q^n_j.
 \]
For convenience, also let $\tilde{\Omega} = \prod_{n=1}^{\infty}
\tilde{\Omega}^n$, $\tilde{\mathbf{P}} = \prod_{n=1}^{\infty} q^n$.
 We will regard the random variables $\Lambda_{j}\left( n\right)$ as defined
 on either $\tilde{\Omega}^n$  or $\tilde{\Omega}$ when it seems helpful.
 By definition $\bar{\gamma}^{n}_{j}(A) = \int\gamma^{n}_{j}(A) \, d q^n_j$
 for $A \in \mathcal{B}$. 
 It is straightforward to show that 
 $ \int f \, d\gamma^{n}_{j}$ is measurable
 for $f \in b\mathcal{B} \cup \mathcal{B}^{+}$.
 It is straightforward to 
 show that $\bar{\gamma}^{n}_{j}$ is a measure, that
 \begin{equation}
 \label{measurable_function_integral_eqn}
 \int\left(\int f \, d\gamma^{n}_{j}\right) \, d q^n_j
 = \int f \, d \bar{\gamma}^{n}_{j},
 \end{equation}
 for all  $n, j$ and each  $f \in b\mathcal{B} \cup \mathcal{B}^{+}$.
For nonnegative bounded $f \in L^{1}\left(  m\right)$,
$\tilde{\mathbf{E}}\left[\int f \left( G_{\alpha}\, \gamma^{n}_{j} \right) \, dm \right]
= \tilde{\mathbf{E}}\left[\int \left( \hat{G}_{\alpha}\,  f \right) \,d \gamma^{n}_{j}\right]
= \int \left( \hat{G}_{\alpha}\,  f \right) \,d \bar{\gamma}^{n}_{j}
= \int f \left( G_{\alpha}\, \bar{\gamma}^{n}_{j}  \right) \, d m$, 
so $\int f \left( G_{\alpha}\, \bar{\gamma}^{n}_{j}  \right) \, d m 
\leq \int f \, d m$.
and it follows that
$G_{\alpha}\,  \bar{\gamma}^{n}_{j}  \leq 1$, $m$-a.e.
Hence $G_{\alpha}\,  \bar{\gamma}^{n}_{j}  \leq 1$, $\mathcal{E}$-q.e.,
by IV.3.3 of \cite{ma_roeckner}, and
$\left\lVert G_{\alpha}\, \bar{\gamma}^{n}_{j} \right\rVert_{\mathcal{E}, \alpha}^2 = \int \left( G_{\alpha}\, 
{ \bar{\gamma}^{n}_{j}}  \right) \, d \bar{\gamma}^{n}_{j} \leq \chi_{n}$.
 
For any $x$ such that $\delta_x \, R_{\alpha} << m$, and 
for each $\beta >0$, let
$g_{\beta}\left(x,  \cdot\right)$ be a
density for $\delta_x \, R_{\alpha}$ with respect to~$m$.
Otherwise let $g_{\alpha}\left(x, \cdot\right) = 0$.
Since~\eqref{kernel_density_cond_eqn} is assumed to 
hold,  we can choose  $g_{\beta}\left(x, \cdot\right)$
 (for example, via the martingale theorem)
so that $g_{\beta}\left(x, y\right)$ 
is jointly measurable in~$x$ and~$y$.
For  any $h \in L^{2}\left( m\right)$, 
\begin{equation}
\label{background_f_x_rep_pot_eqn}
G_{\beta}\, h(x) = \int g_{\beta}\left(x, y\right) h(y) \, m ( d y)
\end{equation}
for $m$-a.e.~$x$.  It follows that also
$\hat{G}_{\beta}\, h(y) = \int g_{\beta}\left(x, y\right) h(x) \, m ( d x)$
for $m$-a.e.~$y$.

\begin{lemma}
\label{choose_good_subsequence_without_energy_convergence_lemma}
Let the assumptions for Theorem~\ref{random_sets_example_theorem}
other than~\eqref{unequal_inner_converges_eqn} hold.
Let a subsequence $\left(n_i\right)$ be given. Then a 
further subsequence $(n_{i_{\ell}})$ can 
be chosen, such that with $\tilde{\mathbf{P}}$-probability one,
 $ \gamma^{n_{i_\ell} } \overset{\mathcal{E}}{\sim} \eta$ 
in the sense of Definition~\ref{limit_equilibrium_measure_short_def}.
\end{lemma}
\medskip \par \noindent
        \mbox{\bf Proof}\ \   
Since $\left\lVert \gamma^{n}_{j} \right\rVert_{\text{tv}} \leq \chi_{n}$,  $ \sum_{j} \left\lVert \gamma^{n}_{j}\right\rVert_{\text{tv}}^2
\leq \chi_{n} \sum_{j} \left\lVert \gamma^{n}_{j}\right\rVert_{\text{tv}}$.
Using independence,
\begin{multline*}
\tilde{\mathbf{E}}\left[ \left(\left\lVert \gamma^{n}\right\rVert_{\text{tv}}
-  \left\lVert\bar{\gamma}^{n}\right\rVert_{\text{tv}}\right)^2\right]
= \sum_{j=1}^{\kappa_{n}}
 \tilde{\mathbf{E}}\left[ \left(\left\lVert \gamma^{n}_{j}\right\rVert_{\text{tv}}
-  \left\lVert\bar{\gamma}^{n}_{j}\right\rVert_{\text{tv}}\right)^2\right]
\leq  4\sum_{j=1}^{\kappa_{n}}
 \tilde{\mathbf{E}}\left[ \left\lVert \gamma^{n}_{j}\right\rVert_{\text{tv}}^2\right]
\leq 4 \chi_{n} \left\lVert \bar{\gamma}^{n}\right\rVert_{\text{tv}}
 \rightarrow 0.
 \end{multline*}
Let the subsequence $\left( n_i \right)$ be given.
Using Borel-Cantelli 
 we can choose $n_{i_\ell}$ so that $\left( \left\lVert \gamma^{n_{i_\ell}}\right\rVert_{\text{tv}}
-  \left\lVert\bar{\gamma}^{n_{i_{\ell}}}\right\rVert_{\text{tv}} \right) \rightarrow 0$,  $\tilde{\mathbf{P}}$-a.e.
Since $\sup_n \left\lVert \bar{\gamma}^{n}\right\rVert_{\text{tv}} < \infty$,  
$\left\lVert\gamma^{n_{i_{\ell}}} \right\rVert_{\text{tv}}$ is bounded in~$\ell$, $\tilde{\mathbf{P}}$-a.e.

By independence, for $i \neq j$ we have
\begin{multline}
\label{envexp_inside_energy_integral_eqn}
\tilde{\mathbf{E}}\left[ \int \left( G_{\alpha}\,  \gamma^{n}_{j}\right) \, d \gamma^{n}_{i}\right]
= 
\tilde{\mathbf{E}}\left[ \int \left( G_{\alpha}\,  \gamma^{n}_{j}\right) \, d \bar{\gamma}^{n}_{i}\right]
= 
\tilde{\mathbf{E}}\left[ \int \left( \hat{G}_{\alpha}\,  \bar{\gamma}^{n}_{i}\right) \, d \gamma^{n}_{j}\right]\\
= 
 \int \left( \hat{G}_{\alpha}\,  \bar{\gamma}^{n}_{i}\right) \, d \bar{\gamma}^{n}_{j}
= 
 \int \left( G_{\alpha}\,  \bar{\gamma}^{n}_{j}\right) \, d \bar{\gamma}^{n}_{i}.
\end{multline}
If $\gamma^{n}_{i} \neq 0$, we have, using Jensen and the fact that $G_{\alpha}\,  \gamma^{n}_{j} \leq 1$, 
$\mathcal{E}$-q.e., that
\begin{multline*}
\left(\int \left( G_{\alpha}\, \gamma^{n}_{j}\right) \, d \gamma^{n}_{i} \right)^2 = 
\left\lVert \gamma^{n}_{i}\right\rVert_{\text{tv}}^2  
\left(\int \left( G_{\alpha}\, \gamma^{n}_{j}\right) \, 
d \tfrac{\gamma^{n}_{i}}{\left\lVert\gamma^{n}_{i}\right\rVert_{\text{tv}}} \right)^2 \\
\leq \left\lVert\gamma^{n}_{i}\right\rVert_{\text{tv}}^2  
\left(\int \left( G_{\alpha}\, \gamma^{n}_{j}\right)^2 \, 
d \tfrac{\gamma^{n}_{i}}{\left\lVert\gamma^{n}_{i}\right\rVert_{\text{tv}}} \right) 
\leq  \chi_{n}\int \left( G_{\alpha}\, \gamma^{n}_{j}\right) \, 
d \gamma^{n}_{i}.
\end{multline*}
When $i \neq j$, \eqref{envexp_inside_energy_integral_eqn} gives
\begin{equation}
\label{square_energy_envexp_estimate_eqn}
\tilde{\mathbf{E}}\left[ \left(\int  \left( G_{\alpha}\, \gamma^{n}_{j}\right) \, d \gamma^{n}_{i} \right)^2 \right]
\leq \chi_{n} 
 \int \left( G_{\alpha}\,  \bar{\gamma}^{n}_{j}\right) \, d \bar{\gamma}^{n}_{i}.
\end{equation}
Using Jensen, the same bound holds if $\gamma^{n}_{i}$ is replaced by $\bar{\gamma}^{n}_{i}$
and/or $\gamma^{n}_{j}$ is replaced by $\bar{\gamma}^{n}_{j}$.

For each $n$ and each $i,j \in \left\{ 1, \ldots, \kappa_{n}\right\}$, 
let 
\[
Y_{i j}(n) = \int \left( G_{\alpha}\,  \gamma^{n}_{j}- G_{\alpha}\,  \bar{\gamma}^{n}_{j} \right) \, 
d \left( \gamma^{n}_{i} - \bar{\gamma}^{n}_{i} \right)
= \int \left( \hat{G}_{\alpha}\,  \gamma^{n}_{i}- \hat{G}_{\alpha}\,  \bar{\gamma}^{n}_{i} \right) \, 
d \left( \gamma^{n}_{j} - \bar{\gamma}^{n}_{j} \right).
\]
Using~\eqref{measurable_function_integral_eqn},
 $\int Y_{i j }(n) \, d q^n_j  = 0 =
\int Y_{i j}(n) \, d q^n_i $.  Let $
Z(n) = \sum_{i \neq j} Y_{i j}(n)$.
Then 
\begin{multline*}
\int Z(n)^2 \, d q^n
= \int \left(\sum_{i \neq j} \sum_{k \neq \ell} 
Y_{i j}(n) Y_{k  \ell}(n)\right) \, d q^n\\
= \sum_{i \neq j} \int \left( Y_{i j}(n)^2 
+ Y_{i j} Y_{ j i} \right) \, d q^n
+ \sum_{i \neq j} \sum_{k,  \ell}\raisebox{.15in}{}^{\ast} \, \,
\int Y_{i j}(n) Y_{ k \ell}(n) \, d q^n,
\end{multline*}
where $\sum^\ast_{k, \ell}$
is the sum over all $k \neq \ell$ such that either~$k$ is different from~$i$ and~$j$
or~$\ell$ is different from~$i$ and~$j$.
Consider a term $\int Y_{i j}(n) Y_{ k \ell}(n) \, d q^n$ where~$k$ is different
from~$i$ and~$j$.  Evaluating the integral as an iterated integral, and integrating with 
respect to $q^n_k$ first, the value is zero.  Thus
\[
\int Z(n)^2 \, d q^n
=  \sum_{i \neq j} \int \left( Y_{i j}(n)^2 
+ Y_{i j} Y_{ j i}\right) \, d q^n.
\]
Using~\eqref{square_energy_envexp_estimate_eqn},
$\sum_{i \neq j} \int Y_{ij}(n)^2 \, d q^n\leq 4 \chi_{n} \int G_{\alpha}\,  \bar{\gamma}^{n} \, d \bar{\gamma}^{n}$.
Since $\int \left( G_{\alpha}\,  \bar{\gamma}^{n}\right) \, d \bar{\gamma}^{n}$
is bounded in~$n$, it follows that
$\int Z(n)^2 \, d q^n\rightarrow 0$ as $n \rightarrow \infty$.
We have $\int \left( G_{\alpha}\,  \gamma^{n}_{j}  \right) \, d \gamma^{n}_{i}
\leq Y_{i j} + \int \left( G_{\alpha}\,  \gamma^{n}_{j}  \right) \, d \bar{\gamma}^{n}_{i}
+ \int \left( G_{\alpha}\,  \bar{\gamma}^{n}_{j}  \right) \, d \gamma^{n}_{i}
\leq Y_{i j} + \left\lVert \bar{\gamma}^{n}_{i}\right\rVert_{\text{tv}} + \left\lVert  \gamma^{n}_{i}\right\rVert_{\text{tv}}$.
\[
\int \left( G_{\alpha}\,  \gamma^{n}\right) \, d \gamma^{n}
= \sum_{i=1}^{\kappa_{n}} \left\lVert \gamma^{n}_{i}\right\rVert_{\text{tv}}
+ \sum_{i \neq j}  \int \left( G_{\alpha}\,  \gamma^{n}_{j}\right) \, d \gamma^{n}_{i}
\leq \gamma^{n}(E)
+ Z(n) + \bar{\gamma}^{n} ( E) + \gamma^{n} (E).
\]
Hence we can refine the subsequence $n_{i_\ell}$ to ensure that
with $\tilde{\mathbf{P}}$-probability one, $
\int \left( G_{\alpha}\,  \gamma^{n_{i_\ell}}\right) \, d \gamma^{n_{i_\ell}}$
is bounded in~$\ell$.

 By IV.3.3 in \cite{ma_roeckner}, there exists 
 a  countable dense subset~$\Gamma$ of
 $D\left(\mathcal{E}\right)$.  By I.4.17 of \cite{ma_roeckner} ,
 we may choose the functions in~$\Gamma$ to be bounded, 
 and we will also choose them to be $\mathcal{E}$-quasi continuous.
 For any fixed $v \in \Gamma$, 
 let $V_n =  \int v \, d \gamma^{n}$, 
so that $\tilde{\mathbf{E}}\left[ V_n\right] = \int v \, d \bar{\gamma}^{n} = 
\mathcal{E}_{\alpha}\left(G_{\alpha}\,  \bar{\gamma}^{n},v\right) \rightarrow 
\mathcal{E}_{\alpha}\left(G_{\alpha}\,  \eta,v\right)
=\int v \,  d \eta$.
Let $V_{n,j} = \int v \, d \gamma^{n}_{j}$.
Using independence, 
\[
\tilde{{\bf Var}}\left(V_n\right)
= \sum_{j=1}^{ \kappa_{n}} \tilde{{\bf Var}}\left(V_{n,j}\right)
 \leq 4 \sum_{j=1}^{\kappa_{n}}
\tilde{\mathbf{E}}\left[ V_{n,j}^2\right] \leq
4 \sum_{j=1}^{\kappa_{n}} \left\lVert v\right\rVert_{\text{sup}}^2 \tilde{\mathbf{E}}\left[ \left\lVert \gamma^{n}_{j} \right\rVert_{\text{tv}}^2 \right]
\leq
4 \left\lVert v\right\rVert_{\text{sup}}^2  \chi_{n} \left\lVert \bar{\gamma}^{n}\right\rVert_{\text{tv}} \rightarrow 0.
\]
 Hence
we
can refine the subsequence  $n_{i_{\ell}}$ so that
with $\tilde{\mathbf{P}}$-probability one, 
$\int v \, d \gamma^{n_{i_\ell}} \rightarrow 
 \int v \, d \eta$ for each $v \in \Gamma$.
 Then with $\tilde{\mathbf{P}}$-probability one, 
 $G_{\alpha}\,  {\gamma^{n_{i_{\ell}}} }
\rightarrow G_{\alpha}\, \eta$,  $\mathcal{E}$-weakly.
By Lemma~\ref{background_vague_gives_E_weak_lemma}, 
$\bar{\gamma}^{n_{i_\ell} } \overset{\mathcal{E}}{\sim} \eta$.
\hfill \mbox{\raggedright \rule{0.1in}{0.1in}}

\begin{lemma}
\label{upper_measure_without_energy_convergence_lemma}
Let the assumptions for Theorem~\ref{random_sets_example_theorem}
other than~\eqref{unequal_inner_converges_eqn} hold.
Let a subsequence $\left(n_{i_\ell}\right)$ be given.
Let a particular environment $\tilde{\omega} \in \tilde{\Omega}$ be such that
the corresponding sequence of measures $\gamma^{n_{i_\ell}}$
has the properties stated in Lemma~\ref{choose_good_subsequence_without_energy_convergence_lemma}.
 Then  $\eta$
is $\alpha$-bounded from above for $\left(  \Lambda\left( n_{i_\ell}\right)( \tilde{\omega})\right)$, 
in the sense of Definition~\ref{limit_equilibrium_measure_short_def}.
\end{lemma}

\medskip \par \noindent
        \mbox{\bf Proof}\ \   
By assumption $\gamma^{n_{i_{\ell}}}\overset{\mathcal{E}}{\sim} \eta$. 
Let $\pi$ be a smooth probability measure. 
By approximating $\pi$ we may assume that $\pi$ has finite energy.
We must verify~\eqref{alpha_upper_sense_eqn}.

Let $\varepsilon > 0$ be given.
For any $\beta \in (0, \infty)$, let $H(\beta) = 
\left\{
\beta
G_{\beta + \alpha}\, G_{\alpha}\, 
{\eta}+ \varepsilon 
\geq G_{\alpha}\,  \eta\right\}$.
Since $\int \left( \hat{G}_{\alpha}\, \pi\right)  \, d \eta < \infty$
and $\int  \left( \hat{G}_{ \alpha}\,  \eta \right) \, d \eta < \infty$, 
\eqref{converge_up_to_excessive_eqn} tells us that for all sufficiently large $\beta$ we have
$\int_{H ( \beta)^c} \left( \hat{G}_{ \alpha}\,  {\pi}\right) \, d  \eta < \varepsilon$
and $\int_{H ( \beta)^c} \left( \hat{G}_{ \alpha}\,  {\eta}\right) \, d  \eta < \varepsilon$.
Choose such a~$\beta$, and let
$\psi^1_\varepsilon = \mathbf{1}_{H(\beta)^c} \eta$.  Then
$\int \left( \hat{G}_{ \alpha}\,  {\pi}\right) \, d  \psi_\varepsilon^1 < \varepsilon$
and $\int \left( \hat{G}_{ \alpha}\,  {\eta}\right) \, d  \psi_\varepsilon^1 < \varepsilon$.
By the domination principle ( Lemma~\ref{background_domination_principle_lemma}), 
$
\beta
G_{\beta + \alpha}\, G_{\alpha}\, 
{\eta}+ \varepsilon
\geq G_{\alpha}\, \left( \eta - \psi^1_\varepsilon \right)$
holds $\mathcal{E}$-q.e.\ on~$E$. 
Thus $
\beta
G_{\beta + \alpha}\, G_{\alpha}\, 
{\eta} + \varepsilon + G_{\alpha}\,  \psi^1_\varepsilon \geq G_{\alpha}\, \eta$, 
$\mathcal{E}$-q.e.

Since $\int \left( \hat{G}_{\alpha}\,  \pi \right) \left(
G_{ \alpha }\,  \eta \right) \, d m < \infty$, by~\eqref{converge_up_to_excessive_eqn}
we have $\lim_{\lambda  \rightarrow \infty}
\int \left( \hat{G}_{\alpha}\,  \pi \right) \left(G_{\alpha}\,  \eta
- \lambda G_{ \lambda + \alpha }\,  
G_{ \alpha }\,  \eta  \right) \, d m = 0$. 
Choose $\lambda \geq \alpha$ such that 
$\beta \int \left( \hat{G}_{\alpha}\,  \pi \right) \left(G_{\alpha}\,  \eta
- \lambda G_{ \lambda + \alpha }\,  
G_{ \alpha }\,  \eta  \right) \, d m < \varepsilon$
and 
$\beta \int \left( \hat{G}_{\alpha}\,  \eta \right) \left(G_{\alpha}\,  \eta
- \lambda G_{ \lambda + \alpha }\,  
G_{ \alpha }\,  \eta  \right) \, d m < \varepsilon$. 
Let $\psi^2_\varepsilon = \beta 
\left(G_{\alpha}\,  \eta
- \lambda G_{ \lambda + \alpha }\,  
G_{ \alpha }\,  \eta  \right) \, m$.
Then $\int  \left( \hat{G}_{\alpha}\,  \pi \right) \, d \psi^2_\varepsilon < \varepsilon$
and
$\int  \left( \hat{G}_{\alpha}\,  \eta \right) \, d \psi^2_\varepsilon < \varepsilon$.
Also $G_{ \beta + \alpha}\,  { \psi^2_\varepsilon}
= \beta G_{ \beta + \alpha}\,  G_{ \alpha}\,  \eta 
- \beta \lambda G_{ \beta + \alpha}\,  { G_{ \lambda  + \alpha}\,  { 
G_{ \alpha}\,  \eta } }$, so 
$\beta G_{ \beta + \alpha}\,  G_{ \alpha}\,  \eta 
= \beta \lambda G_{ \beta + \alpha}\,  { G_{ \lambda  + \alpha}\,  { 
G_{ \alpha}\,  \eta } } + G_{ \beta + \alpha}\,  \psi^2_\varepsilon$.
Hence 
$\beta \lambda G_{ \beta  + \alpha}\,  { G_{\lambda + \alpha}\,  G_{ \alpha}\,  \eta  }
+ \varepsilon + G_{ \alpha}\,  { \left( \psi^1_\varepsilon + \psi^2_\varepsilon\right)}
 \geq G_{ \alpha}\,  \eta $, $\mathcal{E}$-q.e.

For any $A \in \mathcal{B}$ with 
$m\left(A\right) < \infty$, 
and any number $c >0$, for any $h \in L^{2}\left( m\right)$ let
\[
M_{A, c} h (x)  =   \int \mathbf{1}_{A}(x)
\mathbf{1}_{A}(y) ( c\wedge ( \lambda g_{\lambda + \alpha}\left(x, y\right)))
 h(y) \, m( d y).
 \]

As $A \nearrow E$ and $c \nearrow \infty$, 
$M_{ A, c} G_{ \alpha}\,  \eta \nearrow 
\lambda  G_{\lambda + \alpha}\,  G_{ \alpha}\, \eta  
\leq  G_{ \alpha}\, \eta  $, $m$-a.e.
We can choose $A = A_{\varepsilon}$ and $c = c_{\varepsilon}$ such 
that 
$\beta \int \left( \hat{G}_{ \alpha}\,  {\pi}\right) 
\left( \lambda  G_{\lambda + \alpha}\,  G_{ \alpha}\, \eta 
-  M_{ A, c} G_{ \alpha}\,  \eta \right) \, d m < \varepsilon$
and 
$\beta \int \left( \hat{G}_{\alpha}\,  \eta\right) 
\left( \lambda  G_{\lambda + \alpha}\,  G_{ \alpha}\, \eta 
-  M_{ A, c} G_{ \alpha}\,  \eta \right) \, d m < \varepsilon$.

Let $\psi^3_\varepsilon =  \beta \left( \lambda  G_{\lambda + \alpha}\,  G_{ \alpha}\, \eta 
- M_{ A, c} G_{ \alpha}\,  \eta\right) \, m$.
Then $\int \left( \hat{G}_{ \alpha}\,  \pi \right) \, d \psi^3_\varepsilon < \varepsilon$
and $\int \left( \hat{G}_{ \alpha}\,  \eta \right) \, d \psi^3_\varepsilon < \varepsilon$.
Also $G_{ \beta + \alpha}\,  \psi^3_\varepsilon
= \beta  \lambda G_{ \beta + \alpha}\,  G_{ \lambda + \alpha}\,  
G_{ \alpha}\,  \eta  -  \beta G_{ \beta + \alpha}\, 
 M_{A, c} G_{ \alpha}\,  \eta $, 
so $\beta  \lambda G_{ \beta + \alpha}\,  G_{ \lambda + \alpha}\,  
G_{ \alpha}\,  \eta  
= 
\beta G_{ \beta + \alpha}\, 
 M_{A, c} G_{ \alpha}\,  \eta  + G_{ \beta + \alpha}\,  \psi^3_\varepsilon$.
Hence $\beta G_{ \beta + \alpha}\, 
 M_{A, c} G_{ \alpha}\,  \eta  
+ \varepsilon + G_{ \alpha}\,  { \left( \psi^1_\varepsilon
+ \psi^2_\varepsilon + \psi^3_\varepsilon \right)} \geq
G_{ \alpha}\,  \eta$, $\mathcal{E}$-q.e.
Let $\psi_\varepsilon = \psi^1_\varepsilon + \psi^2_\varepsilon + \psi^3_\varepsilon$, 
so that $\beta G_{ \beta + \alpha}\, 
 M_{A, c} G_{ \alpha}\,  \eta  
+ \varepsilon + G_{ \alpha}\,  { \psi_\varepsilon} \geq
G_{ \alpha}\,  \eta$, $\mathcal{E}$-q.e., 
where 
\begin{equation}
\label{psi_varepsilon_bound_eqn}
\int \left( \hat{G}_{ \alpha}\,  \pi \right) \, d \psi_\varepsilon < 3 \varepsilon, \ 
\int \left( \hat{G}_{ \alpha}\,  \eta \right) \, d \psi_\varepsilon < 3 \varepsilon.
\end{equation}

Since $\left\lVert M_{ A, c} \left( G_{\alpha}\, \gamma^{n}
-  G_{ \alpha}\,  \eta  \right)\right\rVert_{\text{sup}}
\leq c \int_{A} \left\lvert G_{\alpha}\, \gamma^{n}
- G_{ \alpha}\,  \eta \right\rvert  \, d m$, we have
$\left\lVert \beta G_{\beta + \alpha}\,  {  M_{ A, c} \left( G_{\alpha}\, \gamma^{n}
-  G_{ \alpha}\,  \eta  \right)} \right\rVert_{\text{sup}}
= \left\lVert \beta R_{\beta + \alpha} \,  {  M_{ A, c} \left( G_{\alpha}\, \gamma^{n}
-  G_{ \alpha}\,  \eta  \right)} \right\rVert_{\text{sup}} \leq  c \int_{A} \left\lvert G_{\alpha}\, \gamma^{n}
- G_{ \alpha}\,  \eta \right\rvert  \, d m$, 
$\mathcal{E}$-q.e.  Hence
\[
 \beta G_{\beta + \alpha}\,  M_{ A, c} G_{\alpha}\, \gamma^{n}  
+ c \int_{A} \left\lvert G_{\alpha}\, \gamma^{n}
- G_{ \alpha}\,  \eta \right\rvert  \, d m + \varepsilon + 
G_{\alpha}\,  \psi_\varepsilon  
\geq  G_{ \alpha}\,  { \eta}, \text{ $\mathcal{E}$-q.e.}.
\] 
For any $j = 1, \ldots, \kappa_{n}$, 
\[
\int \mathbf{1}_{ A }\left( G_{  \alpha}\,  \gamma^{n}_{j}\right) 
\, d m
= \int \left( \hat{G}_{ \alpha}\,  \mathbf{1}_{ A} \right) \, d \gamma^{n}_{j}
= \int \left( \hat{R}_{ \alpha} \,  \mathbf{1}_{ A} \right) \, d \gamma^{n}_{j}
\leq \frac{\left\lVert \gamma^{n}_{j}\right\rVert_{\text{tv}}}{ \alpha} \leq \frac{ \chi_{n}}{\alpha}.
\]
Thus $M_{A, c} 
G_{ \alpha}\,  { \gamma^{n}_{j}}
\leq c \chi_{n}/\alpha$ everywhere
and hence $\beta G_{\beta + \alpha}\,  
M_{A, c} G_{\alpha}\,  \gamma^{n}_{j}  
= \beta R_{\beta + \alpha} \,  
M_{A, c} G_{\alpha}\,  \gamma^{n}_{j}  
\leq c \chi_{n}/\alpha$, $\mathcal{E}$-q.e.

For $j = 1, \ldots, \kappa_{n}$, 
let $\breve{\gamma}^{n}_{j}   = \sum_{r \neq j} \gamma^{n}_{r}$.  
Then  $G_{\alpha}\,  \breve{\gamma}^{n }_{j}   
+ \tfrac{ c \chi_{n}}{ \alpha} \geq \beta \lambda  G_{\beta + \alpha}\, G_{ \lambda + \alpha}\,   {
G_{ \alpha}\,  { \breve{\gamma}^{n }_{j}    } } + \tfrac{ c \chi_{n}}{ \alpha}
\geq \beta G_{\beta + \alpha}\, M_{A, c} 
G_{ \alpha}\,  { \breve{\gamma}^{n }_{j}   } + \tfrac{ c \chi_{n}}{ \alpha}
\geq \beta 
G_{\beta + \alpha}\, M_{A, c} G_{ \alpha}\,  { \gamma^{n }}  $, $\mathcal{E}$-q.e.  It follows that 
\begin{equation}
\label{targ_except_inequal_eqn}
G_{\alpha}\,  \breve{\gamma}^{n }_{j}    + \frac{ c_\varepsilon \chi_{n}}{ \alpha}
+ c_\varepsilon \int_{A} \left\lvert G_{\alpha}\, \gamma^{n}
- G_{ \alpha}\,  \eta \right\rvert  \, d m + \varepsilon + G_{\alpha}\,  \psi_\varepsilon 
\geq G_{ \alpha}\,  \eta \text{  $\mathcal{E}$-q.e.}
\end{equation}

On $\Lambda_{j}\left( n\right)$, 
$G_{\alpha}\,  \gamma^{ n}
= G_{\alpha}\,  \breve{\gamma}^{ n}_{j}   + 1$, $\mathcal{E}$-q.e., and so
on every $\Lambda_{j}\left( n\right)$, 
\begin{equation}
\label{helping_measure_psi_ineq_eqn}
G_{\alpha}\,  \gamma^{ n }  -1 + 
\frac{ c_\varepsilon \chi_{n} }{\alpha } + c_\varepsilon \int_{A}
 \left\lvert G_{\alpha}\, \gamma^{n}
- G_{ \alpha}\,  \eta \right\rvert  \, d m 
+ \varepsilon + G_{\alpha}\,  \psi_\varepsilon  
\geq  G_{ \alpha}\,  { \eta},  \text{ $\mathcal{E}$-q.e.}
\end{equation}
By~\eqref{helping_measure_psi_ineq_eqn}, 
$\left( 1 +  G_{\alpha}\,  \eta  - G_{\alpha}\,  \gamma^{n_{i_\ell}}
 \right)^{+} \leq 2 \varepsilon + G_{\alpha}\,  \psi_\varepsilon$, 
$\mathcal{E}$-q.e. on $\Lambda\left( n\right)$, for large~$\ell$.

It is not hard to show that
$e^{ - \alpha t} G_{ \alpha}\, \psi_\varepsilon \left( X_t\right)$ 
is a supermartingale
with respect to $ \mathbf{P}_{ \pi}$.  Hence
$\mathbf{E}_{\pi} \left[  G_{\alpha}\, \psi_\varepsilon \left( X_0\right) \right]
\geq \mathbf{E}_{\pi} \left[  e^{ - \alpha \tau_{n_{i_\ell}}}
G_{\alpha}\, \psi_\varepsilon \left( X_{\tau_{n_{i_\ell}} }\right) \right]$,
so $\mathbf{E}_{\pi} \left[   e^{ - \alpha \tau_{n_{i_\ell}}}
G_{\alpha}\,  {\psi_\varepsilon} \left( X_{\tau_{n_{i_\ell}} }\right) \right]
\leq \int \left( G_{\alpha}\,  \psi_\varepsilon  \right) \, d \pi
< 3 \varepsilon$ by~\eqref{psi_varepsilon_bound_eqn}.
It follows that~\eqref{alpha_upper_sense_eqn}
holds with $n$ replaced by $n_{i_{\ell}}$
and~$\nu_n$ replaced by $\gamma^{n_{i_{\ell}}}$, 
so $\eta$ is a $\alpha$-bounded from above for the sets $\Lambda\left( n_{i_\ell}\right)$. 

\hfill \mbox{\raggedright \rule{0.1in}{0.1in}}

\medskip \par \noindent
        \mbox{\bf Proof of Theorem~\ref{random_sets_example_theorem}}\ \ 

Lemma~\ref{upper_measure_without_energy_convergence_lemma} 
and Theorem~\ref{convergence_from_limit_equilibrium_measures_theorem}
prove
the statements of Theorem~\ref{random_sets_example_theorem}
for which~\eqref{unequal_inner_converges_eqn}
is not assumed to hold. Suppose now
that~\eqref{unequal_inner_converges_eqn} holds, 
along with the other assumptions for Theorem~\ref{random_sets_example_theorem}.

Let a subsequence $\left( n_i \right)$ be given.
By
Lemma~\ref{choose_good_subsequence_without_energy_convergence_lemma}
 we can choose $n_{i_\ell}$ so that   with $\tilde{\mathbf{P}}$-probability one,  
$\sup_{\ell} \left\lVert\gamma^{n_{i_\ell}}\right\rVert_{\text{tv}} < \infty$
and $\gamma^{n_{i_{\ell}}} \overset{\mathcal{E}}{\sim}
\eta$.
In particular, $\int\left( G_{ \alpha}\,  \eta \right) 
\, d \bar{\gamma}^{n} \rightarrow \int\left( G_{ \alpha}\,  \eta \right) 
\, d \eta$.  Thus by~\eqref{unequal_inner_converges_eqn}, 
\[
\limsup_{n \rightarrow \infty} \left(
\sum_{ r \neq j} \int \left( G_{\alpha}\,  \bar{\gamma}^{n}_{r}\right) 
\, d \bar{\gamma}^{n}_{j} - \int\left( G_{ \alpha}\,  \eta \right) 
\, d \bar{\gamma}^{n} \right) \leq  0.
\]
We have
$ \tilde{\mathbf{E}}\left[ \int \left(G_{\alpha}\,  \eta \right)\, d \gamma^{n}\right]
= \int \left( G_{\alpha}\, \eta \right) \, d \bar{\gamma}^{n}$.
Using that and~\eqref{envexp_inside_energy_integral_eqn}, 
\[
\limsup_{n \rightarrow \infty}
\tilde{\mathbf{E}}\left[ 
\sum_{r \neq j} 
\int \left( G_{\alpha}\,  \gamma^{n}_{s}  \right)\, d \gamma^{n}_{j}
- \int  \left( G_{ \alpha}\,  { \eta} \right) \, d \gamma^{n} \right]
\leq 0.
\]

Let $\varepsilon >0$ be given.
Let $\psi_{\varepsilon}$ be the measure defined in the proof of 
Lemma~\ref{upper_measure_without_energy_convergence_lemma}, 
so that~\eqref{targ_except_inequal_eqn} holds.
Let $A_\epsilon, c_\varepsilon$ be the quantities $A, c$ appearing in~\eqref{targ_except_inequal_eqn}.
Let $e_\epsilon(n) = 
\tfrac{ c_\varepsilon \chi_{n} }{\alpha } + c_\varepsilon \int_{A_\varepsilon}
 \left\lvert G_{\alpha}\, \gamma^{n}
- G_{ \alpha}\,  \eta \right\rvert  \, d m 
+  \varepsilon$.  Then
$G_{\alpha}\,  \breve{\gamma}^{n }_{j}    + e_{\varepsilon}\left(n_{i_\ell}\right)
+ G_{\alpha}\,  \psi_\varepsilon 
\geq G_{ \alpha}\,  \eta$, $\mathcal{E}$-q.e.,  
so 
\begin{multline*}
\tilde{\mathbf{E}}\left[ 
\sum_{j} \int  \left(  G_{\alpha}\,  \breve{\gamma}^{n_{i_\ell}}_{j}   - G_{\alpha}\,  \eta  \right)^{+} 
\, d \gamma^{n_{i_\ell}}_{j}\right]\\
= \tilde{\mathbf{E}}\left[ 
\sum_{j} \int  \left(  G_{\alpha}\,  \breve{\gamma}^{n_{i_\ell}}_{j}   + 
e_\varepsilon\left(n_{i_\ell} \right) +   G_{\alpha}\,  \psi_\varepsilon 
- G_{\alpha}\,  \eta - e_\varepsilon\left(n_{i_\ell} \right) -   G_{\alpha}\,  \psi_\varepsilon   \right)^{+} 
\, d \gamma^{n}_{j}  \right]\\
\leq \tilde{\mathbf{E}}\left[ \sum_{j} \int \left( \left( G_{\alpha}\,  \breve{\gamma}^{n_{i_\ell}}_{j}   + 
e_\varepsilon\left(n_{i} \right) +   G_{\alpha}\,  \psi_\varepsilon 
- G_{\alpha}\,  \eta \right)  + 
e_\varepsilon\left(n_{i_\ell} \right) +   G_{\alpha}\,  \psi_\varepsilon   \right)
\, d \gamma^{n}_{j}  \right]\\
\leq \tilde{\mathbf{E}}\left[ 
\sum_{r \neq j} 
\int \left( G_{\alpha}\,  \gamma^{n}_{s}  \right)\, d \gamma^{n}_{j}
- \int  \left( G_{ \alpha}\,  { \eta} \right) \, d \gamma^{n} \right] 
+ 2 e_\varepsilon\left(n_{i_\ell} \right) \left\lVert\bar{\gamma}^{n}\right\rVert_{\text{tv}}
 +  2 \int \left(  G_{\alpha}\,  \psi_\varepsilon   \right)
\, d \bar{\gamma}^{n}  
\end{multline*}
Since $\limsup_{\ell \rightarrow \infty} e_\varepsilon \left( n_{i_\ell} \right) \leq \varepsilon$, 
and $\int \left(  G_{\alpha}\,  \psi_\varepsilon   \right)
\, d \bar{\gamma}^{n} \rightarrow \int \left(  G_{\alpha}\,  \psi_\varepsilon   \right)
\, d \eta $,
\eqref{psi_varepsilon_bound_eqn} gives
\[
\limsup_{\ell \rightarrow \infty} 
\tilde{\mathbf{E}}\left[ 
\sum_{j} \int  \left(  G_{\alpha}\,  \breve{\gamma}^{n_{i_\ell}}_{j}   - G_{\alpha}\,  \eta  \right)^{+} 
\, d \gamma^{n_{i_\ell}}_{j}\right]
\leq 2 \varepsilon \sup_n \left\lVert\bar{\gamma}^{n}\right\rVert_{\text{tv}}
+ 6 \varepsilon.
\]
Since $\varepsilon > 0$ is arbitrary, 
$\lim_{\ell \rightarrow \infty}
 \tilde{\mathbf{E}}\left[ 
\sum_{j} \int  \left(  G_{\alpha}\,  \breve{\gamma}^{n_{i_\ell}}_{j}   - G_{\alpha}\,  \eta  \right)^{+} 
\, d \gamma^{n_{i_\ell}}_{j}\right] = 0$.
That is, $\lim_{\ell \rightarrow \infty}
 \tilde{\mathbf{E}}\left[ 
\int  \left(  G_{\alpha}\,  \gamma^{n_{i_\ell}} - 1 - G_{\alpha}\,  \eta   \right)^{+}
\, d \gamma^{n_{i_\ell}}\right] = 0$.
By Borel-Cantelli, we can refine the subsequence $n_{i_\ell}$
to ensure that $\int  \left(  G_{\alpha}\,  \gamma^{n_{i_\ell}} - 1 - G_{\alpha}\,  \eta   \right)^{+}
\, d \gamma^{n_{i_\ell}}  \rightarrow 0$ holds with $\tilde{\mathbf{P}}$-probability one.
From now on we deal with an environment $\tilde{\omega}$ such that 
$\int  \left(  G_{\alpha}\,  \gamma^{n_{i_\ell}} - 1 - G_{\alpha}\,  \eta   \right)^{+}
\, d \gamma^{n_{i_\ell}}  \rightarrow 0$.
We will show that 
$\eta$ is $\alpha$-bounded from below
for $ \Lambda\left(  n_{i_{\ell}}\right)$.
 
Let $H_{\ell} = \left\{  
\left(  G_{\alpha}\, \gamma^{n_{i_{\ell}}}
- 1 - G_{\alpha}\, \eta  \right)^{+} \leq \delta_\ell\right\}$, 
where $\delta_\ell \rightarrow 0$ is chosen so 
that  $\gamma^{n_{i_{\ell}}} \left( H_{\ell}^c\right) 
\rightarrow 0$.
Let $\mu_{\ell}^{\prime} = \mathbf{1}_{ H_{\ell}} \gamma^{n_{i_{\ell}}}$,
so that
we have $ G_{\alpha}\, \mu_{\ell}^{\prime} \leq 1 +G_{\alpha}\, \eta
+ \delta_\ell$, $\mu_{\ell}^{\prime}$-a.e.  By the domination principle,
$ G_{\alpha}\, \mu_{\ell}^{\prime} \leq 1 +G_{\alpha}\, \eta
+ \delta_\ell$
holds $\mathcal{E}$-q.e.\ on~$E$. 
Since $\gamma^{n_{i_\ell}} \overset{\mathcal{E}}{\sim} \eta$
and $\left\lVert \gamma^{n_{i_\ell}} - \mu_{\ell}^{\prime}\right\rVert_{\text{tv}} \rightarrow 0$,
it is easy to check from the definition that 
 $\mu_\ell^{\prime} \overset{\mathcal{E}}{\sim} \eta$.

Let $\tau_n = D_{\Lambda\left( n\right)}$.
Let $\pi$ be any smooth probability measure.
Since 
\[
\mathbf{E}_{\pi} \left[
e^{ - \alpha \tau_{n_{i_{\ell}}}} \left( G_{\alpha}\, \mu_\ell^{\prime} 
\left(X_{\tau_{n_{i_{\ell}}}}\right) -G_{\alpha}\, \eta\left(X_{\tau_{n_{i_{\ell}}}}\right)- 1 \right)^{+} \right]
\leq \delta_\ell \pi\left( E\right) \rightarrow 0,
\] 
\eqref{alpha_lower_sense_eqn} holds,
with $\mu_n$ in that equation replaced by $\mu_\ell^{\prime}$.
Thus  with $\tilde{\mathbf{P}}$-probability one, 
$\eta$ is $\alpha$-bounded from below for the sequence $\Lambda\left( n_{i_{\ell} } \right)$.
By Lemma~\ref{upper_measure_without_energy_convergence_lemma},
we already know that
$\eta$ is a  $\alpha$-bounded from above for $\Lambda\left( n_{i_\ell}\right) (\tilde{\omega})$
with $\tilde{\mathbf{P}}$-probability one,
so Theorem~\ref{convergence_from_limit_equilibrium_measures_theorem}
applies.

\hfill \mbox{\raggedright \rule{0.1in}{0.1in}}

\section{Dirichlet problems}
\label{dirichlet_problems_section}
Let $X$ be a right Markov process  with state space~$E$,
cemetery point~$\partial$ and lifetime~$\zeta$.
As in Section~\ref{intro_section},
let~$U$ be an open subset of $E$ and $\sigma = D_{E_{\partial} - U}$, the entrance
time of $E_{\partial} - U$.
Let $\alpha \in [0, \infty)$,  let $\varphi$ a $\mathcal{B}$-measurable function on $E$
and $f$ a $\mathcal{B}$-measurable function on~$U$.
Let $\Lambda\left( n\right)$ be a closed subset of $E$ for  $n=1,2, \ldots$.
If $\alpha = 0$ and $\varphi \neq 0$,  assume that $\Lambda\left( n\right)$
is such that $\mathbf{P}_{\pi}\left( \sigma = \infty\right) = 0$. 
Let $\pi$ be a probability measure on $E$
such that $\mathbf{E}_{\pi} \left[  \int_0^{\sigma} e^{ - \alpha t} \left\lvert f\right\rvert \left( X_t\right) \, d t  \right]< \infty$
and
$\mathbf{E}_{ \pi} \left[   e^{ - \alpha \tau_n \wedge \sigma} \left\lvert\varphi\right\rvert \circ X_{\tau_n \wedge \sigma }  \right] < \infty$
for each~$n$.
Let $u_n$ be the probabilistic solution of the Dirichlet problem on $U - \Lambda\left( n\right)$, 
given by~\eqref{general_prob_Dirichlet_eqn} for $\pi$-a.e.~$x$, 
where $\tau_n = D_{\Lambda\left( n\right)}$.
The following lemma gives conditions under which
stable convergence of stopping times implies convergence for the 
corresponding solutions of the Dirichlet problem.
Similar facts were proved in \cite{baxter_chacon_jain} for the Brownian motion case.


\begin{lemma}
\label{convergence_for_probability_solutions_lemma}
Suppose that there exists a randomized stopping time 
$\tau$ such that $\tau$
has a rate measure and $\tau_n \wedge \zeta \rightarrow \tau \wedge \zeta$, $\mathbf{P}_{\pi}$-stably.
Let $u$ be defined by
\begin{equation}
\label{rate_prob_Dirichlet_eqn}
u(x) = \int  \int_0^{\tau \wedge \sigma} e^{ - \alpha t} f\left( X_t\right) \, d t \, d \mathbf{P}_{x} \, d \lambda^{\ast}
+ \int  e^{ - \alpha \tau \wedge \sigma} \varphi \left( X_{\tau \wedge \sigma} \right) \, d \mathbf{P}_{x} \, d \lambda^{\ast}.
\end{equation}
Assume that $\mathbf{E}_{\pi} \left[
 \sup_{0 \leq t < \sigma}  e^{ - \alpha t} \left\lvert\varphi\right\rvert \circ X_{t } \right] < \infty$, and
$\mathbf{P}_{\pi}\left( t \mapsto \varphi \circ X_t \text{ is continuous on } [0, \sigma) \right) =1$.
If $\varphi$ is nonzero and $\zeta$ is not identically equal to $\infty$, 
assume that the cemetery point $\partial$ is not a limit of points in $\bigcup_n \Lambda\left( n\right)$.
Then $u_n$ converges weakly to $u$,
in the sense that $\int u_n g \, d \pi \rightarrow \int u g \, d \pi$
for any $g  \in L^{2}\left(  \pi\right)$.

If the Markov process $X$ and $\pi$
are such that $ \delta_x \, R_{\alpha} << \pi$
for $\pi$-a.e.~$x$ then $u_n  \rightarrow u$ in $\pi$-measure. 
\end{lemma}

\medskip \par \noindent
        \mbox{\bf Proof}\ \   
Let $Y_t
= g \left( X_0 \right) \int_0^{t}  e^{ - \alpha r} f\left( X_r\right) \, d r$, 
$Z_t = g \left( X_0 \right) e^{ - \alpha t } \varphi \left( X_{t \wedge \sigma} \right)$,
so that $\int Y_{\tau_n \wedge \sigma} \, d \mathbf{P}_{\mu} \, d \lambda^{\ast}
+ \int  Z_{\tau_n }  \, d \mathbf{P}_{\mu} \, d \lambda^{\ast}
= \int u_n g\, d \pi$.
By Lemma~\ref{extended_stable_convergence_lemma}(i), 
we have immediately that
$\int Y_{\tau_n \wedge \sigma} \, d \mathbf{P}_{\mu} \, d \lambda^{\ast} \rightarrow \int   Y_{\tau \wedge \sigma} \, d \mathbf{P}_{\mu} \, d \lambda^{\ast}$.
Suppose $\varphi \neq 0$ and $\zeta$ is not identically equal to $\infty$.
Let $\tau^\ast$ be any $\mathbf{P}_{ \pi}$-stable limit point.
Let $F$ be the closure of $\bigcup_n \Lambda\left( n\right)$.
By Lemma~\ref{cannot_converge_to_lifetime_lemma}
with $B_F =  \Omega$,
$\mathbf{P} \times \lambda^{\ast} \left(  \left\{ X_{\tau^\ast} \notin F,  \tau^{\ast}  < \infty\right\}\right)= 0$.
Thus $\mathbf{P} \times \lambda^{\ast} \left(  \zeta \leq \tau^{\ast} < \infty\right) = 0$. 
Also, since $\tau$ has a rate measure, 
$\mathbf{P}_{\pi} \times \lambda^{\ast} \left(  \zeta \leq \tau < \infty\right) = 0$.
Hence by Lemma~\ref{convergence_to_time_made_infinite_lemma}, 
$\tau_n \rightarrow \tau$, $\mathbf{P}_{\pi}$-stably.
It is easy to check
that $t \mapsto Z_t$ has the
continuity required by Lemma~\ref{extended_stable_convergence_lemma}(i)
at all times except for  $t = \sigma < \infty$ when $\alpha >0$,
and at all times except for $t = \sigma $ when $\alpha = 0$.
Since $\tau$ has a rate measure it follows $\mathbf{P}_{\pi}\left( \tau = \sigma < \infty\right) = 0$,
and if $\alpha = 0$
it is assumed that $\mathbf{P}_{\pi}\left( \sigma = \infty\right) = 0$.
Thus Lemma~\ref{extended_stable_convergence_lemma}(i)
applies, so
$\int Z_{\tau_n} \, d \mathbf{P}_{\mu} \, d \lambda^{\ast} \rightarrow \int   Z_{\tau} \, d \mathbf{P}_{\mu} \, d \lambda^{\ast}$.
The same conclusion holds more easily when $\zeta = \infty$.
This proves the weak convergence.

If $ \delta_x \, R_{\alpha} << \pi$,
then by Lemma~\ref{stable_from_a_point_lemma}
there exists a subsequence $n_k$ 
such that $\tau_{n_k} \wedge \zeta \rightarrow \tau \wedge \zeta$, 
$\mathbf{P}_{x}$-stably, for $\pi$-a.e.~$x$.
By what has already been proved,
$u_{n_k}(x) \rightarrow u(x)$ for $\pi$-a.e.~$x$.
Since $n_k$ could be chosen as a subsequence from any other subsequence, 
it follows that
$u_n  \rightarrow u$ in $\pi$-measure.
\hfill \mbox{\raggedright \rule{0.1in}{0.1in}}

The conclusion of Lemma~\ref{convergence_for_probability_solutions_lemma} 
holds for any limit $\tau$
such that $S^{\tau}_{t}$ is a multiplicative functional, 
with the same proof.
We note that the lemma
requires no smoothness or continuity for $\varphi$ on $E - U$.
Also, when $\varphi = 0$ the proof of weak convergence 
is valid for any limit $\tau$.



Lemma~\ref{convergence_for_probability_solutions_lemma}
deals with the convergence of $u_n$ to~$u$ when $u_n$ and~$u$
are defined probabilistically on an open subset $U$ of $E$, 
for a general Markov process~$X$.
Suppose now that~$X$ is properly associated with 
a Dirichlet form $\mathcal{E}$.
 Let $\varphi$
be an $\mathcal{E}$-quasi-continuous function in $D\left(\mathcal{E}\right)$.
The next lemma gives analytical characterizations of~$u_n$ and~$u$ in terms of $\mathcal{E}$,
for $\alpha >0$ and $f \in L^{2}\left( m\right)$.   It follows that
the probabilistic solutions for the Dirichlet problem
 agree with the analytical solutions which were studied
in  \cite{biroli_mosco}, \cite{dalmaso_decicco_notarantonio_tchou}, \cite{biroli_tchou}, 
\cite{mataloni_tchou}.  The case of $\alpha = 0$ can be dealt with similarly in situations
where an inequality of Poincar\'{e} type holds. 

\begin{lemma}
\label{analytical_from_limit}
Let $\tau_n = D_{ \Lambda\left( n\right)}$.
Let $\alpha >0$, $f \in L^{2}\left( m\right)$, 
and let $\varphi$ be an $\mathcal{E}$-quasi-continuous function in $D\left( \mathcal{E}\right)$. 
If $\varphi$ is nonzero, assume for each~$n$ 
that the cemetery point $\partial$ is not  both a limit of points in $U - \Lambda\left( n\right)$ and a limit of points in 
 $\left( E -  U \right) \cup \Lambda\left( n\right)$. Let
$u_n$ be defined by~\eqref{general_prob_Dirichlet_eqn}, for all $x \in E$
such that the expected values exist.
Then $u_n$ is defined $\mathcal{E}$-q.e.

For $n = 1, 2, \ldots$ let $V_n$ denote the set of all $v \in D\left(\mathcal{E}\right)$
such that $v$ is $\mathcal{E}$-quasi-continuous and 
$v =0$ holds $\mathcal{E}$-q.e.\ on $U^c \cup \Lambda\left( n\right)$.
Then $u_n$ is the $\mathcal{E}$-q.e.\ unique element in $ D\left(\mathcal{E}\right)$
such that $u_n$ is $\mathcal{E}$-quasi-continuous, 
$u_n =\varphi$ holds $\mathcal{E}$-q.e.\ on $U^c \cup \Lambda\left( n\right)$ and
\begin{equation}
\label{general_hitting_resolvent_eqn}
\mathcal{E}_{\alpha}\left(u_n,v\right) = \int  f  v \, d m
\end{equation}
for all $v \in V_n$. 
Furthermore
$\left\lVert u_n\right\rVert_{\mathcal{E}, \alpha}$ is bounded in~$n$.

Let $\eta$ be a measure with finite energy, 
and let $\tau$ be the randomized stopping time
with rate measure~$\eta$. 
Let $V$ be the set of $\mathcal{E}$-quasi-continuous 
elements~$v$ in $D\left(\mathcal{E}\right) \cap L^{2}\left(  \eta\right)$
such that $v = 0$ holds $\mathcal{E}$-q.e.\ on $U^c$.  
Assume that $\int_U \varphi^2 \, d \eta < \infty$
and $\varphi \mathbf{1}_{U} \eta$ has finite energy. 
Let $u$ be defined by~\eqref{rate_prob_Dirichlet_eqn},
for all $x \in E$
such that the expected values exist.
Then $u$ is defined $\mathcal{E}$-q.e. and
$u$ is the $\mathcal{E}$-q.e.\ unique element in $D\left(\mathcal{E}\right)$
with $\int_U u^2 \, d \eta < \infty$ such that
$u$ is $\mathcal{E}$-quasi-continuous, $u =\varphi$ holds
$\mathcal{E}$-q.e.\ on $U^c$ and
\begin{equation}
\label{weak_Poisson_rate_boundary_eqn}
\mathcal{E}_{\alpha}\left(u,v\right) + \int u  v\, d \eta
= \int f v \, d m +  \int  \varphi v \, d \eta
\end{equation}
for all $v \in V$.
The space $V$ can be replaced by the space $V^{\prime}$
consisting of all functions in $V$
which vanish 
$\mathcal{E}$-q.e.\ on the complement of a compact subset of~$U$,
or in the regular case by the space $V^{\prime \prime}$ of functions
in  $V \cap \mathcal{C}_{0} \left( E\right)$
with compact support in~$U$.

Now assume that
the absolute continuity condition~\eqref{kernel_density_cond_eqn} holds,
and that for some probability measure $\pi$ with $m << \pi$ on $U$,
$\tau_n \wedge \zeta \rightarrow \tau \wedge \zeta$, $\mathbf{P}_{\pi}$-stably.
If $\varphi$ is nonzero assume that $m(U) < \infty$.
Then $\left\lVert u_n - u\right\rVert_{L^{2}\left(m\right)} \rightarrow 0$
and $u_n \rightarrow u$, $\mathcal{E}$-weakly.
\end{lemma}

\medskip \par \noindent
        \mbox{\bf Proof}\ \   
We may assume $\varphi \geq 0$.
By equation~\eqref{background_harmonic_via_entrance_expectation_eqn}, $\mathbf{E}_{z} \left[ 
e^{ - \alpha   \tau_n \wedge \sigma} \varphi\left( X_{ \tau_n \wedge \sigma} \right) \right]
= {\varphi}^{\prime}_{U^c \cup \Lambda\left( n\right) , \alpha}$ for $\mathcal{E}$-q.e.\ $z$.
Also $R_{\alpha} \,  \left\lvert f\right\rvert$
is $\mathcal{E}$-quasi-continuous, and hence finite $\mathcal{E}$-q.e.

Let $w_n = u_n - {\varphi}^{\prime}_{U^c\cup \Lambda\left( n\right), \alpha}$.
By equation~\eqref{general_prob_Dirichlet_eqn},
$w_n = \mathbf{E}_{x} \left[  \int_0^{\tau_n \wedge \sigma} e^{ - \alpha t} f\left( X_t\right) \, d t \right]$.
Also by equation~\eqref{background_harmonic_via_entrance_expectation_eqn}
 and the strong Markov
property,
$w_n = 
G_{\alpha}\, f - 
{\left(G_{\alpha}\, f\right)}^{\prime}_{U^c \cup\Lambda\left( n\right), \alpha}$.
It follows that $u_n \in D\left(\mathcal{E}\right)$.
By definition  $u_n = \varphi$ on $U^c \cup\Lambda\left( n\right)$,
and by equation~\eqref{background_reduced_bform_agrees_if_q_e_same_on_A_eqn},  $\mathcal{E}_{\alpha}\left( u_n,v\right)
= \mathcal{E}_{\alpha}\left( G_{\alpha}\, f,v\right)
= \int f v \, d m$ for any $v \in V$.
This proves equation~\eqref{general_hitting_resolvent_eqn}.
If $u_n, u_n^\prime$ are solution of equation~\eqref{general_hitting_resolvent_eqn}
for the same~$\varphi$ then $u_n - u_n^\prime \in V$, 
and hence $\left\lVert u_n- u_n^\prime\right\rVert_{\mathcal{E}, \alpha} = 0$ 
by equation~\eqref{general_hitting_resolvent_eqn}, so the solution is unique.

Since $\left\lVert{\left(G_{\alpha}\, f\right)}^{\prime}_{U^c \cup\Lambda\left( n\right), \alpha}\right\rVert_{\mathcal{E}, \alpha}
\leq K_{\alpha} \left\lVert G_{\alpha}\, f\right\rVert_{\mathcal{E}, \alpha}$
and $\left\lVert{\varphi}^{\prime}_{U^c \cup \Lambda\left( n\right), \alpha}\right\rVert_{\mathcal{E}, \alpha}   
\leq K_{\alpha} \left\lVert \varphi\right\rVert_{\mathcal{E}, \alpha}$,
\begin{equation}
\label{uniform_enorm_bound_eqn}
\left\lVert u_n\right\rVert_{\mathcal{E}, \alpha} \leq \left( 1 + K_{ \alpha} \right) \left\lVert G_{\alpha}\, f\right\rVert_{\mathcal{E}, \alpha}
+  K_{\alpha}  \left\lVert \varphi\right\rVert_{\mathcal{E}, \alpha}
\end{equation}
for all~$n$.

In considering~\eqref{rate_prob_Dirichlet_eqn}, 
since replacing $\eta$ by $\mathbf{1}_{U} \eta$ leaves $\tau \wedge \sigma$ unchanged, 
we may assume without loss of generality that $\eta (U^c) = 0$, i.e.
$\eta = \mathbf{1}_{U} \eta$.
Let  $\left( A_t\right)$ be the positive continuous additive functional 
with Revuz measure~$\eta$ and 
let $\Phi^{\tau}$ be the random measure on $[0, \infty]$ such that 
$\Phi^{\tau}{\left( (t, \infty] \right)} = e^{ - A_t} $. 
Then for each $t$, $\int  e^{ - \alpha t} \mathbf{1}_{\left\{  \tau > t  \right\}  }  f \left( X_t \right) \, d \mathbf{P}_{x} \, d \lambda^{\ast}
= \mathbf{E}_{x} \left[  e^{ - \alpha t} e^{ - A_t} f\left( X_t\right) \right]$
by equation~\eqref{random_measure_integral_survival_function_eqn}
with $Y_t = e^{ - \alpha t} f\left( X_t\right)$.
Also  
by taking $Y_t = e^{ - \alpha( t \wedge \sigma)}   \varphi\left( X_{t \wedge \sigma} \right)$ 
in equation~\eqref{random_measure_integral_survival_function_eqn}
we have $\int  e^{ - \alpha \tau \wedge \sigma} \varphi \left( X_{\tau \wedge \sigma}\right)  \, d \mathbf{P}_{x} \, d \lambda^{\ast}
= \mathbf{E}_{x} \left[  \int_0^\sigma e^{ - \alpha t} e^{ - A_t} \varphi\left( X_t\right) \, d A_t
+ e^{ - \alpha \sigma} e^{ - A_\sigma} \varphi \left( X_{\sigma} \right) \right]$.
Then
equation~\eqref{rate_prob_Dirichlet_eqn} implies 
\begin{equation}
\label{PCAF_expression_u_eqn}
u(x) = \mathbf{E}_{x} \left[  \int_0^{\sigma} e^{ - \alpha t - A_t} f\left( X_t\right) \, d t
+ \int_0^{\sigma} e^{ - \alpha t - A_t} \varphi \left( X_t\right) \, d A_t
+ e^{- \alpha \sigma -A_{\sigma}} \varphi \left( X_{\sigma}\right) \right].
\end{equation}

Let $w(x) =  \mathbf{E}_{x} \left[  \int_0^{\sigma} e^{ - \alpha t - A_t} f\left( X_t\right) \, d t \right]$, 
$g(x) =  \mathbf{E}_{x} \left[ \int_0^{\sigma} e^{ - \alpha t - A_t} \varphi \left( X_t\right) \, d A_t \right]$, 
$h(x) = \mathbf{E}_{x} \left[  e^{ - \alpha \sigma - A_{\sigma}} \varphi \left( X_{\sigma}\right) \right]$, 
so that $u  = w + g + h$.  Clearly $u = \varphi$ holds on $\Lambda$.

We will show that each of the functions $w,g,h$ satisfies an equation similar to equation~\eqref{weak_Poisson_rate_boundary_eqn}.
First we will deal with~$g$. 
Since $\varphi \eta$ has finite energy, 
$g \leq G_{\alpha}\,  \left( \varphi \eta\right) $
by equation~\eqref{background_harmonic_via_entrance_expectation_eqn}.
 It is easy to check that if $\xi$ is a bounded function then $\xi \eta$ has
finite energy.
Let $H_j = \left\{\varphi \leq j, \ G_{ \alpha}\, \left( \varphi \eta\right) \leq j \right\}$, for $j = 1,2, \ldots$.
Let $\varphi_j = \mathbf{1}_{ H_j } \varphi $.  By the domination principle, 
$G_{\alpha}\,  \left( \varphi_j \eta \right) \leq j$ holds $\mathcal{E}$-q.e.
Let $g_j = \mathbf{E}_{x} \left[ \int_0^{\sigma} e^{ - \alpha t - A_t} \varphi_j \left( X_t\right) \, d A_t \right]$.
By equation~\eqref{background_additive_func_pot_energy_eqn}
we know that $g_j \leq G_{\alpha}\,  \left( \varphi_j \eta\right)$, 
and hence $g_j (x)\leq j$ for $\mathcal{E}$-q.e.~$x$.  Also $g_j \nearrow g$, $\mathcal{E}$-q.e.,
by III.3.5 of \cite{ma_roeckner}.

Similarly to equation~(4.1.7) in \cite{chen_fu}, one can show that
$\mathbf{E}_{x} \left[  \int_0^{\sigma} e^{ - \alpha t} \varphi_j\left( X_t\right) \, d A_t \right]
= g_j + \mathbf{E}_{x} \left[  \int_0^{\sigma} e^{ - \alpha t} g_j \left( X_t\right) \, d A_t \right]$ for 
$\mathcal{E}$-q.e.~$x$.  

Using equation~\eqref{background_harmonic_via_entrance_expectation_eqn}, 
\eqref{background_additive_func_pot_energy_eqn}
and the strong Markov property,  we then have
$G_{\alpha}\,  \left( \varphi_j \eta\right) -
 {\left( G_{\alpha}\, \left( \varphi_j \eta\right ) \right) }^{\prime}_{\Lambda, \alpha} 
 = g_j + G_{\alpha}\, \left( g_j \eta\right) - 
{ \left( G_{\alpha}\, \left( g_j \eta\right) \right) }^{\prime}_{\Lambda, \alpha}$.
Hence $g_j = G_{\alpha}\,  \left( \varphi_j \eta\right)  - 
{\left( G_{\alpha}\, \left( \varphi_j \eta\right) \right)}^{\prime}_{\Lambda, \alpha} 
-  G_{\alpha}\, \left( g_j \eta\right) +
{ \left( G_{\alpha}\, \left( g_j \eta\right) \right) }^{\prime}_{\Lambda, \alpha}$.

Let $v \in V$.  Then $\mathcal{E}_{\alpha}\left( G_{\alpha}\,  \left( g_j \eta\right) ,v\right) = 
\int g_j v \, d \eta$,
$\mathcal{E}_{\alpha}\left( G_{\alpha}\,  ( \varphi_j \eta), v\right)
= \int \varphi_j v \, d \eta$, and 
$\mathcal{E}_{\alpha}\left({ v_1}^{\prime}_{\Lambda, \alpha},v\right) = 0$ for any $v_1 \in D\left( \mathcal{E}\right)$.
It follows that 
$\mathcal{E}_{\alpha}\left(g_j,v\right)
= 
\int \varphi_j v \, d \eta - \int g_j v \, d \eta$, 
so $\mathcal{E}_{\alpha}\left(g_j,v\right) +  \int g_j v \, d \eta
= 
\int \varphi_j v \, d \eta$.
Since $g_j \in V$, we have $\mathcal{E}_{\alpha}\left(g_j,g_j\right) + \int g_j^2 \, d \eta
\leq \int \varphi_j g_j \, d \eta \leq \varphi g \, d \eta \leq \int \left( G_{ \alpha}\,  \left( \varphi \eta \right)  
\right) \varphi \, d \eta = 
\mathcal{E}_{\alpha}\left( G_{ \alpha}\,  \left( \varphi \eta \right) , G_{ \alpha}\,  \left( \varphi \eta \right) \right) <\infty$ for all $j$.
Letting $j \nearrow \infty$, using this equation and I.2.12 in \cite{ma_roeckner}, 
it follows that $g \in D\left(\mathcal{E}\right) \cap L^{2}\left(  \eta\right)$ and
$\mathcal{E}_{\alpha}\left(g,v\right) +  \int g v \, d \eta
= 
\int \varphi v \, d \eta$.

The same arguments show that 
$w \in D\left(\mathcal{E}\right) \cap L^{2}\left(  \eta\right)$ and
$\mathcal{E}_{\alpha}\left(w,v\right) +  \int w v \, d \eta
= 
\int f v \, d m$, for all $v \in V$.

Let $h_j(x) = \mathbf{E}_{x} \left[  e^{ - \alpha \sigma - A_{\sigma}} \varphi\wedge j \left( X_{\sigma}\right) \right]$.
One can show that
$\mathbf{E}_{x} \left[ e^{ - \alpha \sigma} \varphi \wedge j \left(X_\sigma\right) \right]
= h_j + \mathbf{E}_{x} \left[  \int_0^{\sigma} e^{ - \alpha t} h_j \left( X_t \right) \, d A_t \right]$, 
which says that
${ \left( \varphi \wedge j \right) }^{\prime}_{ \Lambda, \alpha}
= h_j +  G_{\alpha}\,  (h_j \eta) - { \left( G_{\alpha}\,  {  (h_j \eta)} \right) }^{\prime}_{ \Lambda,  \alpha}$.
Thus $\mathcal{E}_{\alpha}\left(h_j , \right){v} + \int h_j v \, d \eta = 0$ for any $v \in V$.
Clearly $\left\lVert \varphi \wedge j\right\rVert_{\mathcal{E}, \alpha} \leq \left\lVert \varphi\right\rVert_{\mathcal{E}, \alpha}$, 
and since $h_j \leq h$ we have $\left\lVert  G_{\alpha }\,  \left( h_j \eta \right)  \right\rVert_{\mathcal{E}, \alpha}
\leq \left\lVert  G_{\alpha }\,  \left( h \eta \right)  \right\rVert_{\mathcal{E}, \alpha}$.
Thus $\left\lVert h_j \right\rVert_{\mathcal{E}, \alpha}$ is bounded in~$j$.  Hence by I.1.12 of \cite{ma_roeckner},
$h \in D\left( \mathcal{E}\right)$ and $\mathcal{E}_{\alpha}\left(h_j,v\right)\rightarrow \mathcal{E}_{\alpha}\left(h,v\right)$.
Since $\varphi \wedge j - h_j \in V$, 
$\mathcal{E}_{\alpha}\left(h_j,\varphi \wedge j - h_j\right) + \int h_j \left( \varphi \wedge j - h_j\right) \, d \eta = 0$,
and hence $\int h_j^2 \, d \eta \leq c_j + \left( \int h_j^2 \, d \eta \right)^{1/2}
\left( \int \varphi^2 \, d \eta \right)^{1/2}$, where $c_j$ is bounded.   It follows that $\int h_j^2 \, d \eta$
is bounded, and hence that $\int h^2 \, d \eta < \infty$.
Hence  $\mathcal{E}_{\alpha}\left(h,v\right)+ \int h v \, d \eta = 0$.

Adding the equations for $w, g, h$ gives equation~\eqref{weak_Poisson_rate_boundary_eqn}.  
If $u, u^\prime$ are solutions of equation~\eqref{weak_Poisson_rate_boundary_eqn}
for the same $\varphi$ then 
for any $v \in V$, $\mathcal{E}_{\alpha}\left( u - u^\prime,v\right) 
+ \int \left( u - u^\prime\right) v \, d \eta = 0$
by equation~\eqref{weak_Poisson_rate_boundary_eqn}.  Since $u - u^\prime \in V$, 
$u = u^\prime$, and the solution is unique.

Using I.2.12 in \cite{ma_roeckner},
it is straightforward to show that $V^{\prime}$
is dense in $V$ with respect to $\left\lVert\cdot\right\rVert_{\mathcal{E}, \alpha}$-norm, 
and that when $\mathcal{E}$ is regular, $V^{\prime \prime}$
is dense.

Since $\varphi$ is $\mathcal{E}$-quasi-continuous, 
$t \mapsto \varphi \circ X_t$ is continuous with $\mathbf{P}_{\pi}$-probability one
on $[0, \sigma)$ for every smooth $\pi$.
Assume $\pi << m$ and
 that $\tau_n \wedge \zeta \rightarrow \tau \wedge \zeta$, $\mathbf{P}_{\pi}$-stably.
If $\varphi$ is nonzero assume that $m(U) < \infty$.
In proving the final statement of the lemma we may assume
using~\eqref{uniform_enorm_bound_eqn} and IV.4.17 of \cite{ma_roeckner}
that $\varphi$ is bounded.
By Lemma~\ref{convergence_for_probability_solutions_lemma}, 
 $u_n \rightarrow u$ in $m$-measure.
By~\eqref{general_prob_Dirichlet_eqn},  $u_n(x) \leq G_{\alpha}\,  \left\lvert f\right\rvert  
+ \text{ a constant}$, so $\left\lVert u_n - u\right\rVert_{L^{2}\left(m\right)} \rightarrow 0$
by dominated convergence.  Hence 
also $u_n \rightarrow u$, $\mathcal{E}$-weakly, 
by I.2.12 in \cite{ma_roeckner}.
\hfill \mbox{\raggedright \rule{0.1in}{0.1in}}

Comparing~\eqref{general_hitting_resolvent_eqn}
with~\eqref{weak_Poisson_rate_boundary_eqn}, we see that the Dirichlet
boundary condition on $\Lambda\left( n\right)$ has been replaced by a penalty term associated
with the measure~$\eta$, together with an additional source term
if $\varphi$ does not vanish on $\Lambda\left( n\right)$.  For any measure $\eta$, not necessarily finite,
which does not charge 
$\mathcal{E}$-exceptional sets, equation~\eqref{weak_Poisson_rate_boundary_eqn}
with $\varphi = 0$ is said to describe  the relaxed Dirichlet problem for~$u$
with zero boundary conditions,  having
 penalty measure~$\eta$.  The function~$u_n$ defined by~\eqref{general_hitting_resolvent_eqn}
with $\varphi = 0$ can also be regarded as the solution of a relaxed Dirichlet problem, 
namely
\begin{equation}
\label{weak_Poisson_rate_boundary_n_eqn}
\mathcal{E}_{\alpha}\left(u_n,v\right) + \int u_n  \tilde{v} \, d \eta_n
= \int f v \, d m,
\end{equation}
where  the measure~$\eta_n$ is infinite
on all subsets of $\Lambda\left( n\right)$ which have positive capacity,
and is zero otherwise.  This is the setting in which
the convergence of $u_n$ to~$u$ has been studied as a special case of 
convergence for solutions of relaxed Dirichlet problems, 
in the analytical papers cited in Section~\ref{random_holes_subsection}.
  General relaxed Dirichlet problems
are not considered in the present paper, although 
a probabilistic representation of the solution of a general relaxed Dirichlet problem
has been given in \cite{getoor}.

\section{Transformations}
\label{transformations_section}

\subsection{Using Girsanov's Theorem}
\label{using_Girsanov_subsection}
In order to prove stable convergence of stopping times,
it may be possible to use a Girsanov transformation (cf. 7.6.4 in \cite{liptser_shiryaev})
to reduce the problem to the case of a simpler process.
Let $X, Y$ be processes which can be defined on the same sample space
$ \Omega$.  Let $\mathbf{P}$ and $\mathbf{Q}$ be probabilities on $ \Omega$
for $X$ and $Y$ respectively.

Let $\mathcal{G}_{t}$ be $\sigma\left( X_s, \, s \leq t\right)$,
$\mathcal{G} = \sigma \left( X_s, \ s \in [0, \infty) \right)$.
Let $\tau_n$ be a sequence of $\left( \mathcal{G}_{t} \right)$-stopping times, 
and let $\tau$ be a randomized stopping time.
Assume that it has been shown, by any method, that $\tau_n  \rightarrow \tau$, $\mathbf{Q}, \mathcal{G}$-stably.
If $\mathbf{P}$ can be obtained by a Girsanov transformation
from $\mathbf{Q}$, then  $\mathbf{P} << \mathbf{Q}$ on $\mathcal{G}_{t}$
for each~$t$.  
In this case, since $\bigcup_{t} \mathcal{G}_{t}$ is $\mathbf{P}$-dense in $\mathcal{G}$, 
it follows that 
$\tau_{n} \rightarrow \tau $, $\mathbf{P}, \mathcal{G}$-stably.

\subsection{Time changes}
\label{time_change_subsection}
A time change can be used to change the measure~$m$
for the $L^2$-space containing  $D\left(\mathcal{E}\right)$. 
\begin{lemma}
\label{continuous_time_change_lemma}
Assume that $\mathcal{E}$
is regular and symmetric.
Let $D_{e}\left(\mathcal{E}\right)$ denote the extended Dirichlet space
associated with $\mathcal{E}$, and let $\mathcal{E}^e$ 
be the extension of $\mathcal{E}$ to $D_{e}\left(\mathcal{E}\right)$, as defined after 
Definition~1.1.4 and Theorem~1.1.5 of \cite{chen_fu}.
Let $b$ be a locally bounded function in $\mathcal{B}^{+}$ with $b >0$ everywhere, 
and define the positive continuous additive functional~$B$
by $B_t = \int_0^t b \left( X_s \right) \, d s$.
Define the partial inverse $\kappa_{t}$ on $[0, \infty]$ by $\kappa_{t} = \inf \left\{s:  B_s > t\right\}$
for $t < B_{\zeta - }$ and $\kappa_{t} = \infty$ for $t \geq B_{\zeta -}$.
Let $\check{X}_t = X_{\kappa_{t} }$.
Then $\check{X}$ is a Markov  process on~$ \Omega$
with lifetime $\check{\zeta} = B_{\zeta  -}$,  which 
is properly associated with the  Dirichlet form 
$\left( \check{\mathcal{E} } , D\left(\check{\mathcal{E} } \right) \right)$, 
where $D\left( \check{\mathcal{E} } \right) = 
D_{e}\left(\mathcal{E}\right) \cap L^2 \left( E, b m\right)$
and  $\check{\mathcal{E} }\left(f,g\right) = \mathcal{E}^e\left(f,g\right)$
for all $f,g$ in $D\left(\check{\mathcal{E} } \right)$.
A subset of $E$ is $\check{\mathcal{E} } $-exceptional if and only 
if it is $\mathcal{E}$-exceptional.

Let $\Lambda\left( n\right)$ be a sequence of Borel sets.  Let
$\tau_n$ be the entrance or hitting time of $\Lambda\left( n\right)$
by~$X$, and correspondingly let $\check{\tau}_n$ be the entrance 
or hitting time 
of $\Lambda\left( n\right)$ by $\check{X}$.  Then there is 
a $\mathcal{E}$-exceptional set $N$ such that 
if $x \notin N$, then with $\mathbf{P}_{x}$-probability one
we have  $\check{\tau}_n = B_{\tau_n}$
if $\tau_n < \zeta$, $\check{\tau}_n = \infty$ otherwise.
Let $\nu$ be a probability measure such that $\nu (N) = 0$.
 If $\tau_n \wedge \zeta$ converges $\mathbf{P}_{\nu}$-stably 
to $\tau \wedge \zeta$, for some
randomized stopping time $\tau$, 
let $\check{\tau} = B_{\tau}$ for $\tau < \zeta$, $\check{\tau} = \infty$
otherwise.
Then $\check{\tau}_n \wedge \check{\zeta} $ converges $\mathbf{P}_{\nu}$-stably 
to $\check{\tau} \wedge \check{\zeta}$. If $S^{\tau}_{t} = e^{ - A_t}$ for
some positive continuous additive functional~$A$,
then $S^{\check{\tau}}_{t}  = e^{ - \check{A}_t}$, 
where $\check{A}$ is the positive continuous additive functional defined by 
$\check{A}_t = A_{\kappa_{t}  - } $.  The 
Revuz measure $\mu_{ \check{A} }$ for $\check{A}$
with respect to $\check{X}$ 
is equal to the Revuz measure $\mu_{A}$ for~$A$
with respect to~$X$.
\end{lemma}

\medskip \par \noindent
        \mbox{\bf Proof}\ \   
Theorem~5.2.2 of \cite{chen_fu} shows that $\check{X}$
is the Markov process associated with $\check{\mathcal{E} } $.
By Theorem~5.2.8 of \cite{chen_fu}, 
a subset of $E$ is $\check{\mathcal{E} } $-exceptional if and only 
if it is $\mathcal{E}$-exceptional.
It follows from the definitions that
$\check{\tau}_n = B_{\tau_n}$ if $\tau_n < \zeta$, $\check{\tau}_n = \infty$
otherwise. 
Let $C_t = B_{t - }$ for $t \in [0, \infty]$.
On a defining set for $B$, $C_t = B_{t}$ for $t < \infty$.
$\check{\tau_n} \wedge \check{ \zeta} 
= C_{ \tau_n \wedge \zeta}$, $\check{\tau} \wedge \check{ \zeta} 
= C_{ \tau \wedge \zeta}$, and $t \mapsto C_t(\omega)$ is continuous
on $[0, \infty]$ for $\omega $ in a defining set for $B$.  
Lemma~\ref{extended_stable_convergence_lemma} (i) 
implies that $\check{\tau}_n \wedge \check{\zeta} $ converges $\mathbf{P}_{\nu}$-stably 
to $\check{\tau} \wedge \check{\zeta}$.  

It is easy to check that 
$S^{\sigma}_{t}(\omega) = 1 - \sup \left( \left\{ r: \ r \in (0,1), \ \sigma(\omega, r) \leq t\right\} \cup \left\{0\right\} \right)$,
for any randomized stopping time $\sigma$.
Applying this formula to $\sigma = \check{\tau}$ shows that
$S^{\check{\tau}}_{t} = e^{ - \check{A}_{t} }$.
By Lemma~6.2.8 of \cite{fu_symm}, 
$\mu_{\check{A}} = \mu_A$.
\hfill \mbox{\raggedright \rule{0.1in}{0.1in}}

Lemma~\ref{continuous_time_change_lemma} is applied in the proofs of 
Lemmas~\ref{translation_invariant_case_lemma} and \ref{general_basic_example_lemma},
to reduce the convergence argument to the case in 
which $m$ is Lebesgue measure $\lambda_d$.

\subsection{Localization}
\label{localization_subsection}
\begin{lemma}
\label{restriction_of_form_lemma}
Let~$X$ be a Hunt process.\
Let~$U$ be an open subset of~$E$ such that
the cemetery point $\partial$ is not in the closure of~$U$.
Let~$m^{U}$ be the restriction 
of~$m$ to the measurable subsets of~$U$.
Let~$\mathcal{V}$ be the set of $f \in D\left(\mathcal{E}\right)$
such that $f = 0$ holds $\mathcal{E}$-q.e.\ on $U^c$.
Let $\mathcal{V}^{U}$
be the set of functions~$h^{U}$ on~$U$, where $h^U$ is
the restriction to~$U$ of a function~$h$ in~$\mathcal{V}$. 
Let~$\mathcal{E}^{U}$
be the form with domain~$\mathcal{V}^{U}$ such that 
$\mathcal{E}^{U}\left(h^{U},g^{U}\right)
= \mathcal{E}\left(h,g\right)$
for all $h,g \in \mathcal{V}$.
Let~$\zeta^{U}$ be the first exit time for~$U$, i.e.\ 
$\zeta^{U} = D_{ E_{\partial} - U}$.  Let
$X^{U}$ be the process defined by $X^{U}_t
= X_{t}$ for $t < \zeta^{U}$, $X_{t}^{U} =\partial$ otherwise, 
with filtration $\left( \mathcal{F}_{t}^{U} \right)$
equal to the closure of the natural filtration 
associated with $X^U$.
Then:
$\left( \mathcal{E}^{U}, D\left( \mathcal{E}^{U}\right) \right)$
is a quasi-regular 
Dirichlet form on $L^{2}\left( U, m^{U} \right)$, and
$X^{U}$ is a special standard process with lifetime 
$\zeta^{U}$,
which is properly associated 
with  
$\left( \mathcal{E}^{U}, D\left( \mathcal{E}^{U}\right) \right)$.
With an appropriate topology on $U \cup \left\{ \partial\right\}$, $X^{U}$
is a Hunt process.
Let $C$ be the set of functions $g^U$, 
where $g \in D\left(\mathcal{E}\right)$ is such 
that $g = 0$ holds $\mathcal{E}$-q.e. on the complement of a
compact subset of $U$.  Then $C$ is dense in $D\left( \mathcal{E}^{U}\right)$.
Any $\mathcal{E}$-exceptional subset of $U$
is $\mathcal{E}^{U}$-exceptional.
If $G_{\alpha}^{U}$ denotes the potential operator for $\mathcal{E}^{U}$, then 
for any $f \in L^{2}\left( m\right)$ such that $f =0$ on $U^c$, 
\begin{equation}
\label{potential_of_restriction_equation}
G_{\alpha}^{U}\, f^{U}
= G_{\alpha}\, f - {\left(G_{\alpha}\, f\right)}_{U^c, \alpha},
\end{equation} 
and for any measure~$\nu$ on $E$ with finite energy, such that 
$\nu = 0$ on all subsets of $U^c$, if $\nu^{U}$ denotes the restriction 
of $\nu$ to subsets of $U$ then $\nu^{U}$ has finite $ \mathcal{E}^{U}$-energy and
\begin{equation}
\label{potential_of_restriction_measure_equation}
G_{\alpha}^{U}\, \nu^{U}
= G_{\alpha}\, \nu - {\left(G_{\alpha}\, \nu\right)}_{U^c, \alpha},
\end{equation} 
Furthermore, if $\left( A_t\right)_{ t \geq 0}$ 
is a positive continuous additive functional for~$X$ with Revuz measure~$\mu$, 
then $\left( A_{ t \wedge \zeta^{U}}\right)$
is the positive continuous additive functional for~$X^U$ with Revuz measure~$\mu^{U}$, 
where $\mu^{U}$ is the restriction of~$\mu$ to subsets of~$U$.
 \end{lemma}

 The proof is omitted.  It follows from the definitions
and the properties given in \cite{ma_roeckner} and 
Section~\ref{Dirichlet_form_properties_section}.  
We will refer to the Dirichlet form $\left( \mathcal{E}^{U}, D\left( \mathcal{E}^{U}\right) \right)$
described in this lemma
as the restriction of $\left( \mathcal{E}, D\left(\mathcal{E}\right)\right)$
to~$U$.
The process~$X^{U}$ has the same sample space as $X$, with an appropriate change in 
the shift operator.

\paragraph*{Combining localizations.}
 When stable convergence can be proved for
a large enough class of  localized versions of a process, global convergence 
 can be obtained by 
 combining restrictions, as in the following lemma.

\begin{lemma}
\label{localization_lemma}
Let~$X$ be a Hunt process on a separable metric space
such that for $\mathcal{E}$-q.e.~$x$,
$t \mapsto X_t$ is continuous on $[0, \zeta)$
with $\mathbf{P}_{x}$-probability one.
Let $U_\ell$, $\ell = 1,2, \ldots$
be a locally finite open cover for~$E$.
Let $\zeta^{\ell}$ be the first exit time for $U_\ell$, 
i.e. $\zeta^{\ell} = D_{ E_{\partial} - U_\ell}$.
Let $X^{U_\ell }_t
= X_{t}$ for $t < \zeta^{\ell}$, $X_{t}^{U_\ell} =\partial$ otherwise.
Suppose that the absolute continuity condition~\eqref{kernel_density_cond_eqn}
holds.

Let $\tau_n$ be a terminal time
for each~$n$.  
Let~$\eta$ be a smooth  measure on~$E$.
For each~$\ell$, let $\eta^{\ell}$
be the restriction of $\eta$ to $U_\ell$.
Suppose that for each~$\ell$ and any smooth probability measure $\pi^\ell$ on $U^{\ell}$, 
$\tau_{n} \wedge \zeta^{\ell}$ converges 
$\mathbf{P}_{\pi^\ell }$-stably
to $\tau^{\ell} \wedge \zeta^{\ell}$, where 
$\tau^{\ell}$ is the randomized stopping time which has
rate measure~$\eta^{\ell}$ for $X^{U_\ell}$.

Let~$\tau$ be the randomized stopping time with rate measure~$\eta$.
Then for any smooth probability measure $\pi$ on $E$, 
$\tau_n \wedge \zeta \rightarrow \tau \wedge \zeta$, $\mathbf{P}_{\pi}$-stably.
 \end{lemma}

The proof is omitted.  The idea of the proof is the following.
By piecing together the stopping times $\zeta^{\ell}$ one obtains a sequence of stopping times 
$\sigma_k$ such that $\sigma_k \nearrow \zeta$.  Using
Lemmas~\ref{extended_stable_convergence_lemma}
and~\ref{stable_from_a_point_lemma} one can show by induction
that 
$\tau_n \wedge \sigma_k \rightarrow \tau \wedge \sigma_k$, $\mathbf{P}_{\pi}$-stably
as $n \rightarrow \infty$, 
and Lemma~\ref{more_extended_stable_convergence_lemma} then gives the result.

\section{Examples}
\label{example_lemmas_section}
Details for proofs of most of the statements in this section, 
and other examples, are given in \cite{baxter_nielsen_supp}.

Let $\mathcal{E}$ be regular.
Let $\left( \tilde{\Omega}_1, \tilde{\mathbf{P}}_1\right)$ be a probability space.
For each $x \in E$ and $n = 1, 2, \ldots$, 
let $\Gamma^{x}_{n}$ be a map from $\tilde{\Omega}_1$ to the collection of 
closed subsets of $E$. Let $F(E)$ be the space of compact
subsets of $E$ with Hausdorff metric. 
It is assumed that for
any $K \in F(E)$, the map 
$\left(x, \tilde{\omega}_1\right) \mapsto \Gamma^{x}_{n} \left( \tilde{\omega}_1\right) \cap K$
is jointly measurable from $E \times \tilde{\Omega}_1$ into $F(E)$.
The case that $\tilde{\Omega}_1$ is a one-point space, so that $\Gamma^{x}_{n}$ is nonrandom,
is an important special case.

Let $\mu$ be a probability 
measure on $E \cup \left\{ \partial\right\}$, where $\partial$
is the cemetery point.  For each~$n$, let
$\xi_1 (n, \mu), \ldots, \xi_{\kappa_{n}}(n, \mu)$ be iid random variables with distribution $\mu$,
defined on some probability space $\left( \tilde{\Omega}_2, \tilde{\mathbf{P}}_\mu \right)$,
where $\lim_{n \rightarrow \infty} \kappa_{n} = \infty$.
Let $\left( \tilde{\Omega}, \tilde{\mathbf{P}}\right) = \left(  \Omega_1 \times  \Omega_2, \tilde{\mathbf{P}}_1 \times \tilde{\mathbf{P}}_\mu\right)$, 
and let $\Lambda^{\mu}_{i}\left( n\right)$ be the random set $\Gamma^{ \xi_i(n, \mu) }_{n}$, 
where we define $\Gamma^{\partial}_{n} = \emptyset$.
The family $\Lambda^{\mu}_{i}\left( n\right)$ will be said to be the
random center model associated with $\left(\Gamma^{x}_{n} \right), \mu$.
The terminology is intended to suggest that $\Lambda^{\mu}_{i}\left( n\right)$
could be randomly chosen by first selecting the random ``center'' $x = \xi_i (n, \mu)$
and then choosing a possibly random set $\Gamma^{x}_{n}$ near $x$.
Let $\Lambda^{\mu}\left( n\right) = \bigcup_i \Lambda^{\mu}_{i}\left( n\right)$.
We are interested in random center models such that
$D_{ \Lambda^{\mu}\left( n\right) }$ and $ T_{ \Lambda^{ \mu}\left( n\right) }$ converge $\mathbf{P}_{\pi}$-stably
in $\tilde{\mathbf{P}}$-probability.  Random center models were studied
in  \cite{kac}, \cite{rauch_taylor}, \cite{papanicolaou_varadhan}, 
\cite{balzano}, \cite{balzano_notarantonio}.

Let $B_r(x)$ denote the open metric ball in $E$ with center $x$ 
and radius $r$.
It is assumed that the sets $\Gamma^{x}_{n}$ become small, 
meaning that there exists a nonrandom sequence $\varrho_n \in (0, \infty)$
such that $\varrho_n \rightarrow 0$ and such that for $\mu$-a.e.\ $x \in E$,
$\tilde{\mathbf{P}}_1\left( \Gamma^{x}_{n}  \subset B_{\varrho_n}(x) \right)  = 1$.

\begin{lemma}
\label{approximate_mu_in_tv_preserve_stable_convergence_lemma}
Let $\mu$ be a probability measure on $E \cup \left\{ \partial\right\}$.
Suppose that for each $n$ there is a constant $\chi_{n} \in [0, \infty)$ such 
that 
\begin{equation}
\label{uniform_capacity_bound_eqn}
\text{$\sup_n \kappa_{n} \chi_{n} < \infty$ and  }\tilde{\mathbf{P}}_1 \left( 
\mathsf{Cap}_{\alpha}\left(\Gamma^{x}_{n} \right) \leq \chi_{n} \right) = 1  \text{ for $\mu$-a.e.\ $x$}.
\end{equation}
Let $B_\ell \in \mathcal{B}$ be nondecreasing and such that $\mu \left( E - B_\ell \right) 
\rightarrow 0$. Let $\mu_\ell$ be the probability measure on $E \cup \left\{ \partial\right\}$
such that $\mu_\ell \left( K \right) = \mu \left( K \cap B_\ell\right)$ for 
every $K \in \mathcal{B}$.  Let $\tau_n = D_{\Lambda^{\mu}\left( n\right)}$ or $\tau_n =  T_{\Lambda^{\mu}\left( n\right)}$,
and let $\tau_n^\ell = D_{ \Lambda^{\mu_\ell}\left( n\right) }$ if $\tau_n = D_{\Lambda^{\mu}\left( n\right) }$,
$\tau_n^\ell =  T_{ \Lambda^{\mu_\ell}\left( n\right) }$ if $\tau_n =  T_{\Lambda^{\mu}\left( n\right) }$.
Let $\pi$ be a smooth probability measure on $E$, such 
that for each~$\ell$, $\tau_n^\ell$
converges $\mathbf{P}_{\pi}$-stably in $\tilde{\mathbf{P}}$-probability
to a randomized stopping time $\tau^\ell$.
Then $\tau^\ell$  decreases to a limit $\tau$, $\mathbf{P}_{\pi}$-a.e., 
and $\tau_n$
converges $\mathbf{P}_{\pi}$-stably in $\tilde{\mathbf{P}}$-probability
to $\tau$ as $n \rightarrow \infty$.
If $\eta^\ell$ is a rate measure for $\tau^\ell$, and $\eta$ is a smooth 
measure such that $\eta^\ell \nearrow \eta$, 
then $\eta$ is a rate measure for $\tau$.
\end{lemma}

The proof is straightforward, and is given in  \cite{baxter_nielsen_supp}.
Condition~\eqref{uniform_capacity_bound_eqn} gives a uniform bound on the total capacity
of the sets $\Lambda_{j}\left( n\right)$.
When $\mu$ is a smooth measure such that~\eqref{uniform_capacity_bound_eqn} holds, 
and we wish to prove convergence for $D_{ \Lambda^{\mu}\left( n\right) }$ or $ T_{ \Lambda^{\mu}\left( n\right) }$,
Lemma~\ref{approximate_mu_in_tv_preserve_stable_convergence_lemma}
allows us to assume that the measure $\mu$ has finite energy
and compact support in $E$. 

We now consider particular random center models for which $E$ is an open subset $U$ of 
${\mathbb R}^{d}$,  $d \geq 2$.  Let $b_{i j}$, $i,j =1, \ldots, d$, be bounded measurable
functions on~$U$, $d \geq 2$,  with $b_{i j } = b_{ j i}$, such that for some 
constant $e_0 >0$, $\sum_{i j} w_i b_{i j}(x) w^j \geq e_0
\sum_i w_i^2$ for all 
$x \in U$, $ w \in {\mathbb R}^{d}$.  We will denote the matrix function $\left( b_{i j} \right)$
by $b$.
Let $\sigma$ be a positive function on $U$ which is 
measurable, bounded and bounded away from zero.
Let 
$\mathcal{E}^{b, \sigma, U}$ be the Dirichlet form on $L^{2}\left( U, \sigma \lambda_d\right)$ such that
$\mathcal{E}^{b, \sigma, U}\left(f,g\right)  = 
\int_{U} \sum_{i j} b_{ i j} \left( \partial_i f\right) \left( \partial_j g\right) \, d \lambda_d$
for any smooth functions $f,g$ on $U$ which have compact support in $U$, 
where $\lambda_d$ denotes Lebesgue measure on ${\mathbb R}^{d}$
and $D\left(\mathcal{E}^{b, \sigma, U}\right)$ is the closure of the space of 
such functions.
$\mathcal{E}^{b, \sigma, U}$ exists 
by II.2 in \cite{ma_roeckner}.  Let 
$X = X^{b, \sigma, U}$ be the Markov process with lifetime $\zeta$
and cemetery point $\partial$
which is properly associated with $\mathcal{E}^{b, \sigma, U}$.

When studying $X$, it will be convenient to consider to extend $b_{i j}$
and $\sigma$ to all of ${\mathbb R}^{d}$, in such a way that $b_{ij}$
is bounded on ${\mathbb R}^{d}$ and  $\sum_{i j} w_i b_{i j}(x) w^j \geq e_0
\sum_i w_i^2$ for all 
$x,  w \in {\mathbb R}^{d}$
and $\sigma$ is bounded and bounded away from zero on ${\mathbb R}^{d}$.
We can define the Dirichlet form $\mathcal{E}^{b, \sigma,  V}$
for any open subset $V$ of ${\mathbb R}^{d}$ analogously to $\mathcal{E}^{b, \sigma, U}$. 
Let $X^{b, \sigma, V}$ be the process associated
with $\mathcal{E}^{b, \sigma, V}$.
By Lemma~\ref{restriction_of_form_lemma}, we may assume
that the process $X = X^{b, \sigma, U}$ is the restriction of $X^{b, \sigma, {\mathbb R}^{d}}$
to $U$, so that $\zeta$ is the exit time of $U$ by $X^{b, \sigma,  {\mathbb R}^{d} } $,
and $X_t = X^{b, \sigma, {\mathbb R}^{d}}_{t } $ when $t < \zeta$, 
$X_t = \partial$ if $t \geq \zeta$.
We will use this version of $X$ in what follows.

It will be assumed from now on that $\Gamma^{x}_{n}$ 
is compact with $\tilde{\mathbf{P}}_1$-probability one for $\mu$-a.e.\ $x$.  If $\Gamma^{x}_{n}$
is not compact, redefine $\Gamma^{x}_{n} = \emptyset$,
and also extend $\Gamma^{x}_{n}$ to all $x \in {\mathbb R}^{d}$
by setting $\Gamma^{x}_{n} = \emptyset$ for $x \in U^c$.
For any open subset $V$ of ${\mathbb R}^{d}$ and any $x \in V$, let $\psi_{n, \alpha}^{ x, b, \sigma, V}$ denote
the $\alpha$-equilibrium measure for $\Gamma^{x}_{n}$
using $\mathcal{E}^{b, \sigma, V}$.
The map $\left(x, \tilde{\omega}_1\right) \mapsto 
\int f \, d \psi_{n, \alpha}^{ x, b, \sigma, V}$ 
is jointly measurable on $V \times \tilde{\Omega}_1$
for any $f \in b\mathcal{B}$, by regularity.
Similarly the map $\left(x, \tilde{\omega}_1\right) \mapsto 
\int \left( G_{\alpha}\,  { \psi_{n, \alpha}^{x, b, \sigma, V} }\right)  \, d \psi_{n, \alpha}^{ x, b, \sigma, V}$
is jointly measurable.
Define the average measure $\bar{\psi}^{x, b, \sigma, V}_{n, \alpha}$
by  $\bar{ \psi }_{n, \alpha}^{x, b, \sigma, V}(W)
= \tilde{\mathbf{E}}_1\left[ \psi^{x, b, \sigma, V}_{n, \alpha} (W)\right]$.
Using the notation of Theorem~\ref{random_sets_example_theorem}, 
for $X = X^{ b, \sigma, U}$, 
 we have $\gamma^{n}_{j} = \psi_{n, \alpha}^{ \xi_j(n, \mu)}$ and
\begin{equation}
\label{ameas_as_integral_eqn}
\bar{\gamma}^{n} = \kappa_{n} \int \bar{\psi}_{n, \alpha}^{x, b, \sigma, U} \, \mu ( d x).
\end{equation}

We first consider the \emph{translation-invariant} case. 
By translation-invariance we mean here that $b$ is constant on $U$,
and the distribution of the sets $\Gamma^{x}_{n} - x$ is the same
for $\mu$-a.e.\ $x$. 
Translation-invariant models in the Brownian motion setting were considered
in \cite{kac}, \cite{rauch_taylor}, \cite{papanicolaou_varadhan}, 
\cite{balzano}, with the sets
 $\Gamma^{x}_{n}$ equal to nonrandom scaled translates of a fixed compact set.
The next lemma differs from earlier results in some technical 
aspects, since $\mu$ is only required to be smooth, the sets $\Gamma^{x}_{n}$
are allowed to be random with different shapes for each~$n$, 
and the rate at which the sets shrink is only constrained by~\eqref{uniform_capacity_bound_eqn}.  
It seems of interest 
as an example for Corollary~\ref{iid_corollary_random_sets_example_theorem} because
of the simplicity of the proof.

\begin{lemma}
\label{translation_invariant_case_lemma}
Let $\mu$ be a smooth probability measure on $U$.
Suppose that  the random center model for $U$ associated
with $\left( \Gamma^{x}_{n} \right), \mu$ is translation-invariant.
Then $\left\lVert \bar{\gamma}^{n}\right\rVert_{\text{tv}} = \kappa_{n} \left\lVert \bar{ \psi}_{n, \alpha}^{x, b, 1, {\mathbb R}^{d}}\right\rVert_{\text{tv}}$
for $\mu$-a.e.\ $x$.  Suppose that~\eqref{uniform_capacity_bound_eqn} holds and
$\lim_{n \rightarrow \infty} \left\lVert \bar{\gamma}^{n}\right\rVert_{\text{tv}} = c \in [0, \infty)$.
Let $\tau_n = D_{ \Lambda\left( n\right)}$ or $\tau_n =  T_{ \Lambda\left( n\right)}$, using the 
process $X = X^{b, \sigma, U}$.
For any smooth probability measure $\pi$ on $U$, 
$\tau_n \rightarrow \tau$, $\mathbf{P}_{\pi}$-stably in $\tilde{\mathbf{P}}$-probability,
where $\tau$ is the randomized stopping time with rate measure $c \mu$.
\end{lemma}

\medskip \par \noindent
        \mbox{\bf Proof}\ \   
By the definition of translation-invariance,
for each $n$ there is a measure $\pi_n$ such that $\bar{\psi}_{n, \alpha}^{b,1, {\mathbb R}^{d}}(W)
= \pi_n (W - x)$ for $\mu$-a.e.~$x$.  
By~\eqref{ameas_as_integral_eqn},
$\kappa_{n} \left\lVert \pi_n\right\rVert_{\text{tv}} = \left\lVert \bar{\gamma}^{n}\right\rVert_{\text{tv}} \rightarrow c$.

By Lemma~\ref{localization_subsection} 
we may assume without loss of generality that $U = {\mathbb R}^{d}$.
By Lemma~\ref{continuous_time_change_lemma} we may then assume that $\sigma = 1$.
Then $X = X^{ b, 1, {\mathbb R}^{d}}$ is essentially Brownian motion,
and  
the absolute continuity condition~\eqref{kernel_density_cond_eqn}
holds. 

By Lemma~\ref{approximate_mu_in_tv_preserve_stable_convergence_lemma}, 
we can assume that $\mu$
has  finite energy and has compact support, so that with $\tilde{\mathbf{P}}$-probability one
all the sets $\Lambda\left( n\right)$ are contained in a compact set.
Since $\varrho_n \rightarrow 0$, 
$\kappa_{n} \pi_n \rightarrow c \delta_0$ weakly as a sequence of measures.

Since
$\bar{\gamma}^{n} = \kappa_{n} \pi_n \ast \mu$, 
$\bar{\gamma}^{n} \rightarrow c \mu$ weakly as a sequence of measures.
Let $\mu^x$ be the translated measure defined by $\mu^x(B) = \mu( B - x)$.
We have 
\begin{multline*}
\kappa_{n}^2 
\mathcal{E}_{\alpha}\left( G_{\alpha}\,  \pi_n \ast \mu  , G_{\alpha}\,  \pi_n  \ast \mu  \right)
= \kappa_{n}^2\int \left(  G_{\alpha}\,  \left(\pi_n \ast \mu \right) \right)
\, d \left( \mu \ast \pi_n\right)
= \kappa_{n}^2 \int \left( G_{\alpha}\,  \left(\pi_n \ast \mu \right) \right)
\, d \mu^x \, \pi_n ( d x)\\
=  \kappa_{n}^2 \int \left( G_{\alpha}\,  \mu^x \right) \, 
d \left( \pi_n \ast \mu\right) 
\, \pi_n ( d x) 
=   \kappa_{n}^2 \int \left( G_{\alpha}\,  \mu^x \right) \, d \mu^y 
\, \pi_n( d y) 
\, \pi_n ( d x)\\
\leq  \kappa_{n}^2
 \int  \left\lVert G_{\alpha}\,  \mu \right\rVert_{\mathcal{E}, \alpha}^2 
\, \pi_n ( d y) 
\, \pi_n ( d x)
= \kappa_{n}^2 \left\lVert G_{\alpha}\,  \mu \right\rVert_{\mathcal{E}, \alpha}^2 \left\lVert \pi_n\right\rVert_{\text{tv}}^2.
\end{multline*}
Hence $\limsup_n \kappa_{n}^2 
\mathcal{E}_{\alpha}\left( G_{\alpha}\,  \left( \pi_n \ast \mu \right) , G_{ \alpha}\,  \left( \pi_n \ast \mu \right) \right)
\leq  c^2\left\lVert G_{\alpha}\,  \mu \right\rVert_{\mathcal{E}, \alpha}^2
= \mathcal{E}_{\alpha}\left( G_{\alpha}\, \left( c\mu \right)  , G_{ \alpha}\, \left(c \mu \right) \right)$.
Since $\left\lVert\kappa_{n} G_{\alpha}\,  {\left(  \pi_n \ast \mu\right)}\right\rVert_{\mathcal{E}, \alpha}$ is bounded,
$\kappa_{n}G_{\alpha}\, \left( \pi_n \ast \mu \right)  \rightarrow 
G_{\alpha}\,  \left( c\mu \right)$, $\mathcal{E}$-weakly,
 by Lemma~\ref{background_vague_gives_E_weak_lemma}.
Hence  also $\liminf_n \kappa_{n}^2
\mathcal{E}_{\alpha}\left( G_{\alpha}\, \left(  \pi_n \ast \mu \right) , G_{ \alpha}\,  \left( \pi_n \ast \mu \right)  \right)
\geq \mathcal{E}_{\alpha}\left( G_{\alpha}\, \left( c\mu \right)  , G_{ \alpha}\, \left(c \mu \right) \right)$,
so
$\kappa_{n} G_{\alpha}\,  \left( \pi_n \ast \mu \right) \rightarrow 
G_{\alpha}\,  \left( c \mu \right) $ in energy norm.
Thus Corollary~\ref{iid_corollary_random_sets_example_theorem} applies
and gives convergence.

\hfill \mbox{\raggedright \rule{0.1in}{0.1in}}

From a physical standpoint, one can think of the holes $\Lambda_{i}\left( n\right)$ 
is representing fixed defects in some material, but one might also
consider the case of \emph{moving obstacles} in a fluid 
medium (``dust particles'' in \cite{kac}). 
These moving holes would presumably travel slowly in 
comparison to Brownian motion, but if their movement 
is considered it would at least affect the formula
for the limit of the stopping times $\tau_n$.
It seems to be an interesting problem to actually prove convergence of $\tau_n$ 
in the Brownian motion case when the holes are moving.
When all holes move with identical velocity $v(t)$, 
Girsanov's theorem can 
be applied to show that if convergence holds without the movement of the holes, 
then convergence holds for the moving case also. 
 It is more reasonable physically
to consider the case in which the holes
$\Lambda_{i}\left( n\right)$, $i = 1, \ldots, \kappa_{n}$,  move independently.
In this case it is plausible that the $\tau_n$ would still
converge in probability under suitable conditions.  We have no 
result of this sort, however.

One can measure the asymptotic capacity of the sets $\Gamma^{x}_{n}$
in various ways.  For any sequences $t_n, u_n$ of nonzero numbers, let
$t_n \sim u_n$ mean that $\lim_{n \rightarrow \infty} t_n/u_n = 1$.
It is not hard to show using the estimates in \cite{stroock}
that 
for any open subset $V$ of ${\mathbb R}^{d}$ with $x \in V$
and any $\alpha, \beta \in (0, \infty)$,
\begin{equation}
\label{point_set_capacity_sigma_V_B_eqn}
\left\lVert \psi_{n, \alpha}^{x, b, \sigma, V }  \right\rVert_{\text{tv}}
\sim \left\lVert \psi_{n, \beta}^{x, b, 1, {\mathbb R}^{d}}  \right\rVert_{\text{tv}},
\end{equation}
with a corresponding asymptotic equivalence for the average measures
$\left\lVert \bar{\psi}_{n, \alpha}^{x, b, \sigma, V }  \right\rVert_{\text{tv}}$ and 
$\left\lVert \bar{\psi}_{n, \beta}^{x, b, 1, {\mathbb R}^{d}}  \right\rVert_{\text{tv}}$.
Equation~\eqref{point_set_capacity_sigma_V_B_eqn}
 holds whether or not the model is translation-invariant.
In particular it shows that the constant $c$ in
Lemma~\ref{translation_invariant_case_lemma} can be expressed in terms of
$\alpha$-equilibrium measures with respect to $\mathcal{E}^{b, \sigma, U}$
if desired, although these measures may not be as easy to compute.

The same arguments used to show~\eqref{point_set_capacity_sigma_V_B_eqn}
also show that for general coefficients $b^{i j}$ satisfying the stated assumptions, 
 $\alpha$-capacity with respect to 
$\mathcal{E}^{b, \sigma, V}$ is locally comparable to classical capacity.
That is, 
given any compact subset $K$ of $V$ there exists a constant $c^\prime$ 
such that when $\varrho_n < 1$, for all $x$ such that $B_{\varrho_n}(x) \subset K$
and all $n$,
\begin{equation}
\label{domination_by_classical_kernel_eqn}
\left\lVert \psi^{x, b, \sigma, V}_{n, \alpha}\right\rVert_{\text{tv}}
\leq c^\prime \varrho_n^{d-2} \text{ if } d > 2, 
\ 
\left\lVert \psi^{x, b, \sigma, V}_{n, \alpha}\right\rVert_{\text{tv}}
\leq c^\prime  / \left( - \log \varrho_n \right) \text{ if } d = 2.
\end{equation}
When $V = {\mathbb R}^{d}$, \eqref{domination_by_classical_kernel_eqn}
holds for all $x \in {\mathbb R}^{d}$.
If the coefficients $b^{i j}$ happen to be continuous on~$U$,
using~\eqref{point_set_capacity_sigma_V_B_eqn} and
equation~\eqref{alpha_capacity_in_as_an_infimum_eqn}
one also finds easily that for $x \in V$,
\begin{equation}
\label{N_x_limit_using_frozen_coefficients_eqn}
\left\lVert \psi_{n, \alpha}^{x, b, \sigma, V }  \right\rVert_{\text{tv}} \sim
\left\lVert \psi^{x, b^x, 1, {\mathbb R}^{d} }_{n,\alpha}\right\rVert_{\text{tv}}, 
\end{equation}
where $b^x$ is the \emph{constant} matrix function equal to $b(x)$ everywhere.
Following an idea in~\cite{balzano}
and~\cite{balzano_notarantonio}
 one can then relate $ \left\lVert \psi_{n, \alpha}^{x, b, \sigma, U }  \right\rVert_{\text{tv}} $ to 
the classical capacity of the sets $\Gamma^{x}_{n}$.
For a compact subset $K$ of ${\mathbb R}^{d}$, 
let $\mathsf{Q}^{\mathsf{cl}} _{d}\left(K\right)$
be the classical capacity of $K$, where the classical capacity is 
calculated using the potential kernel $\wp_{d}\left(y, z\right)  
= 1/\left( \left\lvert z -y\right\rvert^{ d  -2} ( d - 2) \omega_d\right) $
if $d >2$, $\wp_{d}\left(y, z\right)   = - \log \left\lvert z - y\right\rvert/ \omega_2$
if $d = 2$, 
and 
$\omega_d$ here denotes the surface area of the unit hypersphere in ${\mathbb R}^{d}$,
so that for example $\omega_3 = 4 \pi$.
Then
\begin{equation}
\label{N_x_in_terms_classical_capacity}
\left\lVert \psi_{n, \alpha}^{x, b, \sigma, U }  \right\rVert_{\text{tv}}
\sim 
\left(\sqrt{\det b(x)}\right) \mathsf{Q}^{\mathsf{cl}} _{d}\left(b(x)^{-1/2}  \Gamma^{x}_{n} \right),
\end{equation}
where
$b(x)^{-1/2}$ is the inverse of the positive square root of the matrix $b(x)$
and $b(x)^{-1/2} \Gamma^{x}_{n} $ denotes
the set of points $b(x)^{-1/2} z$, 
$z \in \Gamma^{x}_{n}$.
Also, when $U$ is bounded, one can show that
\begin{equation}
\label{N_x_in_terms_bare_capacity}
\left\lVert  \psi_{n, \alpha}^{x, b, \sigma, U }  \right\rVert_{\text{tv}}
\sim
\left(\sqrt{\det b(x)}\right) \mathsf{Q}^{ U}\left(  b(x)^{-1/2}  \Gamma^{x}_{n} \right),
\end{equation}
where we define  
$\mathsf{Q}^{U}\left(K\right) = \inf \left\{ \mathcal{E}^{I,  1,  U}\left(f,f\right) : \ f \in D\left( \mathcal{E}^{I, 1,  U} \right), 
\ f \geq 1 \text{ q.e. on } K\right\}$,
for any compact subset $K$ of $U$, and $I$ is the $d \times d$ identity matrix.

Given Lemma~\ref{translation_invariant_case_lemma},
one would naturally hope that convergence holds for a more general 
case of the random center model for subsets of ${\mathbb R}^{d}$.  However, 
the easy proof of Lemma~\ref{translation_invariant_case_lemma} 
used the translation-invariance of $\bar{\psi}_{n, \alpha}^{x}$ heavily.
A similar proof, using the analog of translation-invariance, is applicable when $\mathcal{E}$
is associated with the Laplace-Beltrami operator on a homogeneous Riemannian manifold.
In the general case a proof can be given by strengthening
the bound in~\eqref{uniform_capacity_bound_eqn}, 
as in~\eqref{size_random_sets_bound_eqn} below.

\begin{lemma}
\label{general_basic_example_lemma}
Let $\mu$ be a smooth probability measure on $U \cup \left\{ \partial\right\}$.
In the random center model for $U$ associated
with $\left( \Gamma^{x}_{n} \right)$ and $\mu$, 
let $q_n (x) = \left\lVert \bar{\psi}^{x, b, 1, {\mathbb R}^{d}}_{n, \alpha} \right\rVert_{\text{tv}}$, 
for $\mu$-a.e.\ $x$, 
and let $\nu_n = q_n \mu$.
Assume that there exists $\varrho_n \in [0, \infty)$ with
$\varrho_n \searrow 0$, such that 
$\tilde{\mathbf{P}}_1 \left( \Gamma^{x}_{n} \subset B_{\varrho_n}(x) \right) =1$
for $\mu$-a.e.\ $x$,  and 
\begin{equation}
\label{size_random_sets_bound_eqn}
\sup_n \kappa_{n} \varrho_n^{ d - 2}
<\infty \text{ if  $d > 2$},\
\sup_n \kappa_{n} /\left\lvert \log \varrho_n\right\rvert < \infty
\text{ if  $d = 2$}.
\end{equation}
Let $\tau_n = D_{ \Lambda\left( n\right)}$ or $\tau_n =  T_{ \Lambda\left( n\right)}$, using the 
process $X = X^{b, \sigma, U}$.
Assume that $\nu_n$ converges weakly as a sequence of measures
to a finite measure $\eta$.
Then for any smooth probability measure $\pi$ on $U$, 
$\tau_n \rightarrow \tau$, $\mathbf{P}_{\pi}$-stably in $\tilde{\mathbf{P}}$-probability,
where $\tau$ is the randomized stopping time with rate measure $\eta$.
The same conclusion holds using $q_n(x) = \left\lVert \bar{\psi}^{x, b, \sigma, U}_{n, \alpha} \right\rVert_{\text{tv}}$.
\end{lemma}
Girsanov's theorem can be used to extend Lemma~\ref{general_basic_example_lemma}
to examples with drift.

Assumption~\eqref{size_random_sets_bound_eqn} is  a uniform smallness condition 
on the sets $\Gamma^{x}_{n}$, and is satisfied by the iid models in 
 \cite{kac}, \cite{rauch_taylor}, \cite{papanicolaou_varadhan}, 
\cite{balzano}, \cite{balzano_notarantonio}.
This condition is equivalent to the statement that
 $\sup_n \kappa_{n} \left\lVert \bar{\psi}^{x, b, \sigma, U}_{n, \alpha} \right\rVert_{\text{tv}} < \infty$
for some point $x \in U$.  By~\eqref{domination_by_classical_kernel_eqn}, 
\eqref{size_random_sets_bound_eqn} implies that $\chi_{n}$
exists such that~\eqref{uniform_capacity_bound_eqn} holds.
Equations~\eqref{point_set_capacity_sigma_V_B_eqn} 
and~\eqref{domination_by_classical_kernel_eqn}
and Lemma~\ref{approximate_mu_in_tv_preserve_stable_convergence_lemma} show that 
the conclusion of Lemma~\ref{general_basic_example_lemma} also
holds if  $q_n$ is defined by $q_n(x) = \left\lVert \bar{\psi}^{ x, b, \sigma, U}_{n, \alpha}\right\rVert_{\text{tv}}$.

The proof of Lemma~\ref{general_basic_example_lemma}
uses that fact that when $U = {\mathbb R}^{d}$ and $\sigma = 1$, 
a nice potential kernel exists (\cite{stroock}).
Lemmas~\ref{restriction_of_form_lemma}
and~\ref{continuous_time_change_lemma}
are again used to reduce the proof to that setting.

When the state space is a manifold rather 
than a subset of ${\mathbb R}^{d}$,  convergence should still be determined by 
local behavior.  Thus
Lemma~\ref{localization_lemma} allows one to extend
  Lemma~\ref{general_basic_example_lemma} 
to the case of a diffusion on a $d$-dimensional 
Riemannian manifold, $d \geq 2$, 
whose topology has a countable base.
This gives a more general form of Theorem~4.2 of \cite{balzano_notarantonio},
which deals with  the Laplace-Beltrami operator on 
a compact Riemannian manifold with boundary, when 
the sets $\Lambda\left( n\right)$ are unions of iid random geodesic balls.

 \section{Dirichlet form properties}
 \label{Dirichlet_form_properties_section}
Here we summarize facts which are used, with references or proofs.
$X$ is assumed to be as in Section~\ref{assumptions_X_subsection}.
Smooth measures were defined in that section.
The proof of Theorem~2.3.15 in \cite{chen_fu} gives:
\begin{lemma}
\label{smooth_gives_finite_energy_on_most_of_space_lemma}
For any finite smooth measure $\mu$ and any $\varepsilon >0$
 there exists $ F \in \mathcal{B}$ such that 
$\mu\left( E - F \right) < \varepsilon$
 and $ \mathbf{1}_{F} \mu $ has finite energy.
\end{lemma}

Let~$\eta$ be a smooth measure which is the Revuz measure for
the positive continuous additive functional $\left( A_t\right)$,  i.e.\  
for any $f \in \mathcal{B}^{+}$, 
$\lim_{t \downarrow 0} \tfrac{1}{t}\mathbf{E}_{m} \left[ 
\int_0^t f\left(X_s\right) 
\, d A_{s} \right] = \int f d \eta$.
The proof of Theorem~4.1.1 in  \cite{chen_fu}  or 
Theorem~4.1.13 in \cite{oshima} shows
that this equation holds if and only if for all $\alpha \in (0, \infty)$
and all $f, h \in \mathcal{B}^{+}$, 
 $\mathbf{E}_{h m} \left[  \int_0^{\infty} e^{ - \alpha t} f \left( X_{t}\right)
 \, d A_t \right]
= \int f  \left( \hat{R}_{\alpha} \,  h\right) \, d \eta$.
It follows that when~$f \eta$
has finite energy, for $\mathcal{E}$-q.e.~$x$ we have
\begin{equation}
\label{background_additive_func_pot_energy_eqn}
 \mathbf{E}_{x} \left[  \int_{0}^{\infty} e^{ - \alpha t} f\left( X_t \right)
 \, d A_t \right]
= G_{\alpha}\,  \left( f \eta \right)(x).
\end{equation}

$\alpha$-excessive functions are defined in III.1.1 of \cite{ma_roeckner}.
By III.1.2(iii), $G_{\alpha}\,  \mu$ is $\alpha$-excessive 
for any measure $\mu$ with finite energy.
Let $u$ be a function 
with $u \geq 0$  and
$e^{ - \alpha t} p_{t}\, u \leq u$, $m$-a.e.,
for all $t >0$.  Suppose also that $u$ has an $\mathcal{E}$-quasi-continuous 
version~$\tilde{u}$.
Then $e^{ - \alpha t} p_{t}\, u \leq \tilde{u}$
 holds $\mathcal{E}$-q.e.\ on~$E$
 by IV.3.3 (iii) of 
\cite{ma_roeckner}, since $p_{t}\, u$
is $\mathcal{E}$-quasi-continuous
by IV.2.9 of \cite{ma_roeckner}. 
By the right continuity 
of $t \mapsto \tilde{u} \circ X_t$ we have for $\mathcal{E}$-q.e.~$x$ that
$\liminf_{t \searrow 0} e^{ - \alpha t} p_{t}\, u (x)
\geq \tilde{u}(x)$.   Similar facts hold
for $\beta R_{ \beta + \alpha} \,  {u}$. Thus 
for $\mathcal{E}$-q.e.~$x$,
\begin{equation}
\label{converge_up_to_excessive_eqn}
e^{ - \alpha t} p_{ t}\, u (x)
\nearrow \tilde{u}(x) \text{ as } t \searrow 0, \ 
\beta R_{\beta + \alpha} \, u (x) 
\nearrow \tilde{u}(x) \text{ as } \beta \nearrow \infty.
\end{equation}
It is also easy to prove the following.
\begin{lemma}
\label{continuous_excessive_gives_submartingale_lemma}
Let $v \geq 0$ be $\mathcal{E}$-quasi-continuous 
and such that $e^{ - \alpha t} p_{t}\, v \leq v$ holds 
$\mathcal{E}$-q.e. for each $t$ (for example, 
let $\nu$ be $\alpha$-excessive and $\mathcal{E}$-quasi-continuous).
Then $t \mapsto e^{ - \alpha t} v \left( X_t \right)$ 
is a supermartingale with respect to $\mathbf{P}_{x}$
for $\mathcal{E}$-q.e.~$x$.
\end{lemma}

For $A \subset E$, by solving III.3.10 of \cite{ma_roeckner} 
we can define
\emph{reduced functions} on~$A$, as follows.  For any function~$f$
on~$E$ which has an $\mathcal{E}$-quasi-continuous version~$\tilde{f}$, 
let
$ \mathcal{L}_{f, A}$ denote the set of all $w \in D\left(\mathcal{E}\right)$
such that $\tilde{w} \geq \tilde{f}$ holds $\mathcal{E}$-quasi-everywhere on~$A$.
Assuming that
$ \mathcal{L}_{f, A} \neq \emptyset$, 
let~$g$ be the unique element in $ \mathcal{L}_{f, A}$
such that $\mathcal{E}_{\alpha}\left(g,w\right)
\geq \mathcal{E}_{\alpha}\left(g,g\right)$ for all $w \in  \mathcal{L}_{f, A}$.
The function~$g$ is an $\alpha$-excessive member of $D\left(\mathcal{E}\right)$,
and
\begin{equation}
\label{background_reduced_bform_agrees_if_q_e_same_on_A_eqn}  
\mathcal{E}_{\alpha}\left(g,v\right) = 0 \text{ for every } v \in D\left(\mathcal{E}\right) \text{ such that }
v = 0 \text{ holds $\mathcal{E}$-q.e.\ on}~A
\end{equation}
We denote any $\mathcal{E}$-quasi-continuous version
of~$g$ by ${f}_{A, \alpha}$, and refer to
${f}_{A, \alpha}$ as the 
$\mathcal{E}, \alpha$-reduced function 
for~$f$ on~$A$.
If $h \wedge {f}_{A, \alpha}$ 
is an $\alpha$-excessive member of $D\left(\mathcal{E}\right)$ 
(in particular if~$h$ itself
is an $\alpha$-excessive member of $D\left(\mathcal{E}\right)$), 
and $h \geq f$ holds $\mathcal{E}$-q.e.\ on~$A$, 
then 
$h \geq {f}_{A, \alpha}$,  $\mathcal{E}$-q.e. 
on~$E$.  If $f \wedge {f}_{A, \alpha}$ 
is an $\alpha$-excessive member of $D\left(\mathcal{E}\right)$ then 
taking $h = f \wedge {f}_{A, \alpha}$ 
shows
${f}_{A, \alpha} = f$ holds 
$\mathcal{E}$-q.e.\ on~$A$.

For any $f \in D\left(\mathcal{E}\right)$, $\alpha \in (0, \infty)$ and any set $A$, 
let ${f}^{\prime}_{A, \alpha}$ be the unique $g \in D\left(\mathcal{E}\right)$
such that $\tilde{f} = \tilde{g}$ holds $\mathcal{E}$-q.e.\ on~$A$ 
and equation~\eqref{background_reduced_bform_agrees_if_q_e_same_on_A_eqn} holds.
An $\mathcal{E}$-quasi-continuous 
version of ${f}^{\prime}_{A, \alpha}$ 
is used whenever pointwise values are needed.   It is easy to check from the definitions
that $\left\lVert{f}^{\prime}_{A, \alpha}\right\rVert_{\mathcal{E}, \alpha}
\leq K_{\alpha} \left\lVert f\right\rVert_{\mathcal{E},  \alpha}$.  
Also, if $f \wedge {f}_{ A,  \alpha}$
is an $\alpha$-excessive member of $D\left(\mathcal{E}\right)$
then ${f}^{\prime}_{A, \alpha} = {f}_{A, \alpha}$, 
and so ${f}^{\prime}_{A, \alpha}$ is an
$\alpha$-excessive element of $D\left(\mathcal{E}\right)$.

Let~$A$ be a closed subset of $E$ such that the cemetery point
$\partial$ is not in the closure of both $A$ and $E - A$.
Let~$u\in D\left(\mathcal{E}\right)$.
Then 
\begin{equation}
\label{background_harmonic_via_entrance_expectation_eqn}  
{u}^{\prime}_{A, \alpha}(z) = 
\mathbf{E}_{z} \left[ 
e^{ - \alpha D_{A}} \tilde{u}\left( X_{D_{A}} \right) \right]
= 
\mathbf{E}_{z} \left[ 
e^{ - \alpha  T_{A}} \tilde{u}\left( X_{ T_{A}} \right) \right]
\end{equation}
holds for $\mathcal{E}$-q.e.~$z$.
To prove equation~\eqref{background_harmonic_via_entrance_expectation_eqn}, 
we note that since $X$ is a Hunt process, 
there exist open sets $U_k$ with $U_k \searrow A$
and $D_{U_k} \nearrow  D_{A}$, $\mathbf{P}_{z}$-a.e.
It is enough to prove equation~\eqref{background_harmonic_via_entrance_expectation_eqn}
when $u$ is bounded. 
The first equality
can be obtained by
applying V.1.6 of \cite{ma_roeckner}
to the open sets~$U_k$, and then using a convergence argument.  The second equality
can be derived from the first
since $\tilde{u}$ and ${u}^{\prime}_{A, \alpha}$
are $\mathcal{E}$-quasi-continuous and $\varepsilon + D_{A} \circ \theta_\varepsilon 
\rightarrow  T_{A}$ as $\varepsilon \searrow 0$.

\begin{lemma}
\label{martingale_up_to_closed_support_lemma}
Let~$\mu$ a measure with finite energy.  Then 
 $e^{ - \alpha t} G_{\alpha}\, \mu \left( X_t \right)$
is a supermartingale with respect to $\mathbf{P}_{x}$ 
for $\mathcal{E}$-q.e.~$x$.
 Let~$A$ be a closed set
with
$\mu\left(A^c\right) = 0$.
Then  ${\left(G_{\alpha}\, \mu\right)}_{A, \alpha}
= G_{\alpha}\, \mu$. If~$A$ is also such that 
$\partial$ is not in the closure of both $A$ and $E - A$,  then $
G_{\alpha}\, \mu (x) 
= \mathbf{E}_{x} \left[ e^{ - \alpha D_{A}} G_{\alpha}\,  \mu \left( X_{D_{A} }\right)  \right]  $
for $\mathcal{E}$-q.e.~$x$, and $e^{ - \alpha t \wedge 
D_{A}} G_{\alpha}\,  \mu \left( X_{ t \wedge D_{A}} \right)$
is a martingale with respect to $\mathbf{P}_{x}$ for $\mathcal{E}$-q.e.~$x$.
\end{lemma}

\medskip \par \noindent
        \mbox{\bf Proof}\ \   
By Lemma~\ref{continuous_excessive_gives_submartingale_lemma}, 
$e^{ - \alpha t} G_{\alpha}\, \mu \left( X_t \right)$
is a supermartingale with respect to $\mathbf{P}_{x}$ 
for $\mathcal{E}$-q.e.~$x$.

 Let~$A$ be a closed set
with
$\mu\left(A^c\right) = 0$.
Using equation~\eqref{background_adjoint_pot_on_measures_eqn}
and equation~\eqref{background_reduced_bform_agrees_if_q_e_same_on_A_eqn}
one has  $\mathcal{E}_{\alpha}\left(G_{\alpha}\, \mu
- {\left(G_{\alpha}\, \mu\right)}_{A, \alpha},G_{\alpha}\, \mu
- {\left(G_{\alpha}\, \mu\right)}_{A, \alpha}\right) = 0$, 
and hence ${\left(G_{\alpha}\, \mu\right)}_{A, \alpha}
= G_{\alpha}\, \mu$.  
Assume~$A$ is also such that 
$\partial$ is not in the closure of both $A$ and $E - A$.
By equation~\eqref{background_harmonic_via_entrance_expectation_eqn}
with ${u}^{\prime}_{A, \alpha} = u = G_{ \alpha}\,  \mu$, 
$G_{\alpha}\, \mu (x)
= \mathbf{E}_{x} \left[ e^{ - \alpha D_{A}} G_{\alpha}\,  \mu \left( X_{D_{A} }\right)  \right]  $
for $\mathcal{E}$-q.e.~$x$.
Since $e^{ - \alpha t} G_{\alpha}\, \mu \left( X_t \right)$
is a supermartingale with respect to $\mathbf{P}_{x}$ for $\mathcal{E}$-q.e.~$x$, it follows that 
that $e^{ - \alpha t \wedge D_{A}} G_{\alpha}\,  \mu \left( X_{ t \wedge D_{A}} \right)$
is a martingale with respect to $\mathbf{P}_{x}$ for $\mathcal{E}$-q.e.~$x$.
\hfill \mbox{\raggedright \rule{0.1in}{0.1in}}

For $\alpha \in (0, \infty)$ and any $\alpha$-excessive $u \in D\left(\mathcal{E}\right)$, 
by VI.2.1 of \cite{ma_roeckner}
there exists a  measure~$\mu$
with finite energy such that $u = G_{\alpha}\, \mu$.
If $u \leq G_{\alpha}\, \nu$, 
then $\mu\left( E\right) \leq \nu\left( E\right)$ by Lemma~\ref{background_bigger_pot_bigger_mass_lemma}
below.

Let $\mathsf{C}$ be the collection of all closed sets
$A$ with $ \mathcal{L}_{1, A} \neq \emptyset$.
Let $A \in \mathsf{C}$.
Since ${1}_{A, \alpha}$ is an $\alpha$-excessive member of $D\left( \mathcal{E}\right)$, 
${1}_{A, \alpha} \wedge 1$ is also an $\alpha$-excessive member of $ \mathcal{L}_{1, A}$, 
and so ${1}_{A, \alpha} = 1$ holds $\mathcal{E}$-q.e.\ on~$A$.
The
unique measure $\gamma$ with finite energy such that
${1}_{A, \alpha} = G_{\alpha}\, \gamma$
will be referred to as the $\alpha$-equilibrium measure for the set~$A$,
and ${1}_{A, \alpha}$ will be called the $\alpha$-equilibrium potential
for~$A$.  Then $\mathcal{E}_{\alpha}\left({1}_{A, \alpha},{1}_{A, \alpha}\right)
= \int \left( G_{\alpha}\,  \gamma \right) \, d \gamma
= \int 1 \, d \gamma = \gamma \left( E\right)$, 
$\gamma$, so $\gamma$ is a finite measure.
Because $A$ is closed, $\gamma\left(A^c \right) = 0$, 
and $\gamma$ is the unique measure
such that $G_{\alpha}\, \gamma = 1$ holds $\mathcal{E}$-q.e.\ on~$A$ and 
$\gamma \left( A^c \right) = 0$.   Define the $\alpha$-capacity 
of $A$, denoted by $\mathsf{Cap}_{\alpha}\left(A\right)$,  
to be $\gamma( E)$. 

For symmetric $\mathcal{E}$,  $w \in  \mathcal{L}_{1, A}$ implies
$\mathcal{E}_{\alpha}\left(w, w\right)
= \mathcal{E}_{\alpha}\left(w - {1}_{A, \alpha} ,w - {1}_{A, \alpha} \right)
+2 \mathcal{E}_{\alpha}\left({1}_{A, \alpha},w - {1}_{A, \alpha}\right)
+ \mathcal{E}_{\alpha}\left( {1}_{A, \alpha} ,{1}_{A, \alpha} \right)
\geq
 \mathcal{E}_{\alpha}\left( {1}_{A, \alpha},{1}_{A, \alpha} \right)$, 
using the definition of reduction.  Hence in the symmetric case,
\begin{equation}
\label{alpha_capacity_in_as_an_infimum_eqn}
\mathsf{Cap}_{\alpha}\left(K\right) = \mathcal{E}_{\alpha}\left({1}_{A, \alpha},{1}_{A, \alpha}\right) = 
\inf_{w \in  \mathcal{L}_{1, K }} \mathcal{E}_{\alpha}\left(w, w\right).
\end{equation}
We can prove that $\alpha$-capacity is monotone, in the sense that 
if $A_1, A_2 \in \mathsf{C}$ 
 with $A_1 \subset A_2$, 
then $\mathsf{Cap}_{\alpha}\left( A_1\right) \leq \mathsf{Cap}_{\alpha}\left(A_2\right)$.
Hence it is convenient to 
to extend the definition of capacity. If $W$ is a closed set which is a
countable union of sets in $\mathsf{C}$,   define
$\mathsf{Cap}_{\alpha}\left(W\right) = \sup \left\{ \mathsf{Cap}_{\alpha}\left(B\right):\ B \in \mathsf{C}, \ B \subset W\right\}$.
Capacities for the nonsymmetric case are defined differently in III.2.8 of \cite{ma_roeckner}, but have
similar properties.
The proof that $\alpha$-capacity is monotone follows
easily from~\eqref{alpha_capacity_in_as_an_infimum_eqn}
in the symmetric case, and in general 
by the next lemma, which is known as the \emph{domination principle}, 
together with Lemma~\ref{background_bigger_pot_bigger_mass_lemma}.
\begin{lemma}
\label{background_domination_principle_lemma}
Let $v \geq 0$ be $\mathcal{E}$-quasi-continuous 
and such that $e^{ - \alpha t}p_{t}\, v \leq v$ holds 
$\mathcal{E}$-q.e. for each $t$.
Let~$\mu$
be a finite measure with finite energy
such that 
$G_{\alpha}\, \mu \leq v$
holds $\mu$-a.e.\ on $E$. Then $G_{\alpha}\, \mu
\leq  v$ holds $\mathcal{E}$-q.e.\ on~$E$.
\end{lemma}

\medskip \par \noindent
        \mbox{\bf Proof}\ \   
Let $f \in L^{2}\left( m\right)$, $f \geq 0$.
There exists a nondecreasing sequence of compact subsets $A_n$
of $E$, 
such that $G_{\alpha}\, \mu \leq v$
everywhere on $A_n$,  and 
 such that $\int_{A_n} \left( \hat{G}_{\alpha}\, f \right) \, d \mu
\nearrow \int \left( \hat{G}_{\alpha}\, f \right) \, d \mu$.
Let $\mu_n = \mathbf{1}_{ A_n} \, \mu$, 
so that $\mu_n \left(A_n^c \right) = 0$.
Let $W_t = e^{ - \alpha t \wedge D_{ A_n} } v \left( X_{t \wedge D_{A_n} }\right)$.
By Lemma~\ref{continuous_excessive_gives_submartingale_lemma}, 
$W$ is a supermartingale with respect to $\mathbf{P}_{x}$ for $\mathcal{E}$-q.e.~$x$,
so $\mathbf{E}_{x} \left[  W_0 \right] \geq \mathbf{E}_{x} \left[  e^{ - \alpha D_{ A_n} }W_{D_{A_n} }  \right]
\geq \mathbf{E}_{x} \left[  e^{ - \alpha D_{A_n} } G_{\alpha}\, \mu_n \left( X_{D_{A_n}} \right) \right]$.
By Lemma~\ref{martingale_up_to_closed_support_lemma}, 
$G_{\alpha}\, \mu_n (x)
= \mathbf{E}_{x} \left[ e^{ - \alpha D_{A_n}} G_{\alpha}\,  \mu_n \left( X_{D_{A_n} }\right)  \right]$
for $\mathcal{E}$-q.e.~$x$. Hence 
$v(x) \geq G_{\alpha}\,  \mu_n (x)$
for $\mathcal{E}$-q.e.~$x$.
Also $\int f \left( G_{\alpha}\,  \mu_n \right) \, d m
= \int \left( \hat{G}_{\alpha}\, f \right) \, d \mu_n
= \int_{A_n} \left( \hat{G}_{\alpha}\, f \right) \, d \mu
\nearrow \int \left( \hat{G}_{\alpha}\, f \right) \, d \mu
= \int f \left( G_{\alpha}\,  \mu \right) \, d m$.
Thus $\int f \left( G_{\alpha}\,  \mu \right) \, d m
\leq \int f v \, d m$.
Since this is true for every nonnegative $f \in L^{2}\left( m\right)$, 
$G_{\alpha}\, \mu \leq v$
holds $m$-a.e., and so by IV.3.3 of \cite{ma_roeckner}, 
$G_{\alpha}\, \mu \leq v$
holds $\mathcal{E}$-q.e.
\hfill \mbox{\raggedright \rule{0.1in}{0.1in}}

\begin{lemma}
\label{background_bigger_pot_bigger_mass_lemma}
Let $\mu, \nu$ be finite measures with finite energy, 
and $\alpha \in (0, \infty)$ such that $G_{\alpha}\, \mu
\leq G_{\alpha}\, \nu$ holds $m$-a.e.
Then $\mu\left(E\right) \leq \nu\left( E\right)$.
\end{lemma}
\medskip \par \noindent
        \mbox{\bf Proof}\ \   
By V.1.7 there exists an $\mathcal{E}$-quasi-continuous
$f \in D\left( \mathcal{E}\right)$ with 
$f >0$, $\mathcal{E}$-q.e. on $E$.
Then $\hat{G}_{\alpha}\, f
= \hat{R}_{\alpha} \, f >0$, $\mathcal{E}$-q.e.
on $E$.
Let $u_n = \left( n \hat{G}_{\alpha}\, f \right) \wedge 1$.
Since $u_n$ is $\alpha$-coexcessive, 
$\mathcal{E}_{\alpha}\left(G_{\alpha}\, \mu, u_n\right) 
\leq \mathcal{E}_{\alpha}\left(G_{\alpha}\, \nu,u_n\right)$
by III.1.2(iii) of \cite{ma_roeckner}.
Hence $\int u_n \, d \mu
\leq \int u_n \, d \nu$ and $u_n \nearrow 1$ $\mathcal{E}$-q.e.
\hfill \mbox{\raggedright \rule{0.1in}{0.1in}}

\begin{lemma}
\label{background_vague_gives_E_weak_lemma}
Let $\mathcal{E}$ be regular.
Let $\mu_n$, $n = 1, 2, \ldots$ be finite measures
with finite energy, and let $\mu$ be a finite smooth measure.
Suppose $\sup_n \left\lVert G_{\alpha}\, \mu_n\right\rVert_{\mathcal{E}, \alpha}
< \infty$ and $\mu_n \rightarrow \mu$ vaguely
as a sequence of measures.  Then~$\mu$ has finite energy and 
$G_{\alpha}\, \mu_n  \rightarrow 
G_{\alpha}\, \mu  $, $\mathcal{E}$-weakly.

\end{lemma}

\medskip \par \noindent
        \mbox{\bf Proof}\ \   
Let $v$ be $\mathcal{E}$-quasi-continuous and let $w \in D\left( \mathcal{E}\right)$
such that $0 \leq v \leq \tilde{w}$ holds $\mathcal{E}$-q.e.
Let $\left(F_k\right)$ be an $\mathcal{E}$-nest such that 
$F_k$ is compact and the restriction of~$v$ to~$F_k$
is nonnegative and continuous for each~$k$.  Let $v_k(x) = v(x)$
for $x \in F_k$, $v_k(x) = 0$ otherwise.
Then $v_k$ is the limit of a decreasing sequence of functions in $\mathcal{C}_{0} \left( E\right)$.
Hence
$\limsup_{n \rightarrow \infty}
\int v_k \, d \mu_n \leq \int v_k \, d \mu
 \leq \int v \, d \mu$.
Let $z_k \in D\left(\mathcal{E}\right)$ be such that
$z_k = 0$ holds $\mathcal{E}$-q.e.\ on $F_k^c$
and $\left\lVert w - z_k\right\rVert_{\mathcal{E}, \alpha} \rightarrow 0$.
Then $\left\lVert \left\lvert w - z_k\right\rvert \right\rVert_{\mathcal{E}, \alpha} \rightarrow 0$.
Let $c = \sup_n \left\lVert G_{\alpha}\, \mu_n\right\rVert_{\mathcal{E}, \alpha}$.
$\int \left\lvert\tilde{w} - \tilde{z}_k\right\rvert \, d \mu_n
= \mathcal{E}_{\alpha}\left( G_{ \alpha}\,  \mu_n , \left\lvert w - z_k\right\rvert \right)
\leq K_{\alpha} c \left\lVert \left\lvert w - z_k\right\rvert \right\rVert_{\mathcal{E}, \alpha}
\leq K_{\alpha} c \left\lVert  w - z_k \right\rVert_{\mathcal{E}, \alpha}$.
Hence  $\left\lvert \int v \, d \mu_n - \int v_k \, d \mu_n\right\rvert
= \int_{F_k^c} v \, d \mu_n
\leq \int_{F_k^c} \tilde{w} \, d \mu_n
= \int_{F_k^c} \left(\tilde{w} - \tilde{z}_k\right) \, d \mu_n 
\leq  K_{\alpha} c \left\lVert  w - z_k \right\rVert_{\mathcal{E}, \alpha}$, and so
$\limsup_{n \rightarrow \infty}
\int v \, d \mu_n \leq \int v \, d \mu$.
Suppose that 
$\lim_{n \rightarrow \infty} \int \tilde{w} \, d \mu_n = \int \tilde{w} \, d \mu$.
Applying what has been shown for $v$ to $\tilde{w} - v$, 
 $\limsup_{n \rightarrow \infty}
\int (\tilde{w} - v) \, d \mu_n \leq \int (\tilde{w} - v) \, d \mu$, so
$\lim_{n \rightarrow \infty}
 \int v \, d \mu_n = \int v \, d \mu$.

Now let $v$ be any $\mathcal{E}$-quasi-continuous function with $v \geq 0$.
Let  $f \in \mathcal{C}_{0} \left( E\right)$ with $f \geq 0$. 
Then $f \wedge v $ is nonnegative and $\mathcal{E}$-quasi-continuous.
 Since $\mathcal{E}$ is regular, there exists $w \in D\left( \mathcal{E}\right) \cap \mathcal{C}_{0} \left(  E\right)$
with $f \leq w$.  Since $\lim_{n \rightarrow \infty} \int w \, d \mu_n = \int w \, d \mu$,
$\lim_{n \rightarrow \infty} 
\int (f \wedge v) \, d \mu_n = \int (f \wedge v) \, d\mu$.
Assume that $v \in D\left(\mathcal{E}\right)$.  We have
$
\int (f \wedge v) \, d \mu
\leq \limsup_{n \rightarrow \infty}
\int v \, d \mu_n
\leq K_{\alpha} c
\left\lVert v\right\rVert_{\mathcal{E}, \alpha}$.  It follows that
$\int v \, d \mu \leq K_{\alpha} c \left\lVert v\right\rVert_{\mathcal{E},  \alpha}$ for every $v \in D\left(\mathcal{E}\right)$
with $v \geq 0$.  Hence for every $v \in D\left(\mathcal{E}\right)$, 
$\left\lvert\int \tilde{v} \, d \mu \right\rvert
\leq \int \left\lvert \tilde{v} \right\rvert \, d \mu \leq c \left\lVert \, \left\lvert v\right\rvert \,\right\rVert_{\mathcal{E},  \alpha}
\leq c \left\lVert v\right\rVert_{\mathcal{E},  \alpha}$, so~$\mu$ has finite energy.

Since $\mathcal{E}$ is regular, $D\left(\mathcal{E}\right) \cap \mathcal{C}_{0} \left(  E\right)$
is dense in 
$D\left(\mathcal{E}\right)$ with respect to 
$\left\lVert\cdot\right\rVert_{\mathcal{E}, \alpha}$-norm.  By vague convergence, 
for $v \in D\left(\mathcal{E}\right) \cap \mathcal{C}_{0} \left(  E\right)$ we have 
$\lim_{n \rightarrow \infty}
\int v \, d \mu_n = \int v \, d \mu$,
i.e.
$\mathcal{E}_{\alpha}\left(G_{\alpha}\, \mu_n,v\right) \rightarrow 
\mathcal{E}_{\alpha}\left(G_{\alpha}\, \mu,v\right)$.
Since $\sup_n \left\lVert G_{\alpha}\, \mu_n\right\rVert_{\mathcal{E}, \alpha}
< \infty$, this convergence holds for all $v \in D\left( \mathcal{E}\right)$
by a $3 \varepsilon$ argument.
\hfill \mbox{\raggedright \rule{0.1in}{0.1in}}

\begin{lemma}
\label{background_weak_gives_L_1_lemma}
Let $\alpha >0$, 
and let 
$B = \left\{x: \  \delta_x \, R_{\alpha} << m\right\}$.
Let $u_n, u$ be nonnegative $\mathcal{E}$-quasi-continuous
functions in $D\left(\mathcal{E}\right)$
such that $u_n \rightarrow u$, $\mathcal{E}$-weakly.
Assume for each~$n$ 
that $\beta R_{\alpha + \beta} \,  {u_n} \leq u_n$
holds $\mathcal{E}$-q.e. on $B$ for all $\beta >0$.
Then $\liminf_{n \rightarrow \infty} u_n \geq u$ 
holds $\mathcal{E}$-q.e. on $B$.

If~\eqref{kernel_density_cond_eqn}
holds then $\lim_{n \rightarrow \infty} \int \left\lvert u_n - u\right\rvert  h \, d m =0$
for all $h \in L^{2}\left(  m\right)$. 
\end{lemma}

\medskip \par \noindent
        \mbox{\bf Proof}\ \   
$\lim_{n \rightarrow \infty} \int u_n h \, d m = \mathcal{E}_{\alpha}\left( u_n,\hat{G}_{\alpha}\, h \right)
\rightarrow \mathcal{E}_{\alpha}\left(u, \hat{G}_{\alpha}\, h \right)
= \int u h \, d m$ for every $h \in L^{2}\left(  m\right)$.
If $x \in B$ and $\beta > 0$, 
let $f_{x, \beta}$ be a density for $ \delta_x \, R_{\alpha + \beta}$
with respect to $m$.  
For  each $x \in B$ and all $c, \beta >0$, 
let $h_{x, \beta}^c =  \beta f_{ x, \alpha + \beta} \wedge c$.
Then for $x \in B$, 
$\int u_n h_{x, \beta}^c \, d m
\leq 
\int u_n  \beta f_{ x, \alpha + \beta} \, d m
= \int u_n \beta \, d \left(  \delta_x \, R_{\alpha + \beta}  \right)
= \beta R_{\alpha + \beta} \,  u_n (x)$.
Also $\beta R_{\alpha + \beta} \,  u_n (x)
\leq u_n(x)$, $\mathcal{E}$-q.e.  
 Thus 
for $\mathcal{E}$-q.e. $x \in B$, 
$\liminf_{n \rightarrow \infty}
 u_n(x) \geq \lim_{n \rightarrow \infty} \int u_n h_{x, \beta}^c \, d m 
= 
\int u h_{x, \beta}^c \, d m$
for all $c >0$, $\beta >0$.
Letting $c \rightarrow \infty$, 
$\liminf_{n \rightarrow \infty}
 u_n(x) \geq \beta R_{\alpha + \beta} \,  u (x)
= \mathbf{E}_{x} \left[  \int_0^\infty \beta e^{ - (\alpha + \beta) t} u \left( X_t \right)  \right]$.
Letting $\beta \uparrow \infty$, by Fatou
we have
$\liminf_{n \rightarrow \infty}
u_n \geq  u$,  $\mathcal{E}$-q.e. on $B$. 
Now suppose that \eqref{kernel_density_cond_eqn}
holds.  Then $\liminf_{n \rightarrow \infty}
u_n \geq  u$, $m$-a.e.
Let $h \in L^{2}\left(  m\right)$.
Then 
$\lim_{n \rightarrow \infty} \int h
\left(  u_n - u \right)^{-} \, dm = 0$
by dominated convergence.
Since $\lim_{n \rightarrow \infty} \int h \left( u_n - u \right) \, d  m = 0$, 
$\lim_{n \rightarrow \infty} \int h
\left(  u_n  - u \right)^{+} \, dm = 0$ as well.
\hfill \mbox{\raggedright \rule{0.1in}{0.1in}}




%
%

\end{document}